%% file: main.tex
\documentclass[a4paper, notitlepage]{report}
\usepackage[utf8]{inputenc}
\usepackage{amsmath,leftindex, amsfonts, mathtools, amssymb, mathdots}
\usepackage[style=alphabetic, citestyle=alphabetic, maxalphanames=8, maxbibnames=8, url=false]{biblatex}
\usepackage[left=3cm, right=3cm, top = 3.5cm, bottom = 3.5cm]{geometry}
\usepackage[percent]{overpic}
\usepackage{slashed}
\usepackage{csquotes}
\usepackage[english]{babel}
\usepackage{amsthm}
\usepackage{tikz-cd}
\usepackage{subcaption}
\usepackage{wrapfig}
\usepackage{url}
\usepackage{graphicx}
\usepackage{tabularx}
\usepackage{mathtools}
\usepackage{todonotes}
\usepackage{import}
\usepackage{booktabs}
\usepackage{bm}           
\usepackage[all]{xy} 
\usepackage{sectsty}
\usepackage{titlesec}
\usepackage{quiver}
\usepackage{mathabx}
\usepackage{mathrsfs}
\usepackage{tikzit}
\usepackage{amsfonts}
\usepackage{xcolor}
\usepackage{hyperref}
\usepackage{physics}
\usepackage{tcolorbox}
\usepackage{enumerate}
\usepackage{ stmaryrd }
\usepackage{ bbold }
\usepackage{ bbm }
\usepackage{bmpsize}
\usepackage{titling}
\usepackage{leftindex}
\usepackage{tabularray}
\usetikzlibrary{matrix}
\usetikzlibrary{cd,snakes}
\usepackage{comment}

\newcommand{\B}{\mathcal{B}}
\newcommand{\C}{\mathcal{C}}

\mathchardef\mhyphen="2D 
\newcommand{\Modhp}{\mathrm{Mod}\mhyphen}

\makeatletter
\renewcommand{\@chapapp}{Lecture}
\makeatother

\theoremstyle{definition}
\newtheorem{theorem}{Theorem}[section]
\newtheorem{definition}{Definition}[section]
\newtheorem{prop}{Proposition}[section]

\newtheorem{remark}{Remark}[section]
\newtheorem{lemma}{Lemma}[section]
\newtheorem{cor}{Corollary}[section]
\newtheorem{ex}{Example}[section]
\newtheorem{fct}{Fact}[section]
\newtheorem*{csl*}{Categorified Schur's Lemma}

\usepackage[font={small,it}]{caption}
\usepackage{listings}
\usepackage{scalerel}
\usepackage{calligra}
\DeclareMathAlphabet{\mathcalligra}{T1}{calligra}{m}{n}
\usepackage{float}

\DeclareMathOperator{\Hom}{Hom}

\DeclareMathOperator{\End}{End}
\DeclareMathOperator{\Vect}{Vec}
\DeclareMathOperator{\BrMod}{BrMod}
\DeclareMathOperator{\Mod}{Mod}
\DeclareMathOperator{\id}{id}
\DeclareMathOperator{\sVec}{sVec}
\DeclareMathOperator{\br}{br}
\renewcommand{\ev}{\operatorname{ev}}
\DeclareMathOperator{\coev}{coev}
\DeclareMathOperator{\Ann}{Ann}

\numberwithin{equation}{section}

\newcommand{\Z}{\mathbb{Z}}

\newtheorem{observation}{Observation}[section]

\newtheorem{notation}{Notation}[section]
\newtheorem{example}{Example}[section]

\usepackage{euscript}

\newcommand\CB			{\EuScript{B}}
\newcommand\CC			{\EuScript{C}}

\newcommand\CM			{\EuScript{M}}

\newcommand\CZ			{\EuScript{Z}}

\newcommand\Cb			{\mathbb{C}}

\newcommand\Kb          {\mathbb{K}}

\newcommand\vect			{\mathrm{Vec}}
\DeclareMathOperator{\sHBA}{sHBA}

\newcommand*\circled[1]{\tikz[baseline=(char.base)]{
            \node[shape=circle,draw,inner sep=2pt] (char) {#1};}}
\usetikzlibrary{decorations.pathmorphing}

\newtheorem{thm}[equation]{Theorem}
\newtheorem{fact}[equation]{Fact}

\usepackage{dsfont}

\title{Higher tensor categories and their extensions}
\author{STOAT participants}
\date{September 2024 \vskip40pt }

\newcommand{\boxpic}[3]{
	\begin{tikzpicture}[baseline={([yshift=-.5ex]current bounding box.center)}]
	\node at (0,0) {\begin{overpic}[scale=#1]{#2}#3\end{overpic}};
	\end{tikzpicture}
}
\newcommand{\Aut}{\mathrm{Aut}}
\newcommand{\Rep}{\mathrm{Rep}}

\newcommand{\ie}{i.e.\,}

\AtEveryBibitem{\clearlist{language}}
\AtEveryBibitem{\clearfield{issn}}
\begin{document}

\begin{titlepage}

\vspace*{0.5cm}

\large

\begin{center}
{\bf\Large 
Higher tensor categories and their extensions}

\bigskip

{\Large 
Notes from the Scottish Talbot On Algebra and Topology\\
\medskip
September 2024}

\vspace{2cm}

\begin{tabular}{cc}
    {\bf Talks by:} & {\bf Notes by:} \\
     Theodoros Lagiotis & Sean Sanford\\
     Tessa Kammermeier & Lorenzo Riva\\
     Jackson van Dyke & Jannik Gr\"{o}ne\\
     Jack Romo & Diogo Andrade\\
     Markus Zetto & Matthew Cellot\\
     Adri\`{a} Mar\'{i}n Salvador & Matthias Vancraeynest\\
     Alea Hofstetter & Julia Bierent\\
     Leon Liu & Jonathan Davies\\
     Cameron Krulewski & Michail Tagaris\\
     Nivedita & Chetan Vuppulury\\
     Lorenzo Riva & Tudor Caba\\
     Daniel Teixeira & Adrien DeLazzer Meunier\\
     Matthew Yu & Iordanis Romaidis
\end{tabular}

\vspace{2cm}

{\bf Introduction by}

Theo Johnson-Freyd and David Reutter

\vspace{2cm}

{\bf Edited by}

Jennifer Brown, Benjamin Ha\"{i}oun, David Jordan and Patrick Kinnear

\end{center}

\end{titlepage}

\pagenumbering{roman}

\tableofcontents

\chapter*{Preface}
\addcontentsline{toc}{chapter}{Preface}
\import{}{Preface.tex}

\chapter*{Introduction}
\addcontentsline{toc}{chapter}{Introduction}
\pagenumbering{arabic}
\import{}{Intro.tex}

\chapter{Braided Fusion Categories and Drinfeld Centers}\label{lec1}
\import{talks/1.1}{1.1.tex}

\chapter{Higher Fusion Categories}\label{lec2}
\import{talks/1.2}{1.2.higher-fusion-categories.tex}

\chapter{Extension Theory of Fusion \texorpdfstring{$1$}{1}-Categories}\label{lec3}
\import{talks/1.3}{notes_1_3.tex}

\chapter{Braided Module Categories and Half-Braided Algebras}\label{lec4}

\import{talks/2.1}{main.tex}

\chapter{Drinfeld Centers of Fusion \texorpdfstring{$2$}{2}-Categories}\label{lec5}
\import{talks/2.2}{main.tex}

\chapter{Lagrangian Algebras in Drinfeld Centers of Fusion \texorpdfstring{$2$}{2}-Categories}\label{lec6}
\import{talks/2.3}{main.tex}

\chapter{The S-Matrix of a Braided Fusion \texorpdfstring{$2$}{2}-Category}\label{lec7}
\import{talks/3.1}{content.tex}

\chapter{The Nondegeneracy of \texorpdfstring{\(\mathcal{Z}(\Mod{\mathcal{B}})\)}{Z(Mod B)}}\label{lec8}
\import{talks/3.2}{3.2.tex}

\chapter{A Homotopical Classification of Nondegenerate Braided Fusion \texorpdfstring{$2$}{2}-Categories}\label{lec9}
\import{talks/4.1}{STOAT_notes_talk_4_1.tex}

\chapter{Functoriality of Quantum Dimensions}\label{lec10}
\import{talks/4.2}{main.tex}

\chapter{The Klein Invariant}\label{lec11}
\import{talks/4.3}{main.tex}

\chapter{Finishing the Proof}\label{lec12}
\import{talks/5.1}{main.tex}

\chapter{Towards a Classification of Fusion \texorpdfstring{$2$}{2}-Categories}\label{lec13}
\import{talks/5.2}{5.2.tex}

\printbibliography


\end{document}

%% file: Preface.tex
These notes are the result of a week-long workshop held in September 2024 on the subject of extensions of braided fusion categories.

There are now many topic-focused workshops aimed at early-career researchers, and these meetings play an invaluable role in disseminating some of the latest ascendant mathematical ideas and in bringing together young mathematicians with others in their speciality at a key time in their development. An influential example of this style of meeting is the Talbot workshop, named for its original rural venue, where since 2004 graduate students and postdocs have gathered under the guidance of one or two invited mentors to explore a contemporary research topic. The impetus for the current meeting was a similar workshop held annually for over twenty years 
by the topology group of the Universit\'e de Toulouse, in the mountain village of Matemale (initially, in la Llagonne). 
It was after attending the Matemale workshop in 2022 that the current organizers were inspired to try something similar in Scotland. Thus the idea for the Scottish Talbot On Algebra and Topology (STOAT) was born. 

After securing initial funding from ARTIN -- Algebra and Representation Theory in the North, and from the Simons Collaboration on Global Categorical Symmetries, the next question was the topic. This should be a significant contemporary piece of research which could sustain a week's worth of intensive study, and ideally would have dense mathematical connections so that its study would be enriching to participants beyond simply understanding the proof of the main result. Additionally, any topic would also need one or two mentors who should be mathematicians we could be sure would be skilled guides through the material. We were lucky that our first choice of such a topic-and-mentor pairing, namely the paper ``Minimal nondegenerate extensions" and its authors Theo Johnson-Freyd and David Reutter, was available.

With this decided, and having secured Cairngorm Lodge youth hostel near Loch Morlich as our remote location to sequester ourselves for the week, we were able to advertise the workshop. In the end 31 participants from universities in Scotland, England, Canada, France, Germany, Italy, Portugal and the U.S. joined us in the Cairngorms for the week of 23 - 27 September 2024. Having brought up a blackboard from Edinburgh, we were all set for an immersive week of mathematics.

The format of the week was roughly 3 talks per day, given by participants, to a syllabus which was prepared in advance by the mentors. The idea was to leave ample time for follow-up discussions. Each day started with a short talk from one of the mentors, and in the evenings the mentors also gave some informal supplementary talks, which served as mathematical detours from the main path of the syllabus. The main talks also had volunteer note-takers drawn from the participants, and those notes form the basis of this volume.

We gratefully acknowledge the support of the Simons Foundation via the Simons Collaboration on Global Categorical Symmetries (Award no. 888988), and the Isaac Newton Institute's Network Support Scheme (EPSRC grant EP/V521929/1) via the ARTIN network for financially supporting the workshop. We also thank Cairngorm Lodge for their hospitality. We thank the workshop speakers for their carefully-prepared talks, and the note takers for their diligent transcriptions. Finally, we express our sincere thanks to David Reutter and Theo Johnson-Freyd for their thoughtful guidance before, during and after the workshop.

\begin{flushright}
Jennifer Brown, Benjamin Ha\"{i}oun, David Jordan and Patrick Kinnear
\end{flushright}

%% file: Intro.tex


Fusion and braided fusion $1$-categories, formally introduced to mathematicians in \cite{MR2183279} based on physics going back to \cite{MR1002038}, 
are a valuable tool for studying technical questions in traditional representation theory. They arise as the representation categories of finite groups (in characteristic not dividing the order), semisimple Hopf algebras, semisimple weak and quasi Hopf algebras, and more general types of quantum groups. 
Fusion $1$-categories are interesting objects of study in their own right in ``categorified noncommutative algebra.'' Fusion $1$-categories feel somewhat like semisimple Frobenius algebras, and somewhat like finite-dimensional $C^*$-algebras; they also have a feeling all of their own stemming from their categorical nature. Questions about the structure and classification of fusion $1$-categories turn out to be extremely rich (see e.g.\ \cite{MR3966762} for a survey), and their answers draw on techniques as diverse as algebraic number theory, functional analysis, and homotopy theory. Fusion $1$-categories also play deep roles in the construction of interesting  invariants of manifolds and in the analysis of physical systems.


These lectures are about fusion and braided fusion $2$-categories, as introduced in \cite{DR18}. Like their $1$-categorical cousins, fusion $2$-categories are both interesting objects in their own right and valuable tools for studying lower-categorical representation theory. To illustrate this perspective, and to provide a structure and destination, these lectures follow and expand upon the paper \cite{JFR}, which used fusion $2$-categories to settle a long-open question about fusion $1$-categories: every slightly degenerate braided fusion $1$-category admits a minimal nondegenerate extension.

To introduce this question, let us sojourn a few categorical levels down. Fix a commutative ring $R$ and a (typically noncommutative) $R$-algebra $A$, i.e.\ a noncommutative ring $A$ equipped with a homomorphism $R \to Z(A)$. Write $A^e$ for the $R$-algebra $A \otimes A^{\mathrm{op}}$ and $\operatorname{End}_R(A)$ for the $R$-algebra of all endomorphisms of the underlying $R$-module of $A$; the actions of $A$ on itself by left and right multiplication supply a canonical algebra homomorphism $A^e \to \operatorname{End}_R(A)$. The algebra $A$ is called \emph{Azumaya} (over $R$) when $A$ is finitely generated projective and faithful as an $R$-module and this canonical map is an isomorphism. Equivalently, there exists an algebra $B$ so that $A\otimes B$ is Morita equivalent to $R$, i.e. $A$  represents a class in the Brauer group $\operatorname{Br}(R)$ of $R$.
  Let us (somewhat abusively) call $A$ \emph{finite} when it is finitely generated projective as both an $R$-module and as an $A^e$-module. An Azumaya algebra is necessarily finite, and a finite algebra is Azumaya exactly when the inclusion $R \to Z(A)$ is an isomorphism. Thus, for finite algebras, we could equivalently replace ``Azumaya'' with ``nondegenerate,'' meaning that the commutator on $A$ is as nondegenerate as possible. 
Suppose that $A$ is finite and $R \to Z(A)$ is an injection (as happens, for example, when $R$ is a field and $A \neq 0$). Then $\operatorname{End}_R(A)$ is Azumaya, and so every such $A$ does extend i.e.\ embed into an Azumaya algebra. Can you find a smaller Azumaya algebra $B$ containing $A$? For example, the centralizer of $A$ inside $\operatorname{End}_R(A)$ is a copy of $A^{\mathrm{op}}$; for an Azumaya extension $A \hookrightarrow B$, how small can the centralizer $Z(A \hookrightarrow B) := \{b \in B : [a,b]=0 \, \forall a \in A\}$ be? It is certainly not smaller than $Z(A)$; let us say that $A \hookrightarrow B$ is \emph{minimal} when this bound $Z(A) \leq Z(A \hookrightarrow B)$ is saturated. 
It is an interesting exercise to study the problem of obstructions and choices to finding minimal Azumaya extensions. For example, when $A$ is simple (so that $Z(A)$ is a field), its minimal 
Azumaya extensions, if they exist, are naturally a torsor for the kernel of the map $\operatorname{Br}(R) \to \operatorname{Br}(Z(A))$; in general, there can be obstructions parameterized by certain Galois cohomology groups.

The organizing question of these lectures is to take the question of minimal Azumaya extensions, and make the following modifications:
\begin{align*}
    \text{commutative ring} \quad & \leadsto \quad \text{symmetric monoidal category} \\
    \text{base ring $R$} \quad & \leadsto \quad \mathbf{Vec}_{\mathbb C} \\
    \text{noncommutative $R$-algebra} \quad & \leadsto \quad \text{braided monoidal $\mathbb C$-linear category} \\
    \text{linear endomorphisms $\operatorname{End}_R(-)$} \quad & \leadsto \quad \text{Drinfel'd centre $\mathcal Z$}\\
    \text{centre $Z$} \quad & \leadsto \quad \text{M\"uger centre $\mathcal Z_2$}  \\
    \text{finite} \quad & \leadsto \quad \text{fusion} 
\end{align*}
In other words, the question becomes: suppose that you are given a braided fusion category $\mathcal A$; can you find a nondegenerately braided extension $\mathcal{A} \hookrightarrow \mathcal{B}$ which is \emph{minimal} in the sense of saturating the obvious bound $\mathcal Z_2(\mathcal A) \leq \mathcal Z_2(\mathcal A \hookrightarrow \mathcal B)$?

The answer is sometimes No. For example, it turns out that if $\mathcal{Z}_2(\mathcal A) = \mathbf{Rep}_{\mathbb C}(G)$ for some finite group $G$, then a minimal nondegenerate extension exists if and only if a certain obstruction class in the group cohomology $\mathrm{H}^4_{\mathrm{gp}}(G; \mathbb{C}^\times)$ vanishes, and moreover every class does arise from some $\mathcal A$. Through a detailed study of how finite groups act on braided fusion categories, the general case can be reduced to a sequence of always-realizable group-cohomological obstructions together with one final question: what about when $\mathcal Z_2(\mathcal A) = \mathbf{sVec}$, the symmetric monoidal category of super vector spaces (with the Koszul sign rule)?

It turns out that in this case, the answer is always Yes, but that this is a highly nontrivial statement. The reason it is nontrivial is that the cohomology group replacing $\mathrm{H}^4_{\mathrm{gp}}(G; \mathbb{C}^\times)$ is nontrivial --- it turns out to be a copy of $\mathbb{Z}/2\mathbb{Z}$ --- and that unlike in the $\mathbf{Rep}_{\mathbb C}(G)$ case, not every cohomology class is realized by some $\mathcal A$. This is unlike many obstruction theory problems: rather than showing that there is no room for an obstruction, one must  compute the actual obstruction and show that in all cases it does vanish.
It also turns out that fusion and braided fusion $2$-categories are an immensely valuable tool to organize and explain these statements: they clarify who these obstructions are and where they live, and they provide vital computational control. At the most basic level, this is because if $\mathcal A$ is a braided fusion $1$-category, then its 2-category $\operatorname{Mod}(\mathcal A)$ of module categories is an example of a fusion $2$-category, and many $1$-categorical questions about $\mathcal A$ are controlled by $2$-categorical questions about $\operatorname{Mod}(\mathcal A)$.

We believe that the utility of fusion $2$-categories to understand questions about fusion $1$-categories is not limited to the minimal nondegenerate extensions question, and the purposes of these lectures is to develop the methods of fusion $2$-categories in general. Fortunately for our pedagogy, the resolution of the minimal nondegenerate extensions question involved essentially every tool that has so far been developed in fusion $2$-category theory. Thus we chose that question as an organizational framework for the lectures.

The lectures are structured as follows. The first three lectures are recollectory. Lecture~\ref{lec1}, by Theodoros Lagiotis, reviews the basic definitions of fusion and braided fusion $1$-categories and their centres, in order to anticipate the categorification of these notions to come. A fusion $1$-category is in particular a semisimple category, and Lecture~\ref{lec2}, by Tessa Kammermeier, explains what the $2$-categorical analog of ``semisimple'' is. A vital tool in fusion (higher) category theory is the ability to recognize and build extensions of categories graded by a finite group; these methods are introduced by Jackson van Dyke in Lecture~\ref{lec3}.

Fusion $2$-categorical methods begin in earnest in Jack Romo's Lecture~\ref{lec4}, which is the first of five lectures on the Drinfeld centres of fusion $2$-categories. First, Romo addresses the following question: if $\mathcal{A}$ is a monoidal category, then any algebra object $A \in \mathcal{A}$ provides an example of an $\mathcal{A}$-module category; if $\mathcal{B}$ is a braided monoidal category, what type(s) of algebras in and modules of $\mathcal{B}$ take advantage of the braiding? In Lecture~\ref{lec5}, Markus Zetto defines the Drinfeld centre of a fusion $2$-category; the connection with Romo's lecture is that when the fusion $2$-category in question is $\operatorname{Mod}(\mathcal{B})$, its Drinfeld centre $\mathfrak{Z}(\operatorname{Mod}(\mathcal{B}))$ can be modeled in terms of the ``half-braided'' algebras introduced by Romo. Drinfeld centres of fusion higher categories are interesting if you are already inclined to be interested in fusion higher categories. Their application to our question is explained by Adri\`a Mar\'in Salvador in Lecture~\ref{lec6}: a braided fusion $1$-category $\mathcal B$ admits a minimal nondegenerate extension if and only if there is an equivalence of braided fusion $2$-categories $\mathfrak{Z}(\operatorname{Mod}(\mathcal{B})) \cong \mathfrak{Z}(\operatorname{Mod}(\mathcal{Z}_2(\mathcal{B})))$. So the question becomes: recognize when a braided fusion $2$-category is $\mathfrak{Z}(\operatorname{Mod}(\mathcal{Z}_2(\mathcal{B})))$. Towards this end, Alea Hofstetter introduces in Lecture~\ref{lec7} a certain (rectangular) matrix determined by any braided fusion $2$-category $\mathfrak{C}$, which is called the (framed) S-matrix and which encodes commutator-type information of the braided fusion $2$-category; the rows of this matrix are indexed by ``connected components'' of $\mathfrak{C}$, and the columns are indexed by simple objects in a symmetric monoidal $1$-category called $\Omega\mathfrak{C}$ which in the case $\mathfrak{C} = \mathfrak{Z}(\operatorname{Mod}(\mathcal{B}))$ recovers $\Omega\mathfrak{C} = \mathcal{Z}_2(\mathcal{B})$.
 Then Lecture~\ref{lec8}, by Leon Liu, proves that the S-matrix of $\mathfrak{Z}(\operatorname{Mod}(\mathcal{B}))$ is invertible (and in particular square); for example, when $\mathcal{Z}_2(\mathcal{B}) = \mathbf{sVec}$, we learn that $\mathfrak{Z}(\operatorname{Mod}(\mathcal{B}))$ has precisely two components.

 Lectures~\ref{lec9}--\ref{lec11} switch from this algebraic flavour of higher category theory, full of matrices and half-braided algebra objects, to a more homotopical approach. First, Cameron Krulewski's Lecture~\ref{lec9} answers the question: what are all the braided fusion $2$-categories $\mathfrak{C}$ with nondegenerate S-matrix and $\Omega\mathfrak{C} = \mathbf{sVec}$? This is done by first showing that any such category necessarily arises from linearizing a homotopy 4-type $X$ whose only nontrivial homotopy groups are $\pi_4 X = \mathbb{C}^\times$ and $\pi_3 X = \pi_2 X = \mathbb{Z}/2\mathbb{Z}$; the axioms ``the S-matrix is nondegenerate'' and ``$\Omega\mathfrak{C} = \mathbf{sVec}$'' are encoded as constraints on the Postnikov extension data of $X$, and the second step is to use spectral sequence methods to classify the spaces satisfying these constraints. It turns out that there are precisely two options for $\mathfrak{C}$, referred to herein as ``$\mathcal S$'' and ``$\mathcal T$''. We will win if we can come up with an indicator that distinguishes them. This indicator is constructed in Lecture~\ref{lec10} by Nivedita and Lecture~\ref{lec11} by Lorenzo Riva. Specifically, Nivedita and Riva interpret a Klein bottle as ``the trace of $180^\circ$ rotation of a figure-8'' and use this to construct a so-called \emph{Klein invariant} of self-dual objects; a mix of algebraic topology and category theory show that this invariant is well-defined and provide methods for its computation. The reason this is useful is that it turns out that all objects in $\mathcal S$ and $\mathcal T$ are self-dual, and that every object in the nontrivial component of $\mathcal S$ has nonnegative Klein invariant, whereas every object in the nontrivial component of $\mathcal T$ has nonpositive Klein invariant.

 These ingredients are brought together by Daniel Teixeira in Lecture~\ref{lec12}. Suppose that $\mathcal{B}$ is a braided fusion 1-category with $\mathcal{Z}_2(\mathcal{B}) = \mathbf{sVec}$. On one hand, the previous lectures show that $\mathfrak{Z}(\operatorname{Mod}(\mathcal{B}))$ is equivalent either to $\mathcal{S}$ or to $\mathcal{T}$, and that we win if we can prove that $\mathfrak{Z}(\operatorname{Mod}(\mathcal{B})) = \mathcal{S}$. On the other hand, the previous lectures supply a very explicit model of $\mathfrak{Z}(\operatorname{Mod}(\mathcal{B}))$ in terms of half-braided algebras in $\mathcal{B}$. Using this model, Teixeira constructs an explicit object $L^-$ in the nontrivial component of $\mathfrak{Z}(\operatorname{Mod}(\mathcal{B}))$. The proof is finished by computing the Klein invariant of $L^-$ and showing that it is strictly positive.
 
 To conclude our lectures, Matthew Yu describes the landscape of fusion $2$-categories more broadly. Fusion $2$-categories equivalent to $\operatorname{Mod}(\mathcal{B})$ are those with a single connected component; given a $\mathcal{B}$ together with some extra group theoretic information, Yu explains how to construct fusion $2$-categories with multiple components. It turns out that Yu's construction is the most general: the lectures end with a survey of the main result of \cite{class2fusion}, showing that this construction gives a complete classification of fusion $2$-categories. Using these results, one should be able, in principle, to completely analyze the question of minimal nondegenerate extensions of an arbitrary braided fusion $1$-category. The details have not been done, and we invite the interested reader to get involved!

\begin{flushright}
Theo Johnson-Freyd and David Reutter
\end{flushright}

%% file: talks/1.1/1.1.tex
Talk by Theodoros Lagiotis, notes by Sean Sanford and Theodoros Lagiotis.

\section{Introduction}

This workshop is dedicated to studying braided fusion categories, and in particular `minimal, non-degenerate extensions' thereof. This problem is motivated by the more general and difficult problem of the classification of (braided) fusion categories.

The definitions and techniques we develop in this talk were used in \cite{DMNO}, in an attempt to simplify this classification problem. However, to answer the problem of classification of minimal nondegenerate extensions of braided fusion categories we will have to categorify the notions we discuss in this talk. So, one could say that `this is not even their final form'.


Fix an algebraically closed field $k$ of characteristic 0.

We will quickly recall the definition of a fusion category:

\begin{definition}
    A category $\C$ is called \textbf{fusion}, if it is a finite semisimple tensor category.

\begin{itemize}
    \item The term `tensor category' is a bit overloaded. For us, it means $k$-linear abelian rigid monoidal category, such that the bifunctor $\otimes\colon\C\boxtimes\C\to\C$ is bilinear on morphisms, and with $\text{End}_\C(1)=k$.
    \item Rigid means that every object $x\in\C$ has both left and right duals (Left: $x^\lor,ev_x,coev_x$), satisfying the duality axioms. 
    \item Finite means $\C\cong A\text{-mod}^{fd}$ for a finite dimensional algebra $A$. 
    \item In \cite{EGNO}, semisimplicity is included in the definition of 'tensor'.
\end{itemize}
\end{definition}

\begin{ex}
    If $G$ is a finite group, then $\text{Vec}_G$, the category of finite dimensional $G$-graded vector spaces is fusion.
    The simple objects $\delta_g\in \text{Vec}_G$ can be described as follows: $(\delta_g)_x= k$, if $x=g$ and 0 otherwise. The monoidal structure is given by $\delta_g\otimes\delta_h\cong \delta_{gh}$.
\end{ex}

An interesting additional structure on a tensor category, that appears in many examples, is that of a braiding:

\begin{definition}
     A braiding on a tensor category $\C$ is a choice of a family of natural isomorphisms $\beta_{x,y}\colon x\otimes y\to y\otimes x$ for every pair of objects $x,y\in\C$, satisfying coherence axioms.
\end{definition}

\begin{ex}
    The fusion category $\text{Vec}_G$ cannot be braided if $G$ is not abelian. One can see this, as a braiding on $\text{Vec}_G$ is in particular a choice of isomorphisms $\delta_g\otimes\delta_h\cong \delta_h\otimes\delta_g$, and this cannot be true if $gh\neq hg$.
\end{ex}

In fact, from a fusion category $\C$, we can obtain a braided fusion category by taking its so called 'Drinfeld center' $\mathcal{Z}(\C)$. The construction is reminiscent of that of the centre of a ring.

\begin{definition}
    The Drinfeld centre $\mathcal{Z}(\C)$ of the (in general just monoidal) category $\C$, is defined as follows:
    \begin{itemize}
        \item Its objects are pairs $(z,b)$, where $z\in\C$ and 
        $$b_x\colon x\otimes z\to z\otimes x, \quad x\in \C$$
        is a natural isomorphism that satisfies the hexagon.

        Explicitly, the following diagram has to be commutative for all $x,y\in \C$:

\[\begin{tikzcd}
	& {x\otimes (z \otimes y)} && {(x\otimes z) \otimes y} \\
	{x\otimes (y \otimes z)} &&&& {(z\otimes x) \otimes y} \\
	& {(x\otimes y) \otimes z} && {z\otimes (x \otimes y)}
	\arrow["{\alpha_{x,z,y}^{-1}}", from=1-2, to=1-4]
	\arrow["{b_x \otimes \text{id}_y}", from=1-4, to=2-5]
	\arrow["{\text{id}_x\otimes b_y}", from=2-1, to=1-2]
	\arrow["{\alpha_{x,y,z}^{-1}}"', from=2-1, to=3-2]
	\arrow["{b_{x\otimes y}}"', from=3-2, to=3-4]
	\arrow["{\alpha_{z,x,y}^{-1}}"', from=3-4, to=2-5]
\end{tikzcd}\]

        \item A morphism from $(z,b)$ to $(z^\prime,b^\prime)$ is a morphism $f\in \Hom(z,z^\prime)$ such that $(f\otimes\id_x)\circ b_x=b^\prime_x\circ(\id_x\otimes f)$.
    \end{itemize}
\end{definition}

The monoidal structure on $\mathcal{Z}(\C)$ is given as follows: If $(z,b), (z^\prime,b^\prime)\in\mathcal{Z}(\C)$, then define $$(z,b)\otimes(z^\prime,b^\prime):=(z\otimes z^\prime, \tilde{b}),$$

where $\tilde{b}_z$ is given by the following commutative diagram:

\[\begin{tikzcd}
	{x\otimes (z \otimes z^\prime)} && {(x\otimes z) \otimes z^\prime} && {(z\otimes x) \otimes z^\prime} \\
	{(z\otimes z^\prime) \otimes x} && {z\otimes (z^\prime \otimes x)} && {z\otimes (x \otimes z^\prime)}
	\arrow["{\alpha^{-1}}", from=1-1, to=1-3]
	\arrow["{\tilde{b}_x}"', from=1-1, to=2-1]
	\arrow["{b_{x} \otimes \id_{z^\prime}}", from=1-3, to=1-5]
	\arrow["\alpha", from=1-5, to=2-5]
	\arrow["{\alpha^{-1}}", from=2-3, to=2-1]
	\arrow["{\id_z \otimes b^\prime_{x} }", from=2-5, to=2-3]
\end{tikzcd}\]

The natural isomorphisms $b$ make $\mathcal{Z}(\C)$ braided.

Since $\C$ is assumed to be fusion, we have the following:

\begin{fct}
    If $\C$ is fusion, then $\mathcal{Z}(\C)$ is braided fusion.
\end{fct}

\begin{ex}
 $\mathcal{Z}(\text{Vec}_G)\cong D(G)-\text{mod}$, where $D(G)$ denotes the Drinfeld double of $G$.
\end{ex}

\begin{ex}
    For $\C$ a modular fusion category, it is true that $\mathcal{Z}(\C)\cong \C\boxtimes\C^{rev}$.
\end{ex}

The goal of this talk is to explain an `if and only if' criterion for detecting Drinfeld centres.
I.e, given a braided fusion category $\B$, when is it the Drinfeld centre of a fusion category $\C$.

As we will see in future talks, this criterion is important, as it (or a `higher' analog of it) relates to obstructions to minimal non-degenerate extensions.


\section{Non-degeneracy and Drinfeld centres}

We will now work our way towards the definition of this non-degeneracy condition for braided fusion categories.

To do this, we will now review some other classical notions.

Fix $\B$ a braided fusion category with braiding $\beta$, and $\text{Irr}(\B)$ its set of isomorphism classes of simple objects.

Note that for any object $x\in \B$, we have the following natural isomorphism $\psi$:

\includegraphics[width=11cm]{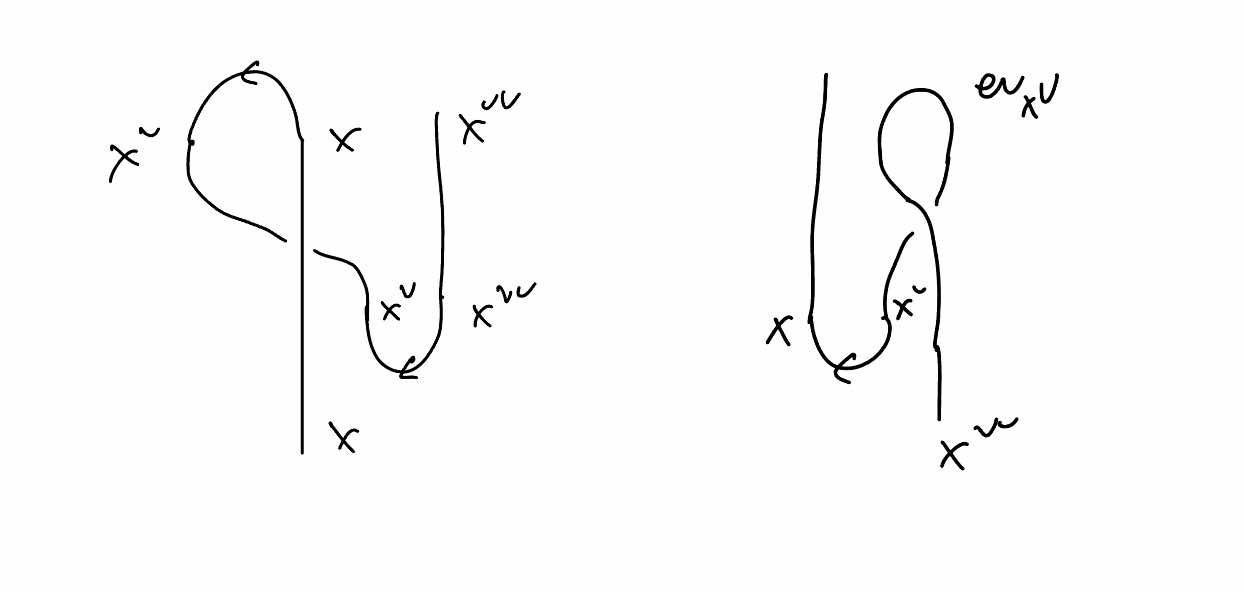}

\begin{fct}
    The natural isomorphism $\psi$ gives rise to a pivotal structure if and only if $\B$ is symmetric.
\end{fct}

\begin{definition}
    \textbf{Traces}: Given $x\in\B$, $f\in \text{End}(x)$ and using $\psi\colon x\to (x^\lor)^\lor$,
    we can define $\text{Tr}^L(f)$ and $\text{Tr}^R(f)$:

\includegraphics[width=13cm]{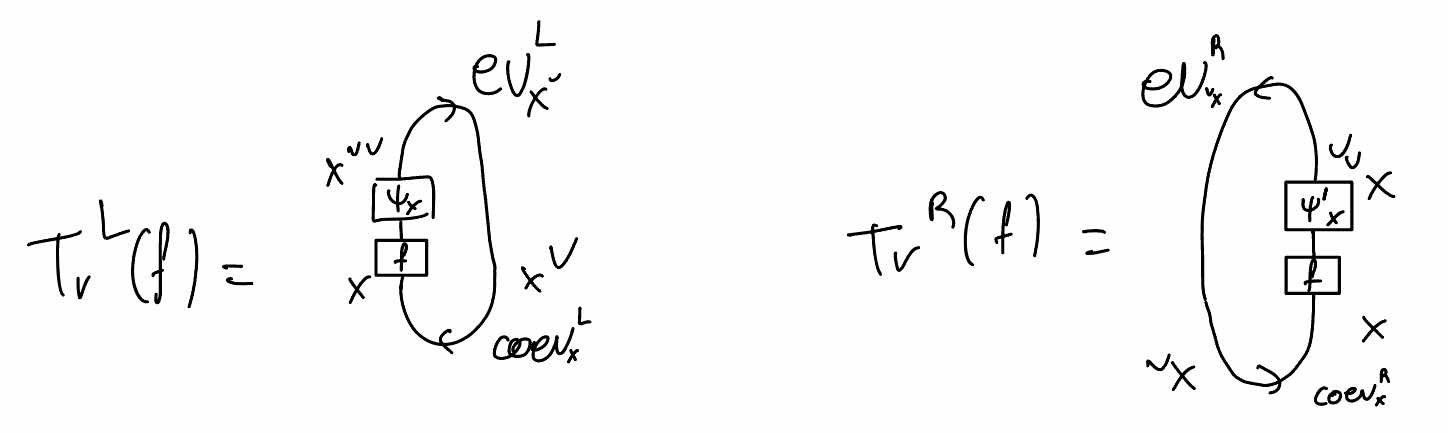}

    In particular, we also define the dimensions $\text{dim}^L(x):=\text{Tr}^L(\id_x)$ and $\text{dim}^R(x):=\text{Tr}^R(\id_x)$.
\end{definition}

\begin{definition}
    Define the matrix $S=(s_{xy})_{x,y\in \text{Irr}(\B)}$ by

$$s_{xy}:=\frac{\text{Tr}^R\otimes\text{Tr}^L(\beta_{y,x}\circ\beta_{x,y})}{\text{dim}^R(x)\text{dim}^L(y)}$$

    
\end{definition}

\begin{definition}
    A braided fusion category $\B$ is said to be non-degenerate if the corresponding matrix $S$ is invertible.
\end{definition}

\begin{fct}
If $\C$ is a fusion category, then $\mathcal{Z}(\C)$ is non-degenerate.
\end{fct}

\section{Detecting Drinfeld centres}

Now we proceed to give the necessary definitions to reach a result that allows us to detect Drinfeld centers. We will be concerned with certain algebra objects in our (braided) fusion category. We follow \cite{DMNO} closely.

Let $\C$ be a fusion category. We can define what it means for $A\in \C$ to be an (associative) algebra with a unit. This just means that we have morphisms $m\colon A\otimes A\to A$ and $u\colon k\to A$ satisfying conditions (associativity, unitality).

We can also make sense of the notion of a (right) module over an algebra $A\in\C$. That is an object $M\in \C$ equipped with a morphism $\alpha\colon M\times A\to M$, that satisfies the axioms of the action.

\begin{definition}
    An algebra $A\in \C$ is called \textbf{separable} if the multiplication map $m\colon A\otimes A\to A$ splits as a morphism of $A$-bimodules.

    In other words, if there exists a map $m^\prime\colon A\to A\otimes A$, such that $m\circ m^\prime=\id_A$.
\end{definition}

For an algebra $A\in \C$, we will denote by $\C_A$, $\leftindex_A{\C}$, $\leftindex_A{\C}_A$ the categories of right, left and bi- A modules respectively.

\begin{fct}
    The following are equivalent:
    \begin{itemize}
        \item $A$ is separable;
        \item The category $\C_A$ is semisimple;
        \item The category $\leftindex_A{\C}$ is semisimple;
        \item The category $\leftindex_A{\C}_A$ is semisimple.
    \end{itemize}
\end{fct}

Now, if we work with a braided fusion category $\B$, we can also make sense of the notion of a commutative algebra in $\B$:

\begin{definition}
   An algebra $A\in \B$ is called \textbf{commutative} if $m\circ\beta_{A,A}=m$.
\end{definition}

\begin{definition}
    A commutative separable algebra $A\in \B$ is called connected if $\Hom(1,A)\cong k$.
\end{definition}

\begin{definition}
    A (right) $A$-module $M\in \C_A$ is called \textbf{local} (or dyslectic) if the following diagram commutes:

\[\begin{tikzcd}
	{M\otimes A} && {M\otimes A} \\
	& M
	\arrow["{\beta_{A,M}\circ\beta_{M,A}}", from=1-1, to=1-3]
	\arrow["\alpha"', from=1-1, to=2-2]
	\arrow["\alpha", from=1-3, to=2-2]
\end{tikzcd}\]

    We will denote the category of local $A$-modules by $\C_A^0\subset\C_A$.
\end{definition}

\begin{definition}
    A \textbf{Lagrangian} algebra $A\in\B$ is a connected commutative separable algebra, such that $\C_A^0\cong \text{Vec}_k$.
\end{definition}

Now we can state the main theorem about detecting Drinfeld centres:

\begin{theorem}
    Let $\B$ be a non-degenerate braided fusion category and $A\in \B$ a Langrangian algebra. Then $\B\cong \mathcal{Z}(\B_A)$.
\end{theorem}

%% file: talks/1.2/1.2.higher-fusion-categories.tex
Talk by Tessa Kammermeier, notes by Lorenzo Riva.

\emph{Notetaker's comment:} I will be using the notation of \cite{JFR}, in particular Remark 2.1. Composition of top morphisms is denoted by $\cdot$ or juxtaposition, and composition of $1$-morphisms in a $2$-category is denoted by $\circ$; horizontal composition of $2$-morphisms (along composable $1$-morphisms) will also be denoted by $\circ$. The symbol $\boxplus$ denotes a direct sum in an arbitrary $2$-category.

\section{$2$-categorical semisimplicity}

\begin{definition}
	A $2$-category $\mathcal{C}$ is \emph{semisimple} if
	\begin{enumerate}
		\item it is \emph{locally semisimple}, i.e. if $\mathcal{C}(X,Y)$ is semisimple for all $X, Y \in \mathcal{C}$;
		\item it has duals for $1$-morphisms, i.e. if for every $1$-morphism $f$ in $\mathcal{C}$ there are adjunctions ${}^\ast f \dashv f \dashv f^\ast$;
		\item and it is additive (\ref{dfn:additive-2-cat}) and idempotent (Karoubi) complete ($\dagger \dagger$).
	\end{enumerate}
\end{definition}

Note that local semisimplicity is relative to a field $\mathbb{k}$, which we assume to be algebraically closed and characteristic zero. The $1$-categorical notions of sums and idempotents have to be refined slightly in order to apply correctly to $2$-categories.

\subsection{Additivity}

Recall that a $1$-category $C$ is \emph{additive} if for all finite sets of objects $\{X_i \in C \mid i \in I\}$ there exists an object $X \in C$ together with inclusions $\iota_i : X_i \to X$ and projections $\rho_i : X \to X_i$ for all $i \in I$ such that
\begin{equation*}
	\rho_i \iota_j = \delta_{i,j} \mathrm{id}_{X_i}, \quad \sum_{i \in I} \iota_i \rho_i = \mathrm{id}_X.
\end{equation*}
In our $2$-categorical setting we only have to ensure that those equations are coherent:

\begin{definition} \label{dfn:additive-2-cat}
	A $\mathbb{k}$-linear $2$-category $\mathcal{C}$ is \emph{additive} if for all finite sets of objects $\{X_i \in C \mid i \in I\}$ there exists an object $X \in C$ (unique up to isomorphism) together with inclusions $\iota_i : X_i \to X$ and projections $\rho_i : X \to X_i$ for all $i \in I$, and coherent isomorphisms
\begin{equation*}
	\rho_i \circ \iota_j \cong \delta_{i,j} \mathrm{id}_{X_i}, \quad \bigoplus_{i \in I} \iota_i \circ \rho_i \cong \mathrm{id}_X,
\end{equation*}
where the latter decomposition occurs in the $\mathbb{k}$-linear $1$-category $\mathcal{C}(X,X)$.
\end{definition}

\subsection{Idempotent completeness}

\begin{definition}
	A \emph{monad} in a $2$-category $\mathcal{C}$ is an algebra object $E \in \mathcal{C}(A,A)$. That is, $E$ is equipped with a unit map $u : \mathrm{id}_A \to E$ and a multiplication map $\mu : E \circ E \to E$ that is coherently associative and unital. A monad $E$ is \emph{separable} if there is a map $\Delta : E \to E \circ E$ satisfying
	\begin{enumerate}
		\item $(\Delta \circ \mathrm{id}_E) \cdot (\mathrm{id}_E \circ \mu) = \Delta \cdot \mu = (\mathrm{id}_E \circ \Delta) \cdot (\mu \circ \mathrm{id}_E)$, i.e. $\Delta$ is a map of $E$-$E$-bimodules,
		\item and $\mu \cdot \Delta = \mathrm{id}_E$, i.e. $\Delta$ is a section of $\mu$.
	\end{enumerate}
	If $\mathcal{C}$ is $\mathbb{k}$-linear then $u, \mu, \Delta$ should be appropriately $\mathbb{k}$-linear as well. 
\end{definition}

Note that the data of a separable monad is the same as that of a split Frobenius algebra object.

\begin{definition}
	An adjunction $f \dashv f^\ast$ is \emph{separable} if the counit $f \circ f^\ast \to \mathrm{id}_Y$ has a section.
\end{definition}

\begin{fct}
	A separable adjunction $f \dashv f^\ast$ induces a separable monad $E = f^\ast \circ f$.
\end{fct}

Recall that, in a $1$-category $C$, an idempotent $e : A \to A$ \emph{splits} if it can be written as $e = i \circ r$, where $r : A \to B$ retracts onto $i : B \to A$ in the sense that $r \circ i = \mathrm{id}_B$. The idempotents in our $2$-categorical setting are the separable monads and the retraction/inclusion pairs are the separable adjunctions:

\begin{definition} \label{dfn:idempotent-complete-2-cat}
	A $2$-category is \emph{idempotent complete} if every separable monad is isomorphic to one induced by a separable adjunction.
\end{definition}

Note that locally idempotent complete $2$-categories admit \emph{idempotent completions}. Given $\mathcal{C}$, its idempotent completion $\mathcal{C}^\triangledown$ can be explicitly constructed by taking the separable monads in $\mathcal{C}$ as objects, the internal bimodules as $1$-morphisms, and the bimodule morphisms as $2$-morphisms. There is a canonical functor $\mathcal{C} \to \mathcal{C}^\triangledown$ which sends $X \in \mathcal{C}$ to the separable monad $\mathrm{id}_X$; when $\mathcal{C}$ is already idempotent complete then this functor is an equivalence.

\section{$2$-categorical simplicity}

\begin{definition}
	An object $A \in \mathcal{C}$ of a semisimple $2$-category $\mathcal{C}$ is \emph{simple} if $\mathrm{id}_A \in \mathcal{C}(A,A)$ is simple (in the $1$-categorical sense).
\end{definition}

Recall that, for semisimple $1$-categories, a non-zero morphism between two simple objects is an isomorphism. This is called Schur's Lemma. We need a weaker version:

\begin{csl*}
	In a semisimple $2$-category, the composition of two non-zero morphisms between simple objects is again non-zero.
\end{csl*}

\begin{proof}
	Let $f : A \to B$, $g : B \to C$ be two non-zero $1$-morphisms where $A$, $B$, and $C$ are simple. By semisimplicity $f$ has a right adjoint $f^\ast$, and the counit $\varepsilon : f \circ f^\ast \to \mathrm{id}_B$ maps into a simple object (since $B$ is simple). But $\varepsilon$ is non-zero, or otherwise the snake equations wouldn't hold, which implies that it has a section $\sigma$. Now assume
	\begin{align*}
		g \circ f = 0 & \implies g \circ f \circ f^\ast = 0 \\
		& \implies \mathrm{id}_g \circ \varepsilon = 0 & \text{(since the domain is $0$)} \\
		& \implies \mathrm{id}_g = (\mathrm{id}_g \circ \varepsilon) \cdot (\mathrm{id}_g \circ \sigma) = 0 & \text{(since $\varepsilon \cdot \sigma = \mathrm{id}_{\mathrm{id}_B}$)} \\
		& \implies g = 0,
	\end{align*}
	and so $g \circ f$ is non-zero.
\end{proof}

This version of Schur's Lemma, together with the existence of adjoints, implies that ``admitting a non-zero morphism'' is an equivalence relation on the set of simple objects. The set of \emph{components} of $\mathcal{C}$,
\begin{equation*}
	\pi_0 \mathcal{C} := \frac{\{\text{simple objects of $\mathcal{C}$}\}}{A \sim B \iff \exists f : A \to B, f \neq 0},
\end{equation*}
is the set of equivalence classes of this relation. Note that two isomorphic simple objects are in the same component but not necessarily viceversa.

\begin{definition}
	We say that a semisimple $2$-category if \emph{finite} if it is locally finite semisimple, i.e. each hom-$1$-category has finitely many simple objects, and $\pi_0 \mathcal{C}$ is finite.
\end{definition}

\begin{fct}
	If $C$ is a multifusion $1$-category then $\Modhp C$ (the $2$-category of finite semisimple $C$-module categories, $C$-module functors, and natural transformations) is a finite semisimple $2$-category. Indeed, any finite semisimple $2$-category $\mathcal{C}$ is equivalent to $\Modhp C$ for some multifusion $1$-category $C$.
\end{fct}

\begin{proof}[Proof sketch for the second statement]
	If $S$ is the set of simple objects of $\mathcal{C}$ we can form $X = \bigboxplus_{Y \in S} Y$. Then $C = \mathcal{C}(X,X)$ is a multifusion $1$-category and the functor $\mathcal{C}(X,-) : \mathcal{C} \to \Modhp C$ is an equivalence.
\end{proof}

Note that $\Modhp C$ is also equivalent to the Morita $2$-category of $C$ whose objects are the separable algebras in $C$, $1$-morphisms are bimodules, and $2$-morphisms are compatible maps in $C$.

\section{(Multi)fusion $2$-categories}

\begin{definition}
	A \emph{multifusion $2$-category} is a finite semisimple $2$-category equipped with a $\mathbb{k}$-linear monoidal structure under which it is rigid, i.e. every object has a left and a right dual. It is said to be \emph{fusion} when its unit object is simple.
\end{definition}

Time for some examples!

\begin{ex}
	Let $G$ be a finite $2$-group. The finite semisimple $2$-category $\mathrm{Fun}(G, 2\mathrm{Vect}) =: 2\mathrm{Vect}[G]$ of $G$-graded $2$-vector spaces has a convolution product turning it into a fusion $2$-category.
\end{ex}

\begin{ex}
	Let $G$ be a finite $2$-group viewed as a $2$-type and let $\omega \in Z^4(G,\mathbb{k}^\times)$ be a topological $4$-cocycle. $\omega$ now corresponds to a map $G \to K(\mathbb{k}^\times, 4)$. We can now form the pullback
\begin{center}
\begin{tikzcd}
\widetilde{G} \arrow[d] \arrow[r] & {PK(\mathbb{k}^\times, 4)} \arrow[d] \\
G \arrow[r, "\omega"]             & {K(\mathbb{k}^\times, 4)}           
\end{tikzcd}
\end{center}	

where $PK(\mathbb{k}^\times, 4)$ is a path space. $\widetilde{G}$ the has the structure of a $3$-group with a trivial $\pi_1$ action on $\pi_3=\mathbb{k}^\times$. We may now think of this $3$-group as a monoidal $2$-category enriched in $\mathbb{k}^\times$-sets. We can base change from $\mathbb{k}^\times$ to $\mathbb{k}$ turning it into a $\mathbb{k}$-linear monoidal $2$-category and then take the semisimple completion of its local Cauchy-completion to form the fusion $2$-category $2\mathrm{Vect}^\omega[G]$ of twisted $G$-graded $2$-vector spaces. 
\end{ex}

\begin{ex}
	Let $C$ be a braided fusion $1$-category. Then $\Modhp C$ is a fusion $2$-category. The monoidal structure is defined with the brading on $C$ which gives rise to a monoidal structure on the one-object $2$-category $BC$. Idempotent completing $BC$ yields the category $BC^\nabla$ which is equivalent to $\Modhp C$.
\end{ex}

\begin{fct}
	Let $\mathcal{C}$ be a fusion $2$-category with one component, meaning that $\lvert \pi_0(\mathcal{C}) \rvert = 1$. Then $\mathcal{C} \simeq \Modhp \mathcal{C}(I,I)$ as fusion $2$-categories.
\end{fct}

%% file: talks/1.3/notes_1_3.tex
Talk by Jackson van Dyke, notes by Jannik Gr\"one.

The goal of this talk is to introduce the homotopy theoretic view on extensions of braided fusion categories. Instead of tackling this directly, we will first look at the ``decategorified'' versions, building up the table below by dimension and observing how the obstruction theory changes along this table.

\begin{center}
\begin{tabular}{lll}
    \toprule
    dim & object & obstruction lives in \\
    \midrule
    1 & $\mathcal{H}$ (Hilbert space) & $H^2(BG;\mathbb{C}^\times)$ \\
    2 & $A$ (algebra) & $H^3(BG;Z(A)^\times)$ \\
    3 & $\mathcal{C}$ (fusion category) & $H^4(BG;\mathbb{C}^\times)$ \\
    4 & $\mathcal{B}$ (BFC) & $H^5(B^2 G;\mathbb{C}^\times)$ \\
    \bottomrule
\end{tabular}
\end{center}

The homotopy-theoretic approach to questions like these was first spelled out in~\cite{ENO10} and most of the results presented here can be found in~\cite{DN21}, in particular in Section 8. For a refresher on the homotopy theory of extensions and Postnikov systems see e.g.~\cite[Section 4.3]{Hat00}.

\section{Thrice decategorified -- Vector Spaces}
Fix a vector space $\mathcal{H}$ (which the physics-inclined among us may view as the Hilbert space of some quantum system) and a group $G$ together with an assignment of operators $g\mapsto T_g \in \mathrm{End}(\mathcal{H})$ which are projectively compatible with the group law, i.e.\ satisfy $T_g\circ T_h = c(g,h)T_{gh}$, where $c(g,h)$ is an element of the underlying field. Checking the compatibility for both bracketings of $ghk$ for $g,h,k\in G$ tells us that $c(g,h)c(gh,k) = c(h,k)c(g,hk)$, which tells us that $c$ represents a class in $H^2(BG;\mathbb{C}^\times)$.

Note that this type of situation is common in physics since the ability to write systems in terms of Hilbert spaces comes at the cost of having an unfixable global phase which has no physical significance.
\begin{example}
    Consider the square-integrable functions on $\mathbb{R}^n$, $L^2(\mathbb{R}^n)$, which encode the wave functions (in e.g.\ position space) of some $n$-dimensional quantum system. This system will naturally come with a unitary action of the symplectic group, which will be projective. We could instead restrict our attention to the pure states and divide the phase out, i.e.\ we consider a map $\mathrm{Sp}_{2n}(\mathbb{R})\to U(L^2(\mathbb{R}^n)) / U(1)$.
\end{example}
\begin{example}
    The best-known example for the other possible approach is probably provided by the appearance of the spin group in quantum mechanics. The cocycle $c$ described above as appearing from a projective action of a group $G$ classifies a central extension $\hat G\to G$. For the rotational group $\mathrm{SO}(n)$ (with $n\ge 3$) the only non-trivial such extension is given by $\mathrm{Spin}(n)$. Instead of considering projective representations of $\mathrm{SO}(n)$ we can study honest representations of $\mathrm{Spin}(n)$.
\end{example}

We will interpret the class $[c]$ as an obstruction to having an honest group action. This situation will be prototypical for the remainder of this talk: We will start with maps (of certain kinds) from an (appropriately interpreted) group into some other algebraic object and ask whether or not this map is has some desired feature, which here means being an honest action, while below we will be concerned with the existence of lifts to larger spaces. We will find that these have cohomological obstructions.

\section{Twice decategorified -- Algebras}
Moving up the categorical dimension we next consider the case of algebra extensions from a point of view that will generalise nicely.
\begin{definition}
    A \emph{categorical group} is a monoidal groupoid with all monoidal inverses.
\end{definition}
\begin{example}
    Given a monoidal category $\mathcal{C}$, its maximal subcategory of invertible objects and invertible morphisms (denoted $\mathcal{C}^\times$) is a categorical group.
\end{example}
Our primary example will be the following: Let $A$ be an algebra (in some ambient category) and denote by $\mathrm{Bim}(A)$ the category of bimodules. Then $\mathrm{Bim}(A)^\times$ (which is the same as the automorphisms of $A$ in its Morita category) is defined to be the \emph{Brauer-Picard categorical group} of $A$ and will be denoted by $\mathbf{BrPic}(A)$.

We will see below that extensions of $A$ are classified by functors from other categorical groups into $\mathbf{BrPic}(A)$. Let us therefore consider what such a functor must look like. First we will need to give a group homomorphism on the (isomorphism classes of) objects, which in the case of $\mathbf{BrPic}(A)$ consist of isomorphism classes of invertible modules (these form the ordinary Brauer-Picard group of $A$, $\mathrm{BrPic}(A)$). This part is clearly unobstructed. We will then need to find a coherent way of dealing with the morphisms, which will lead to cohomological obstruction, similar to the first case. The space of morphisms $\pi_1 \mathbf{BrPic}(A)$ is given by $Z(A)^\times$, i.e.\ the invertible elements in the centre of $A$.

\begin{definition}
    A \emph{$G$-graded associative algebra} is an algebra object $E$ together with a splitting $E=\bigoplus_{g\in G} E_g$ that is compatible with the multiplication (i.e.\ the multiplication map sends $E_g\otimes E_h$ into $E_{gh}$). We also call $E$ a \emph{$G$-extension of $E_e$}, where $E_e$ is the neutral component.

    A $G$-extension is further called \emph{faithful}, if $E_g\neq 0$ for all $g$ and \emph{strong} if the multiplication $E_g\otimes E_h \rightarrow E_{gh}$ is a surjection.
\end{definition}

The problem of finding an extension of an algebra $A$ is the same as mapping into $\mathbf{BrPic}(A)$:

\begin{theorem}
    Strong $G$-extensions of an algebra $A$ are classified by strong monoidal functors $G \rightarrow \mathbf{BrPic}(A)$.
\end{theorem}
Here and later we consider $G$ as a 1- or 2-category by adding only identity morphisms as higher cells.
\begin{proof}[Proof Sketch]
    We first describe how a functor $F:G\to \mathbf{BrPic}(A)$ gives rise to a $G$-extension of $A$. Remember from above that $F$ is on objects given by a group homomorphism $f:G\to \mathrm{BrPic}(A)$. For the functor to be monoidal we must also have isomorphism $\mu(g,h): f(h)\otimes f(g) \xrightarrow{\sim} f(gh)$, which must satisfy a coherence condition:
    \[\begin{tikzcd}
        f(g)\otimes f(h) \otimes f(k) \ar[d] \ar[r] & f(g)\otimes f(hk) \ar[d]\\
        f(gh) \otimes f(k) \ar[r] & f(ghk)
    \end{tikzcd}\]
    We may now set $E:= \bigoplus_{g\in G} f(g)$ on the level of objects. The $\mu$ from above defines the multiplication.

    For the inverse direction we give the functor $G\to\mathrm{Bim}(A)$ as $g\mapsto E_g$ (where $E$ is a $G$-extension of $A$, i.e.\ $E_e=A$). The module $E_g$ will automatically be invertible given that the extension is strong.
\end{proof}

The coherence condition above, when written out, defines a class in the group cohomology $H^3(BG;Z(A)^\times)$. We can interpret this as a topological extension problem. Consider the diagram
\[
    \begin{tikzcd}
        B^2Z(A)^\times \ar[r] & B\mathbf{BrPic}(A) \ar[d] \\
        BG \ar[ur, dashed] \ar[r] & B\pi_0\mathbf{BrPic}(A) \ar[r]& B^3Z(A)^\times
    \end{tikzcd}
\]

The vertical arrow along the centre is the Postnikov tower of the homotopy type associated to $\mathbf{BrPic}(A)$. We begin with a map $BG\to B\pi_0\mathbf{BrPic}(A) = B\mathrm{BrPic}(A)$. This is just the group homomorphism from the proof. The extension to $B\mathbf{BrPic}(A)$ has fiber $B^2Z(A)^\times$ over $B\mathrm{BrPic}(A)$. It is obstructed by an induced map $BG\to B^3Z(A)^\times$, which comes about as the composition of our chosen map $BG\to B\mathrm{BrPic}(A)$ with the $k$-invariant $k:B\mathrm{BrPic}(A)\to B^3Z(A)^\times$ of the fibration $B\mathbf{BrPic}(A)\to B\mathrm{BrPic}(A)$. The homotopy type of $B^3Z(A)^\times$ is by definition that of the Eilenberg-MacLane space $K(Z(A)^\times, 3)$. Therefore the obstruction lives in $[BG,K(Z(A)^\times, 3)] \cong H^3(BG;Z(A)^\times)$.

\section{Once decategorified -- Fusion Categories}
Given $\mathcal{C}$ a fusion category, we define $\mathbf{BrPic}(\mathcal{C})$ to be the 2-categorical group of invertible bimodule categories over $\mathcal{C}$ with 1-morphism invertible bimodule functors and 2-morphisms invertible bimodule natural transformations. We can once again enquire about extensions of $\mathcal{C}$ graded by a group.

\begin{definition}
    A fusion category $\mathcal{D}$ is called \emph{$G$-graded} if it can be written as sum $\mathcal{D} = \bigoplus_{g\in G} \mathcal{D}_g$ such that the grading is compatible with the monoidal product, that is for $x\in\mathcal{D}_g$ and $y\in\mathcal{D}_h$ we have $x\otimes y\in\mathcal{D}_{gh}$. We again call the grading \emph{faithful} if no component is trivial. 
    We call $\mathcal{D}$ a \emph{$G$-graded extension of $\mathcal{D}_e$}.
\end{definition}
\begin{remark}
    We omit the notion of strong here because all faithfully graded fusion categories are also strongly graded, making the notion superfluous for our discussion.
\end{remark}

There is classification of $G$-graded extensions entirely analogous to the case of algebras:
\begin{theorem}
    Faithful $G$-extensions of a fusion category $\mathcal{C}$ are classified by functors $G\to\mathbf{BrPic}(\mathcal{C})$.
\end{theorem}

To understand what these functors look like we can again approach them by building them layer by layer. During the building we encounter obstructions corresponding to coherence conditions which will again have a topological interpretation. Consider the diagram
\[
    \begin{tikzcd}
        B^3\mathbb{C}^\times \ar[r]\ar[d] & B^2\mathcal{Z}(\mathcal{C})^\times \ar[r]\ar[d] & B\mathbf{BrPic}(\mathcal{C}) \ar[d]\\
        * \ar[r] & B^2\pi_1\mathbf{BrPic}(\mathcal{C}) \ar[r]\ar[d] & B\mathbf{Aut}_\mathrm{br}(\mathcal{Z}(\mathcal{C})) \ar[d]\\
        & * \ar[r] & B\pi_0\mathbf{BrPic}(\mathcal{C}) \ar[d]\ar[r,dashed] & B^4\mathbb{C}^\times\\
        & & *
    \end{tikzcd}
\]
The columns of this diagram are Postnikov towers, the rows Whitehead towers. Oftentimes the fibration $B\mathbf{Aut}_\mathrm{br}(\mathcal{Z}(\mathcal{C}))\to B\pi_0\mathbf{BrPic}(\mathcal{C})$ has a section, so extending a map $BG\to B\pi_0\mathbf{BrPic}(\mathcal{C})$ to $BG\to B\pi_{\le1}\mathbf{BrPic}(\mathcal{C})$ will be unobstructed. The lift to the top-level may then still be obstructed. By standard arguments an extension by $B^3\mathbb{C}^\times$ is obstructed by a class in $H^4(BG;\mathbb{C}^\times)$ corresponding to a map $BG\to B\pi_0\mathbf{BrPic}(\mathcal{C}) \to B^4\mathbb{C}^\times$.

Algebraically the obstruction comes from a coherence condition as in the algebra case. In direct analogy we get that the square becomes a cube, describing the various ways for creating a fourfold product.

\section{Braided Fusion Categories}
Let us finally come to the extensions we are most interested in. Let $\mathcal{B}$ be a braided fusion category. We will work by analogy to the decategorified cases. Some of the definitions (such as braided module categories) will be made precise in the second day talks. The Picard 2-categorical group is $\mathbf{Pic}(\mathcal{B}):= (\mathcal{B}-\mathrm{Mod}_\mathcal{C})^\times$. Note that, in difference to the previous section, we take modules rather than bimodules.

Similarly define the braided Picard categorical group as $\mathbf{Pic}_\mathrm{br}(\mathcal{B}):= (\mathcal{B}-\mathrm{BrMod}_\mathcal{C})^\times$

\begin{remark}
    In the case that $\mathcal{B}$ is a Drinfeld centre $\mathcal{Z}(\mathcal{C})$ of some fusion category $\mathcal{C}$, we have $\mathbf{Pic}(\mathcal{B}) \simeq \mathbf{BrPic}(\mathcal{C})$.
\end{remark}

The truncation of $\mathbf{Pic}(\mathcal{B})$ to $\pi_{\le 1}$ maps into the braided automorphisms $\mathbf{Aut}_\mathrm{br}(\mathcal{B})$ by composition of the induction functors for a given module category. The fibre is given by $\pi_{\le1}\mathbf{Pic}_\mathrm{br}(\mathcal{B})$, the truncation of the braided Picard group. It can be understood to ``measure'' the degeneracy of the braiding in that it is related to the categorical Picard group of the symmetric centre of $\mathcal{B}$. As before these categorical groups control graded extensions of $\mathcal{B}$.

\begin{theorem}
    Central $G$-extensions of $\mathcal{B}$ are classified by functors $BG\to B\mathbf{Pic}(\mathcal{B})$.
\end{theorem}

\begin{theorem}
    Braided $G$-extensions of $\mathcal{B}$ (for an abelian group $G$) are classified by functors $B^2G\to B^2\mathbf{Pic}_\mathrm{br}(\mathcal{B})$.
\end{theorem}

\begin{remark}
    If $\mathcal{B}$ is non-degenerate there are no interesting braided extension (i.e.\ all extensions are products).
\end{remark}

The obstruction to a braided extension can be derived to be a class in $H^5(B^2 G;\mathbb{C}^\times)$ in the same fashion as before. It will be our goal to understand this obstruction for non-degenerate braided extensions of slightly degenerate braided fusion categories and show that it indeed vanishes.

%% file: talks/2.1/main.tex
Talk by Jack Romo, notes by Diogo Andrade.

\section{Half-braided algebras and bimodules}

One of our overarching goals is to understand the Drinfeld 2-center of $\Sigma \mathcal{B}$. Instead of unpacking the lengthy definition of a Drinfeld 2-center, we will use our understanding of the 1-categorical Drinfel'd center and half-braidings and use it as motivation for our definitions. 

The main result in this section witnesses a 2-categorical equivalence between the 2-category of separable half-braided algebras in $\mathcal{B}$ and the 2-category of braided $\mathcal{B}$-module categories. Its usefulness comes from the fact the latter 2-category, as we will subsequently prove, provides a model for the Drinfel'd center of $\Sigma\mathcal{B}$.

As a supplement to the main result, we revisit the notion of braided module category using factorization algebras and tangle diagrams. Furthermore, we also revisit, without any attempt at rigour, the notion of surface diagrams for monoidal 2-categories, in order to enhance the reader's intuition for our main definition, that of a half-braided algebras in a braided monoidal category.

\subsection{A whirlwind tour of braided module categories}

From the point-of-view of higher algebra, braided module categories are an immediate and easy object to pin down: they are the ${\mathsf E}_2$-modules of ${\mathsf E}_2$-algebras in $\mathsf{Cat}$. Instead of attempting to give a very rigorous explanation of what this means in terms of operads, we can instead use the result \cite[Theorem 12]{Gin14}, which tells us that these are \textit{stratified locally constant factorization algebras over $\mathbb{R}^{2}_{\ast}$} i.e. the open 2-disk with a single 0-dimensional stratum in the origin.

\subsubsection{Little disks interpretation}
Recall that, broadly speaking, we think of $\mathsf{E}_2$-algebras in $\mathsf{Cat}$ as symmetric monoidal functors from a symmetric monoidal $\infty$-category $\mathrm{Disk}_2$, whose objects are disjoint copies of $\mathbb{R}^2$ and morphisms are embeddings, higher morphisms capturing isotopies of embeddings, and isotopies of isotopies and so on. 

\begin{figure}[h!]
\centering
\includegraphics[scale=0.5]{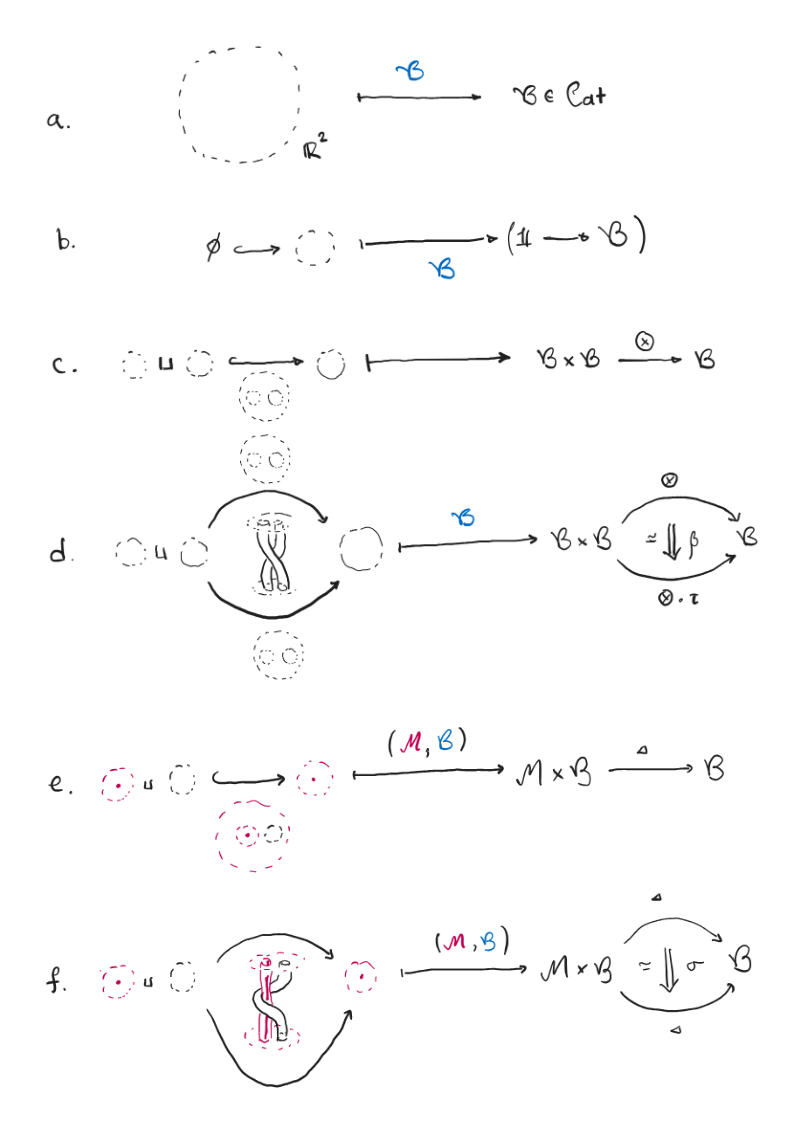}
\caption{A schematic depiction of the data associated to an $\mathsf{E}_2$-algebra in $\mathsf{Cat}$and an $\mathsf{E}_2$-module over it. \textbf{a.} describes the underlying category $\mathcal{B}$ of the $\mathsf{E}_2$-algebra, where \textbf{b.} ensures that the category $\mathcal{B}$ comes equipped with a distinguished object $\mathbf{1}_{\mathcal{B}}$, which will be the unit of the monoidal product $\otimes$ induced by \textbf{c.} The isotopy \textbf{d.} which double-braids the two disks around induces a braiding isomorphism $\beta$. \textbf{e.} captures the $\mathsf{E}_1$-module structure -- a right $\mathcal{B}$-module action in this case. \textbf{f.} is an automorphism of the action functor, which we will henceforth refer to as the $\mathcal{B}$-module braiding. It's an extra structure on top of the right module structure.}
\label{fig:Figure2}
\end{figure}

Braided monoidal categories are thus $\mathsf{E}_2$-algebras in $\mathsf{Cat}$. For the reader who is less familiar with these ideas, note that, an $\mathsf{E}_2$-algebra $\mathfrak{B}$ in $\mathsf{Cat}$ is a symmetric monoidal functor 
\begin{equation}
\mathfrak{B}\colon \mathrm{Disk}_2 \longrightarrow \mathsf{Cat},
\end{equation} 
it is thus described by what it assigns to a small collection of objects and morphisms in $\mathrm{Disk}_2$. We summarize this data in Figure \ref{fig:Figure2}.

Using the aforementioned equivalence between $\mathsf{E}_2$-modules of $\mathsf{E}_2$-algebras and locally constant factorization algebras over $\{(0,0)\} \subset \mathbb{R}^2$, we can also try and establish the categorical information encapsulating the data of an $\mathsf{E}_2$-module of an $\mathsf{E}_2$-algebra in $\mathsf{Cat}$. We again summarize this data in Figure \ref{fig:Figure2}. With this intuition coming from factorization algebras, we are better prepared to read the definition of braided module category.
\begin{definition}
A \textit{braided module category} over the braided monoidal category $(\mathcal{B},\beta)$ is a module category $(\mathcal{M},\gamma)$ over $\mathcal{B}$ together with a natural automorphism $\sigma$ of the functor 
\begin{equation}
\otimes \colon \mathcal{M} \times \mathcal{B} \longrightarrow \mathcal{M}
\end{equation}
satisfying $\eta_{m,\mathbf{1}_{\mathcal{B}}}=\mathrm{id}_m$ for all $m\in\mathcal{M}$ and such that the following diagrams commute for $m\in\mathcal{M}$ and $b,c\in\mathcal{B}$:
\[
\xymatrix{
(m\otimes b)\otimes c \ar[rr]^{\eta_{m\otimes b,c}} \ar[d]_{\gamma_{m,c,b}} & & (m\otimes b)\otimes c  
\\
m\otimes (b\otimes c) \ar[d]_{\mathrm{id}_m \otimes \beta_{b,c}} & & m\otimes (c\otimes b) \ar[u]_{(\gamma_{m,c,b}^{-1})} \\ 
m\otimes (c\otimes b) \ar[d]_{(\gamma)^{-1}_{m,c,b}} & &  m\otimes (b\otimes c) \ar[u]_{\mathrm{id}_m \otimes \beta_{b,c}} \\
(m\otimes c)\otimes b \ar[rr]_{\eta_{m,c}\otimes \mathrm{id}_b} & & (m\otimes b)\otimes c \ar[u]_{\gamma_{m,b,c}} \\ 
}
\]

and 

\[
\xymatrix{
(m\otimes b)\otimes c \ar[rr]^{\eta_{m,b}\otimes \mathrm{id}_c} \ar[dd]_{\gamma_{m,b,c}} & & (m\otimes b)\otimes c \ar[rr]^{\eta_{m\otimes b, c}} & & (m\otimes b)\otimes c \\
& & & & 
\\
m\otimes (b \otimes c)  \ar[rrrr]_{\eta_{m,b\otimes c}} & & & & m\otimes (b \otimes c). \ar[uu]_{\gamma^{-1}_{m,b,c}}
}
\]
\end{definition}

\subsubsection{Graphical calculus perspective and B-tangles}

As is heavily suggested by our previous discussion of locally constant stratified factorization algebras, there is an interpretation of the axioms for braided module categories using tangles rather than little disks, with both interpretations adding valuable intuition to an otherwise ``dry" set of coherence equations.

 In this section we define what we mean by $\mathsf{B}$-\textit{tangles}, and exhibit a theorem, proved in \cite{Bro13}, which says that the free braided $\mathsf{Tang}$-module category generated by a single object, is $\mathsf{BTang}$; where $\mathsf{Tang}$ is the free braided monoidal category generated by a single object.
 
Once again, our point with this subsection is not to give very detailed definition. While you can find a very thorough definition of both $\mathsf{Tang}$ and $\mathsf{BTang}$ in \cite{Bro13} and \cite{Tur16}, we will synthesize them by saying that $\mathsf{Tang}$ has as objects signed points in $\mathbb{R}^2$ and morphisms are oriented tangles in $\mathbb{R}^2\times [0,1]$ which are transverse to the boundary and agree with signature of their source and target points. To define $\mathsf{BTang}$, we use the same definition but replace $\mathbb{R}^2$ with the punctured $\mathbb{R}^2$. In the diagrams below, one should always think there is a red strand implicit, even when the tangles do not braid around it.

The following theorems are generalizations of the Reidemeister theorem and they allows to present braided module structures via generators and relations. These will be helpful in the next subsection.

\begin{theorem}\cite[Thm 5.1]{Bro13}
The ${\sf BTang}$ is generated by the following 1-morphisms which are elementary diagrams (together with every possible orientation)
\begin{figure}[h!]
\includegraphics[scale=0.3]{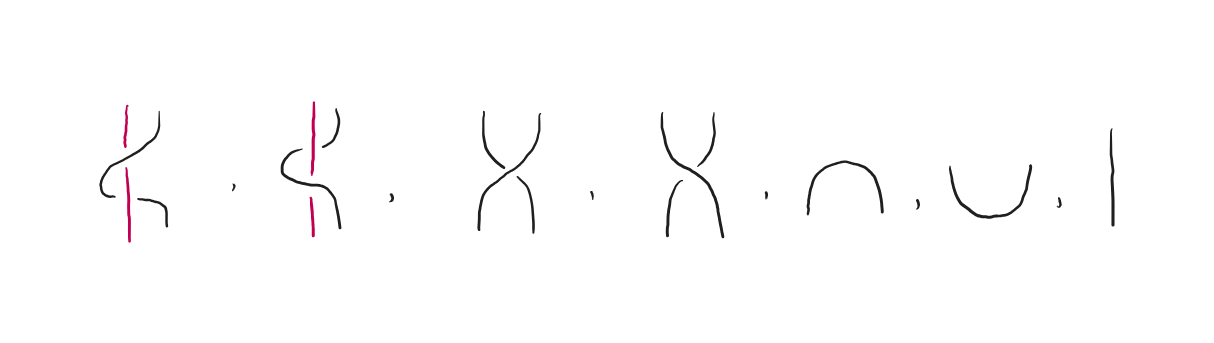}
\centering
\caption{Generators for $\mathsf{BTang}$.}
\label{fig:Figure3}
\end{figure}

together with the following relations:
\begin{enumerate}
\item Planar isotopy.
\item The Turaev moves:\cite[Chapter I.3.2]{Tur16}

\item The reflection relation:
\begin{figure}[h!]
\includegraphics[scale=0.3]{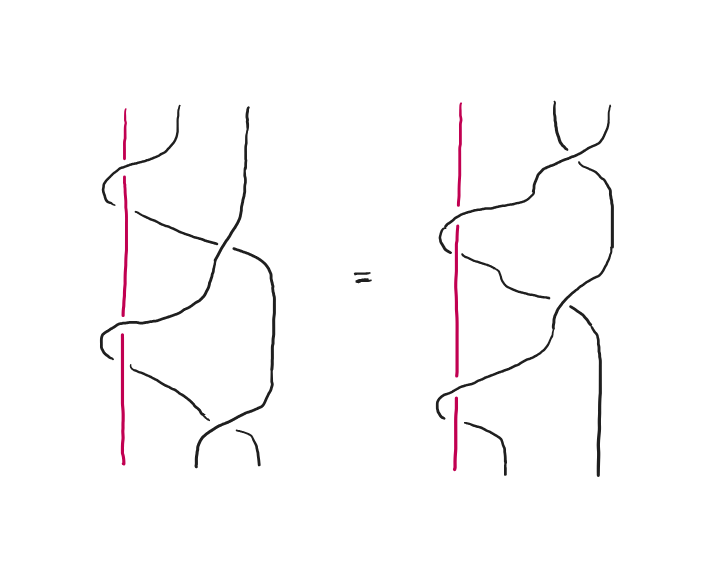}
\centering
\caption{The reflection relation.}
\label{fig:Figure4}
\end{figure}

\item The following relations and their ``upside down version" : (Figure \ref{fig:Figure5})
\begin{figure}[h!]
\includegraphics[scale=0.3]{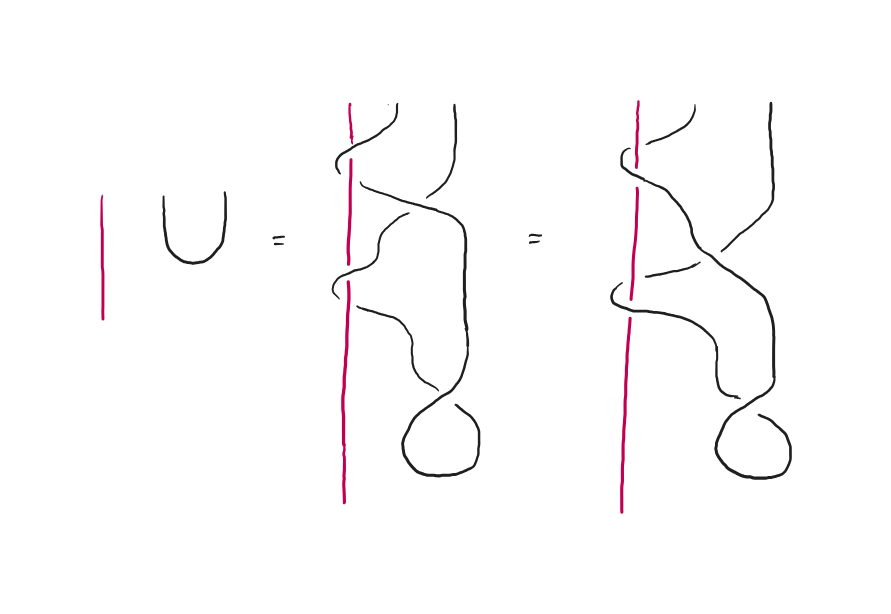}
\centering
\caption{The remaining relations for $\mathsf{BTang}$.}
\label{fig:Figure5}
\end{figure}

\end{enumerate}
\end{theorem}

\begin{theorem}\cite[Thm 5.2]{Bro13}
The pair $(\mathsf{BTang},\tau)$ is a strict braided module category over $\mathsf{Tang}$.
\end{theorem}

\begin{theorem}\cite[Thm 5.3]{Bro13}\label{thm.bm.gen.pres}
Let $(\mathcal{M},\gamma)$ be a strict braided module category $\mathcal{C}$ over $M\in\mathcal{M}$. There exists a unique functor 
\begin{equation}
G\colon {\sf Tang} \longrightarrow \mathcal{M}
\end{equation}
extending $F$ and such that $G(\bullet) = M$ and $G(\tau_{\bullet,+}) = \gamma_{M,V}$.
\end{theorem}

\subsubsection{Surface diagrams and monoidal 2-categories}
In this section we review \textit{en passant} the theory of surface diagrams for monoidal 2-categories. We follow the conventions and notation in \cite[2.1.2]{DR18}.

\begin{figure}[h!]
\includegraphics[scale=0.35]{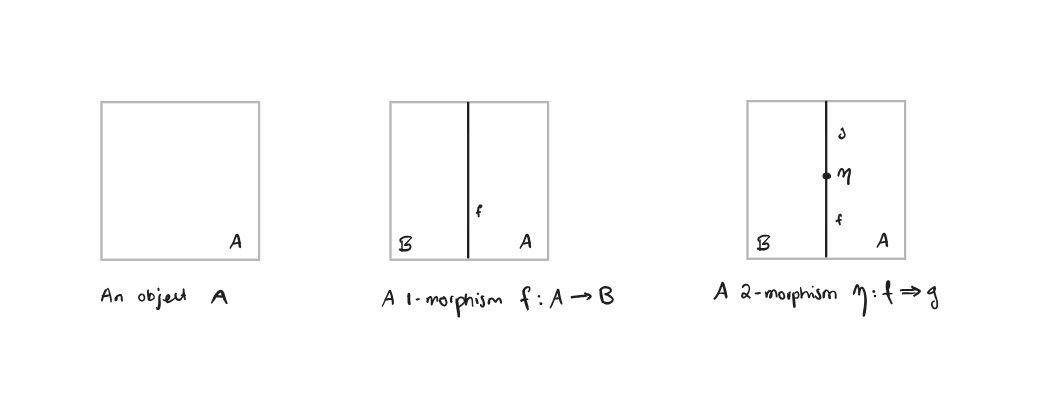}
\centering
\caption{Planar regions encode objects. Defect boundaries separating planar regions depict 1-morphisms, and codimension-2 defects separating two defect boundaries depict 2-morphisms.}
\label{fig:Figure6}
\end{figure}

We depict objects, 1-morphisms, and 2-morphisms in a 2-category using the diagrammatics suggested by Figure \ref{fig:Figure6}.

\begin{figure}[h!]
\includegraphics[scale=0.3]{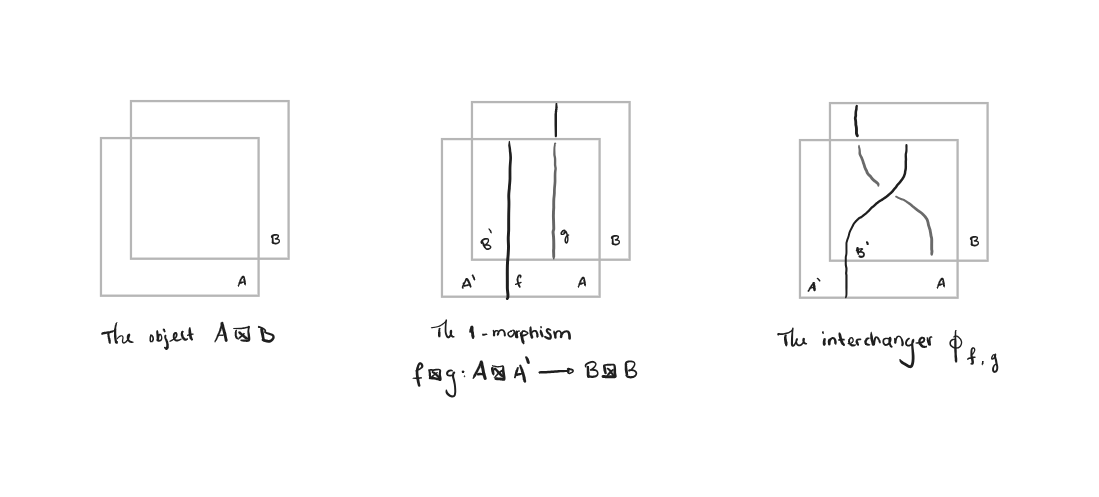}
\centering
\caption{The monoidal structure is encoded by placing sheets next to each other. The interchanger is now structure, rather than an equation, and it is akin to the braiding for monoidal 1-categories.}
\label{fig:Figure7}
\end{figure}

In order to encode monoidal structures, we make use of surface diagrams as embedded sheets in a 3-dimensional copy of Euclidean space. Figure \ref{fig:Figure7} showcases our conventions.

These diagrams are read bottom to top, right to left and front to back, and we can use them to understand otherwise extremely complicated-looking equations in monoidal 2-categories. In particular, we can use it to study the (separable) Morita monoidal 2-category $\Sigma \mathcal{B}.$ In this case, a separable algebra $A\in\mathcal{B}$ is an object in $\Sigma\mathcal{B}$ which we refer to as $[A]$, which we can thus also depict as a labelled sheet (left in Figure \ref{fig:Figure8}):
\begin{figure}[h!]
\includegraphics[scale=0.3]{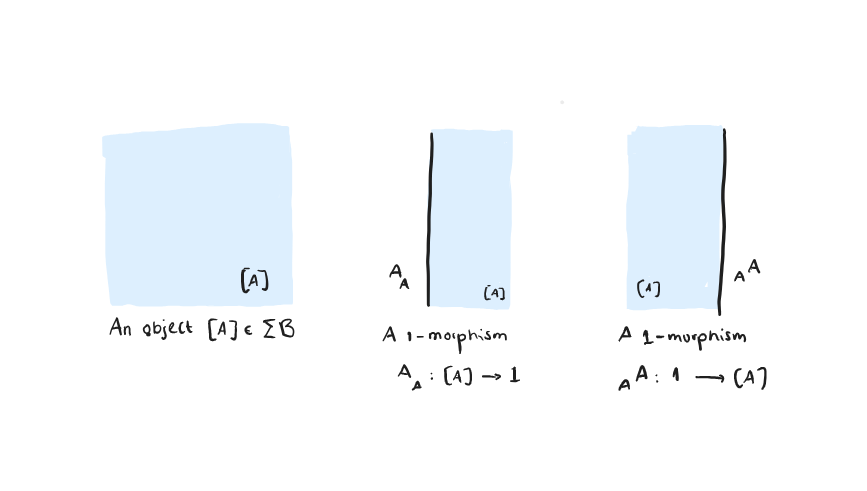}
\centering
\caption{\textbf{Left}: The separable algebra $A$ as an object in $\mathrm{Mod}$-$\mathcal{B}$. \textbf{Middle and right:} The regular right and left module structures on $A$ depicted as 1-morphisms in $\mathrm{Mod}$-$\mathcal{B}$.}
\label{fig:Figure8}
\end{figure}

\begin{figure}[h!]
\includegraphics[scale=0.3]{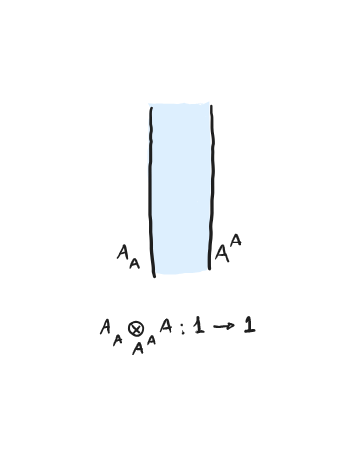}
\centering
\caption{The endomorphism of the unit given by composing the two 1-morphisms.}
\label{fig:Figure9}
\end{figure}

On the other hand, for every separable algebra $A$ there are adjoint 1-morphisms $A_A \colon [A] \rightarrow \mathbf{1}_{\mathcal{B}}$ and $_{A}A\colon \mathbf{1}_{\mathcal{B}}\rightarrow A$ induced by considering the regular left and right $A$-modules, respectively. These are depicted in the middle and right diagrams in Figure \ref{fig:Figure8}. Moreover, we can compose these two to get an interesting endomorphism of the unit. Its diagrammatic representation is captured by the Figure \ref{fig:Figure9}.

From the data contained in these 1-morphisms and their respective diagrams, we can recover the original algebra structure i.e., its multiplication map. Note that the separability condition affords us the following unit and counit maps (Figure \ref{fig:Figure10}):
\begin{figure}[h!]
\includegraphics[scale=0.3]{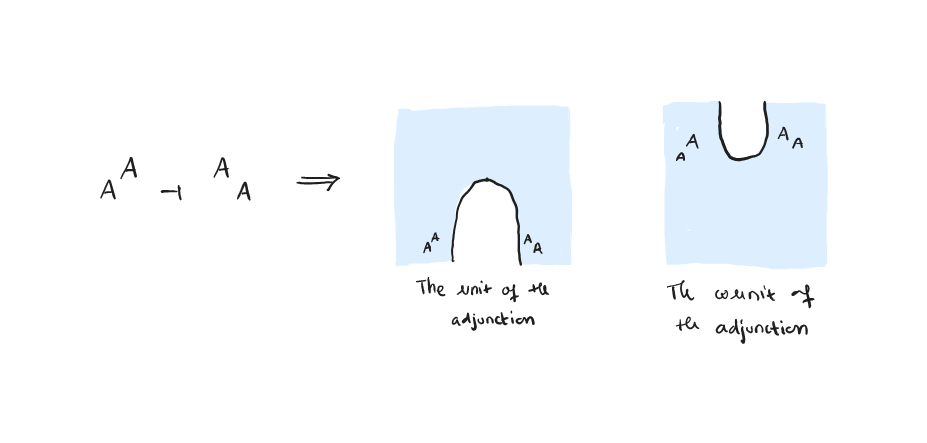}
\centering
\caption{We can use the separability of $A$ to introduce cups and caps into our diagrams.}
\label{fig:Figure10}
\end{figure}
which we thus use to construct a 2-morphism 
\begin{equation}
(A_A \underset{A}{\otimes} \phantom{.}_{A}A) \otimes (A_A \underset{A}{\otimes} \phantom{.}_{A}A) \longrightarrow (A_A \underset{A}{\otimes} \phantom{.}_{A}A)
\end{equation}
whose diagrammatic representation is the following pair of pants.
\begin{figure}[h!]
\includegraphics[scale=0.3]{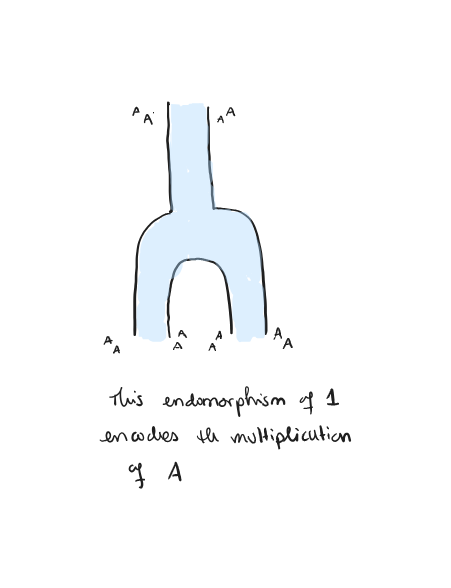}
\centering
\caption{We can use the unit and counit caps to construct a 2-morphism which captures the original algebra structure of $A$.}
\label{fig:Figure11}
\end{figure}

\subsection{Half-braided algebras}

We are now much better placed to begin our discussion of half-braided algebras in a braided monoidal category $\mathcal{B}$. Half-braided algebras in $\mathcal{B}$ are, in particular, associative algebras in $\mathcal{B}$ --- we can, and will, whenever it best suits our intuition, think of them as objects in the monoidal 2-category $\Sigma \mathcal{B}$. The advantage of doing this is access to the calculus of surface diagrams.

Recall that whatever half-braided algebras are, they should be thought of as objects in $\Sigma\mathcal{B}$, together with a 2-half-braiding, which we will refrain from defining, instead using it as a \textit{motif}.

\begin{definition}
A \textit{half-braided algebra} $(A,\gamma)$ in a braided fusion 1-category $\mathcal{B}$ is a unital associative algebra object $A$ whose underlying object is equipped with a half-braiding $(\gamma_b\colon b\otimes A \rightarrow A\otimes b)_{b\in \mathcal{B}}$ fulfilling the axioms depicted in Figure \ref{fig:Figure12}.
\begin{figure}[h!]
\includegraphics[scale=0.3]{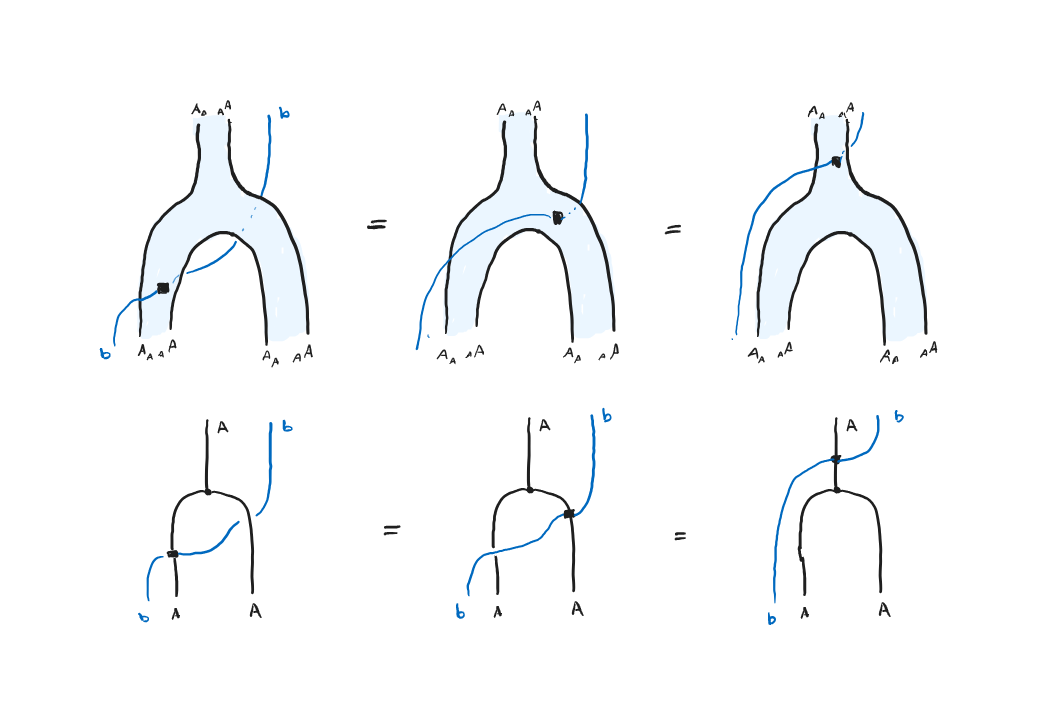}
\centering
\caption{The axioms for a half-braided algebra showcased using both string and surface diagrams. The idea is that the surface diagrammatic is much more intuitive than the string one.}
\label{fig:Figure12}
\end{figure}

\end{definition}
We depict the algebra multiplication of $A$ with a black dot and the half-braiding morphism with a black square.
\begin{definition}
A \textit{homomorphism of half-braided algebras}
\begin{equation}
f \colon (A,\gamma) \longrightarrow (B,\zeta)
\end{equation} is a unital algebra homomorphism $f\colon A \rightarrow B$ which intertwines the half-braiding 
\begin{equation}
\mu_b\cdot (\mathrm{id}_b \otimes f) = (f\otimes \mathrm{id}_b)\cdot \gamma_b.
\end{equation}
\end{definition}

\begin{remark}
We will say that a half-braided algebra is \textit{separable} if it is separable as an associative algebra.
\end{remark}

\begin{observation}
Using the braiding coming from $\mathcal{B}$, we can define a half-braided algebra structure on every algebra object $A\in \mathcal{Z}_2(\mathcal{B})$ in the M\"uger center.
\end{observation}
\begin{observation}
Half-braided algebra structures on the trivial algebra $\mathbf{1}_{\mathcal{B}}$ correspond to monoidal natural automorphisms of $\mathcal{B}$. 
\end{observation}

\begin{lemma}\label{lemma:hba.bimod}
Let $m$ be a left, respectively right, module of an algebra $A\in\mathcal{B}$. Then, any half-braiding 
\begin{equation}
\{\gamma_{b} \colon b \otimes A \rightarrow  A\otimes b\}_{b\in\mathcal{B}}
\end{equation} 
on $A$ satisfying the axioms of a half-braided algebra, equips the object $m\in\mathcal{B}$ with a half-braiding 
\begin{equation}
\{\gamma_{b,m} \colon b\otimes m \rightarrow m \otimes b \}_{b\in \mathcal{B}}
\end{equation}
defined as in Figure \ref{fig:Figure13}, respectively.
\begin{figure}[h!]
\centering
\includegraphics[scale=0.3]{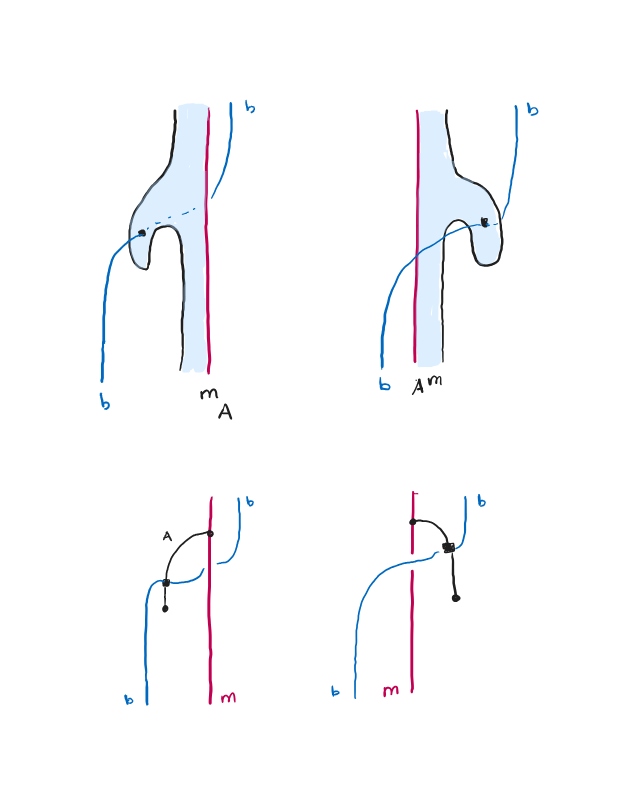}
\centering
\caption{The above diagrams capture the notion of half-braided module induced from the half-braided algebra.}
\label{fig:Figure13}
\end{figure}
For left modules (and analogously for right modules), this extends to a $\mathcal{B}$-module functor 
\begin{equation}
A\text{-}\mathrm{Mod}_{\mathcal{B}} \longrightarrow \mathcal{Z}(\mathcal{B}),
\end{equation}
where $\mathcal{B}$ acts on the right on the Drinfel'd center $\mathcal{Z}(\mathcal{B})$ via the (braided) monoidal functor 
\begin{align*}
\mathcal{B} &\longrightarrow \mathcal{Z}(\mathcal{B}), \\
b &\longmapsto (b,\mathrm{br}_{-,b}).
\end{align*}
\end{lemma}
\begin{proof}
This is immediate from the definition..
\end{proof}

\begin{definition}
A half-braided bimodule $(A,\gamma) \rightarrow (B,\zeta)$ between half-braided algebras $(A,\gamma)$ and $(B,\zeta)$ is a unital $B$-$A$-bimodule $_{B}m_A$ in $\mathcal{B}$ for which the action $B\otimes m \otimes A \rightarrow m$ is compatible with the half-braidings $\gamma$ and $\zeta$ as is depicted in Figure \ref{fig:Figure14}.
\begin{figure}[h!]
\centering
\includegraphics[scale=0.3]{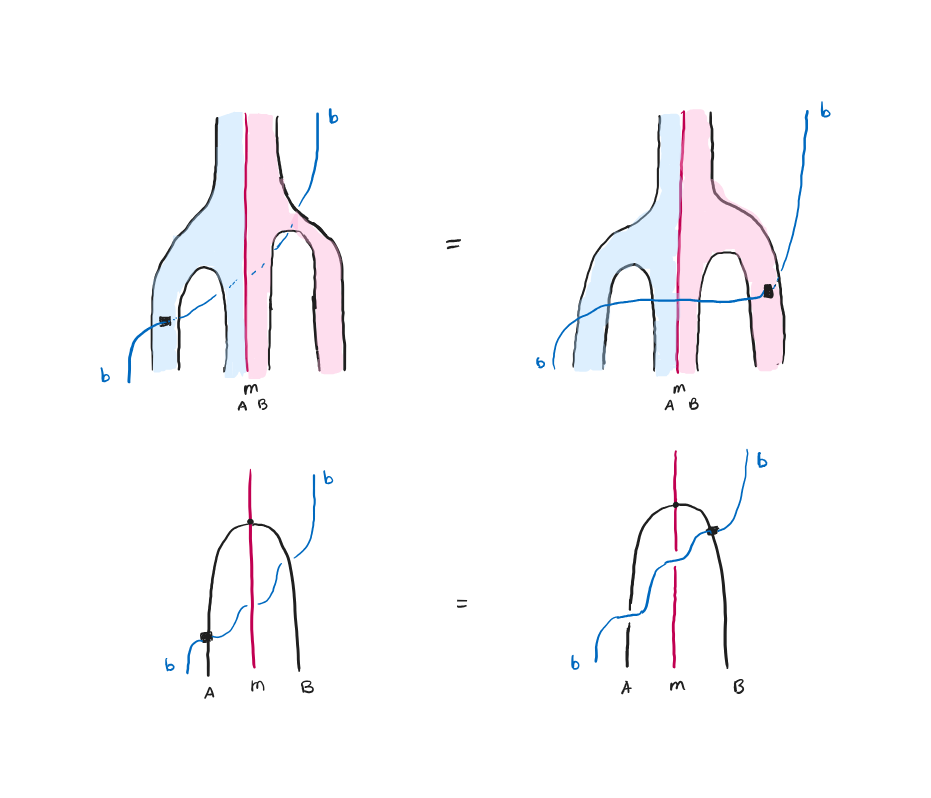}
\caption{The above diagrams codify the equations defining the compatibility for a half-braided bimodule.}
\label{fig:Figure14}
\end{figure}
\end{definition}

\begin{remark}
Being a half-braided bimodule, contrary to being a half-braided algebra, is a property, not an additional structure! In particular, morphisms of half-braided bimodules are just morphisms of bimodules without additional compatibility conditions.
\end{remark}

The next lemma tells us that finding half-braidings for a given associative algebra $A\in \mathcal{B}$ which make it a half-braided algebra, are in bijection with structures turning the right $\mathcal{B}$-module category $A$-$\mathrm{Mod}_{\mathcal{B}}$ into a \underline{braided} $\mathcal{B}$-module category:
\[\{\gamma\colon (A,\gamma) \text{ is half-braided}\} \simeq \{\tau\colon (A\text{-}\mathrm{Mod}_{\mathcal{B}},\tau) \text{ is a braided }\mathcal{B}\text{-}\text{module structure}\}. \]

\begin{lemma}\label{lemma.bijection}
Let $A$ be an algebra in $\mathcal{B}$ with associated right $\mathcal{B}$-module category $A$-$\mathrm{Mod}_{\mathcal{B}}$. For any half-braiding $\gamma$ on $A$ inducing the structure of half-braided algebra, the natural isomorphism 
\begin{equation}
\{ \sigma_{_{A}m,b}:= \gamma_{b, A_{m}} \cdot \beta_{m,b} \colon m\otimes b \rightarrow m\otimes b \}_{b\in\mathcal{B},_{A}m\in A\text{-}\mathrm{Mod}_\mathcal{B}}
\end{equation} 
defines a $\mathcal{B}$-module braiding on $A$-$\mathrm{Mod}_{\mathcal{B}}$. This construction defines a bijection between the set of half-braidings on $A$ which make it into an half-braided algebra, and the set of $\mathcal{B}$-module braidings on $A$-$\mathrm{Mod}_{\mathcal{B}}$.
\end{lemma}

\begin{proof}
Instead of giving a long proof checking that $\sigma$ is a braiding, we make use of Theorem \ref{thm.bm.gen.pres} so that we only have to check the reflection equation (Figure \ref{fig:Figure4}) and the remaining relations (Figure \ref{fig:Figure5}); all other relations, namely isotopy invariance and the Turaev moves, are immediate from the fact that $\mathcal{B}$ is a braided monoidal category. Figure \ref{fig:Figure15} checks the reflection equation:

\begin{figure}[h!]
\includegraphics[scale=0.3]{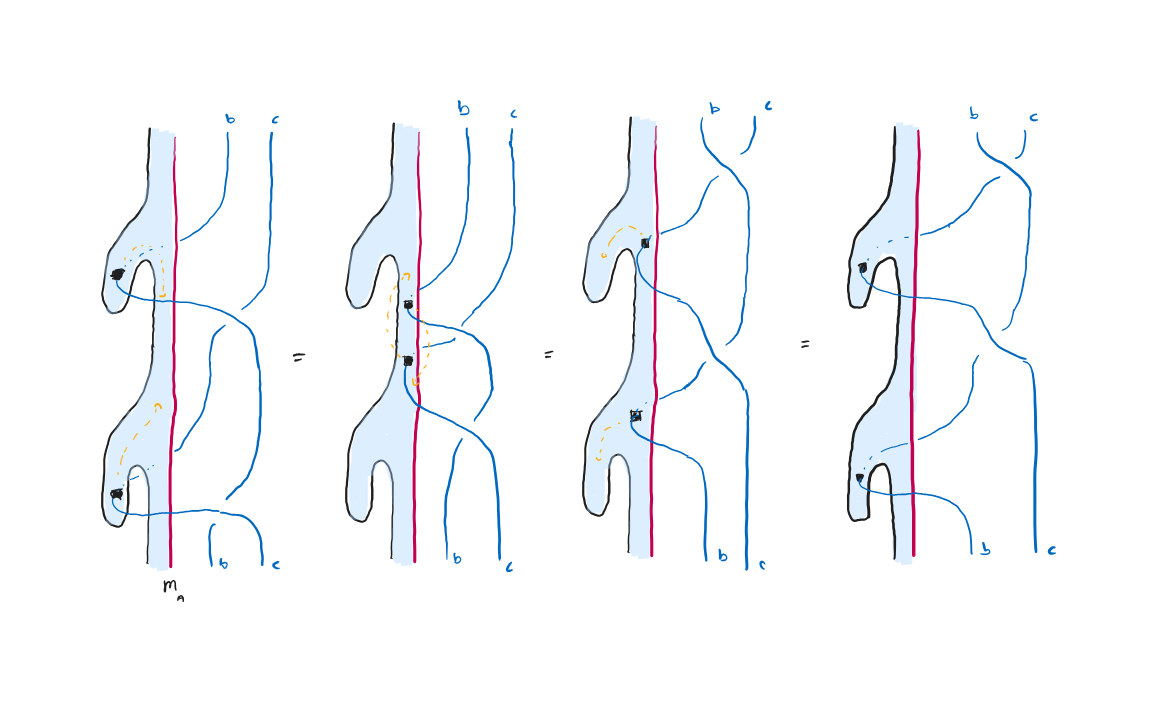}
\centering
\caption{The movie proving that $\sigma$ satisfies the reflection relation.}
\label{fig:Figure15}
\end{figure}
The remaining relations are left for the reader to check.
\end{proof}

\begin{lemma}\label{lemma.bimodules.hb}
Let $_{A}m_{B}$ be a bimodule between half-braided algebras $(A,\gamma)$ and $(B,\zeta)$, then the $\mathcal{B}$-module functor 
\begin{equation}
_{B}m_{A}\underset{A}{\otimes} - \colon A\text{-}\mathrm{Mod}_{\mathcal{B}} \longrightarrow B\text{-}\mathrm{Mod}_{\mathcal{B}}
\end{equation} 
is a braided module functor if and only if the bimodule $_{B}m_A$ is half-braided.
\end{lemma}

\begin{theorem}
The 2-category $\mathrm{sHBA}(\mathcal{B})$ of separable half-braided algebras in $\mathcal{B}$ is equivalent to $\mathrm{BrMod}$-$\mathcal{B}$, the 2-category of braided $\mathcal{B}$-module categories, $\mathcal{B}$-module functors and $\mathcal{B}$-module natural transformations, and hence to the Drinfel'd centre $\mathcal{Z}(\Sigma\mathcal{B})$
\end{theorem}
\begin{proof}
The 2-functor is defined as follows:
\begin{align*}
\mathrm{sHBA}(\mathcal{B}) &\longrightarrow \mathrm{BrMod}\text{-}\mathcal{B} \\
\\
(A,\gamma) &\longmapsto (A\text{-}\mathrm{Mod}_{\mathcal{B}},\sigma) \\
\\
_{B}m_{A} &\longmapsto \big( \phantom{}_{B}m_{A}\underset{A}{\otimes} - \colon A\text{-}\mathrm{Mod}_{\mathcal{B}} \rightarrow B\text{-}\mathrm{Mod}_{\mathcal{B}}\big)\\
\\
(_{B}m_{A} \rightarrow \phantom{}_{B}n_{A}) &\longmapsto \big(f\underset{A}{\otimes}- \colon _{B}m_{A}\underset{A}{\otimes} - \Longrightarrow \phantom{}_{B}n_{A}\underset{A}{\otimes}-\big),
\end{align*}
where $\sigma$ is the induced $\mathcal{B}$-module braiding induced from the half-braiding $\gamma$, which we introduced in Lemma \ref{lemma.bijection}. Finally, to check that this 2-functor witnesses an equivalence of 2-categories we take note of the previously shown equivalence 
\begin{equation}\label{sigmaB.equiv}
(-)\text{-}\mathrm{Mod}_{\mathcal{B}}\colon \Sigma\mathcal{B} \xrightarrow{~ \simeq ~} \mathrm{Mod}\text{-}{\mathcal{B}}.
\end{equation} It can be used to establish fully-faithfulness of our 2-functor. Essential surjectivity can be derived by noting that, given an arbitrary $\mathcal{B}$-module category $(\mathcal{M},\sigma)$, its underlying $\mathcal{B}$-module is equivalent to $A$-$\mathrm{Mod}_{\mathcal{B}}$ for some separable associative algebra $A \in \Sigma\mathcal{B}$ --- this follows from essential surjectivity of (\ref{sigmaB.equiv}). Moreover, using the bijection established in Lemma \ref{lemma.bijection}, we note that the braiding $\sigma$ can be used to construct an half-braided algebra structure on $A$ whose induced braiding on $A$-$\mathrm{Mod}_{\mathcal{B}}$ is $\sigma$. 
\end{proof}

%% file: talks/2.2/main.tex
Talk by Markus Zetto, notes by Matthew Cellot.

\section{Motivation}
In the 1-categorical setting, the operations of taking the Drinfel'd and Müger center fit into the following diagram:
\begin{center}
\begin{tikzcd}
    \mathrm{Cat} \arrow[r, "\mathrm{End}"] & \mathrm{Alg}_{E_1}(\mathrm{Cat}) \arrow[r, "\mathcal{Z}"] & \mathrm{Alg}_{E_2}(\mathrm{Cat}) \arrow[r, "\mathcal{Z}_2"] & \mathrm{Alg}_{E_3}(\mathrm{Cat}) = \mathrm{CAlg(Cat)}
\end{tikzcd}
\end{center}
\paragraph{Remark} We also could have started with the category of pointed categories $\mathrm{Alg}_{E_0} (\mathrm{Cat})$ and point-preserving endomorphisms instead of $\mathrm{Cat}$.
\paragraph{Warning} None of these constructions are functorial, at least not in a straightforward way. For instance, $\mathrm{End}(\mathcal{C}) = \operatorname{Fun}(\mathcal{C}, \mathcal{C})$ has a simultaneously co- and contravariant dependence on $\mathcal{C}$.

The above assignments admit 2-categorical analogs:
\begin{center}
\begin{tikzcd}
    \mathrm{Cat_2} \arrow[r, "\mathrm{End}"] & \mathrm{Alg}_{E_1}(\mathrm{Cat_2}) \arrow[r, "\mathcal{Z}"] & \mathrm{Alg}_{E_2}(\mathrm{Cat_2}) \arrow[r, "\mathcal{Z}_2"] & \mathrm{Alg}_{E_3}(\mathrm{Cat_2}) \arrow[r, "\mathcal{Z}_3"] & \mathrm{CAlg(Cat)}
\end{tikzcd}
\end{center}
In this talk we focus on the construction of the Drinfel'd center $\mathcal{Z}$ which assigns a braided fusion 2-category to a fusion 2-category. 

\paragraph{Notation}The assignment $\mathcal{C} \mapsto \Sigma(\mathcal{C})$ from multifusion 1-categories to pointed semisimple 2-categories is defined by first taking the one-object delooping $B\mathcal{C}$ and then taking the Karoubi-completion $B\mathcal{C}^\nabla$. Alternatively, we can describe $\Sigma (\mathcal{C})$ as the “Morita 2-category” $\mathrm{Morita}^\mathrm{sep}(\mathcal{C})$ of separable algebra objects in $\mathcal{C}$ and their bimodule objects. $\Sigma(\mathcal{C})$ is also equivalent to the 2-category $\mathrm{Mod}\text-\mathcal{C}$ of finite semisimple right $\mathcal{C}$-module categories.

\section{Description \#1: abstract definition}

The Drinfel'd center $\mathcal{Z}(A)$ of a monoidal 2-category $\mathcal{A}$ is the monoidal 2-category $\mathrm{End}_{{}_{\mathcal{A}} \! \operatorname{BMod}_\mathcal{A}}(\mathcal{A})$ 
of endomorphisms of $\mathcal{A}$ thought of as an $\mathcal{A}$-bimodule.

If $F$ denotes such an endomorphism, and say $b_0 := F(1)$, then
\begin{equation*}
    F(a)\otimes b_0 = F(a\otimes 1) \cong F(1\otimes a) = b_0 \otimes F(a)
\end{equation*}
which motivates the second description.

\section{Description \#2: half-braidings}

The objects of the Drinfel'd center $\mathcal{Z}(\mathcal{A})$ of a monoidal 2-category $\mathcal{A}$ are pairs $(A, h)$ where $A$ is an object of $\mathcal{A}$ and $h$ is a half-braiding isomorphism
\begin{equation*}
    h_{\lvert X} : A\boxtimes X \to X\boxtimes A
\end{equation*}
that depends 2-naturally and monoidally on $X\in \mathcal{A}$. Explicitly, $h$ is
\begin{enumerate}
    \item 2-natural, \textit{i.e.}, for any 1-morphism $\alpha : X\to Y$ in $\mathcal{A}$, we have the data of an isomorphism 
\begin{equation*}
    h_{\lvert X} \circ (\alpha\boxtimes \mathrm{id}) \cong (\mathrm{id}\boxtimes \alpha) \circ h_{\lvert X}
\end{equation*}
that is compatible with composition and unit,
    \item monoidal, \textit{i.e.}, we have isomorphisms
\begin{equation*}
    h_{\lvert X \boxtimes Y} \cong (\mathrm{id} \boxtimes h_{\lvert Y}) \circ (h_{\lvert X} \boxtimes \mathrm{id})
\end{equation*}
that are compatible with the tensor product and are natural.
\end{enumerate}

The 1-morphisms of $\mathcal{Z}(\mathcal{A})$ from $(A, h)$ to $(A',h')$ are pairs $(f,\gamma)$ where $f : A \to A'$ is a 1-morphism in $\mathcal{A}$, and $\gamma$ is an isomorphism
\begin{equation*}
    \gamma_{\lvert X}: (\mathrm{id} \boxtimes f) \circ h_{\lvert X} \cong h'_{\lvert X} \circ (f\boxtimes \mathrm{id})
\end{equation*}
that depends naturally and monoidally on $X\in \mathcal{A}$.

The 2-morphisms of $\mathcal{Z}(\mathcal{A})$ from $(f, \gamma)$ to $(f', \gamma')$ are 2-morphisms in $\mathcal{A}$ from $f$ to $f'$ that are compatible with $\gamma$ and $\gamma'$. 

\section{Description \#3: right braided module categories}

For any braided fusion 1-category $\mathcal{B}$, we denote by $\mathrm{BrMod}\text-\mathcal{B}$ the 2-category of braided right $\mathcal{B}$-module categories.

\paragraph{Theorem. }\textit{There is an equivalence of 2-categories: 
\begin{equation*}
    \mathcal{Z}(\Sigma(\mathcal{B})) \cong \mathrm{BrMod}\text-\mathcal{B}.
\end{equation*}}

\noindent \textit{Proof. }The proof was sketched on the blackboard. 

\paragraph{Examples:}
\begin{enumerate}
    \item A braided monoidal functor $F : \mathcal{B} \to \mathcal{C}$ induces a braiding on $\mathcal{C} \in \mathrm{Mod}\text-\mathcal{B}$: 
    \begin{equation*}
        b\triangleright c \mapsto \mathrm{br}^\mathcal{C}\mathrm{br}^\mathcal{C}F(b)\otimes c.
    \end{equation*}
    \item If $F = \mathrm{id}$, then $\mathcal{B}\in \mathrm{BrMod}\text-\mathcal{B}$ using $(\mathrm{br}^\mathcal{B})^{\circ 2}$ 
\end{enumerate}

\paragraph{Corollary.}\textit{If $\mathcal{B}$ is a braided fusion 1-category, then $\Omega\mathcal{Z}(\Sigma(\mathcal{B})) \cong \mathcal{Z}_2(\mathcal{B})$.}\\

\noindent \textit{Proof. }Consider $F : (\mathcal{B}, \mathrm{br}^{\circ 2}) \to (\mathcal{B}, \mathrm{br}^{\circ 2})$. $b_0 = F(1)$, such that for all $b\in \mathcal{B}$, $b_0 \otimes \mathrm{br}_{1,b}^{\circ 2} \cong \mathrm{br}_{b_0, b}^{\circ 2}$. Therefore $b_0 \in \mathcal{Z}_2(\mathcal{B}).$ \\

Recall that if $\mathcal{C}$ is a fusion 1-category, then $\mathcal{Z}(\mathcal{C})$ is a braided fusion 1-category.

\paragraph{Proposition.}\textit{If $\mathcal{B}$ is a braided fusion 1-category, then $\mathcal{Z}(\Sigma(\mathcal{B}))$ is a braided fusion 2-category.}\\

\noindent \textit{Proof. }We know from the previous theorem that $\mathcal{Z}(\Sigma(\mathcal{B})) \cong \mathrm{BrMod}\text-\mathcal{B}$. Therefore, using the properties of $\mathrm{Mod}\text-\mathcal{B}$, we have that $\mathcal{Z}(\Sigma(\mathcal{B}))$ is additive, has adjoints and is idempotent complete. Moreover, given that the functor $\mathrm{BrMod}\text-\mathcal{B}\to \mathrm{Mod}\text-\mathcal{B}$ is locally fully faithful, $\mathcal{Z}(\Sigma(\mathcal{B}))$ is locally finite semisimple. We will see later that $\mathcal{Z}(\mathcal{B})$ is a generator of $\mathcal{Z}(\Sigma(\mathcal{B}))$, and that $\mathcal{Z}(\Sigma(\mathcal{B}))$ is braided.

\section{Description \#4: half-braided algebras}

Let $\mathcal{B}$ denote a braided fusion 1-category. Recall from the previous talk:

\paragraph{Theorem.}\textit{The 2-category $\mathrm{sHBA}(\mathcal{B})$ of separable half-braided algebras is equivalent to $\mathrm{BrMod}\text-\mathcal{B}$.}\\

\begin{figure}[h]
        \centering
        \includegraphics[width=0.5\linewidth]{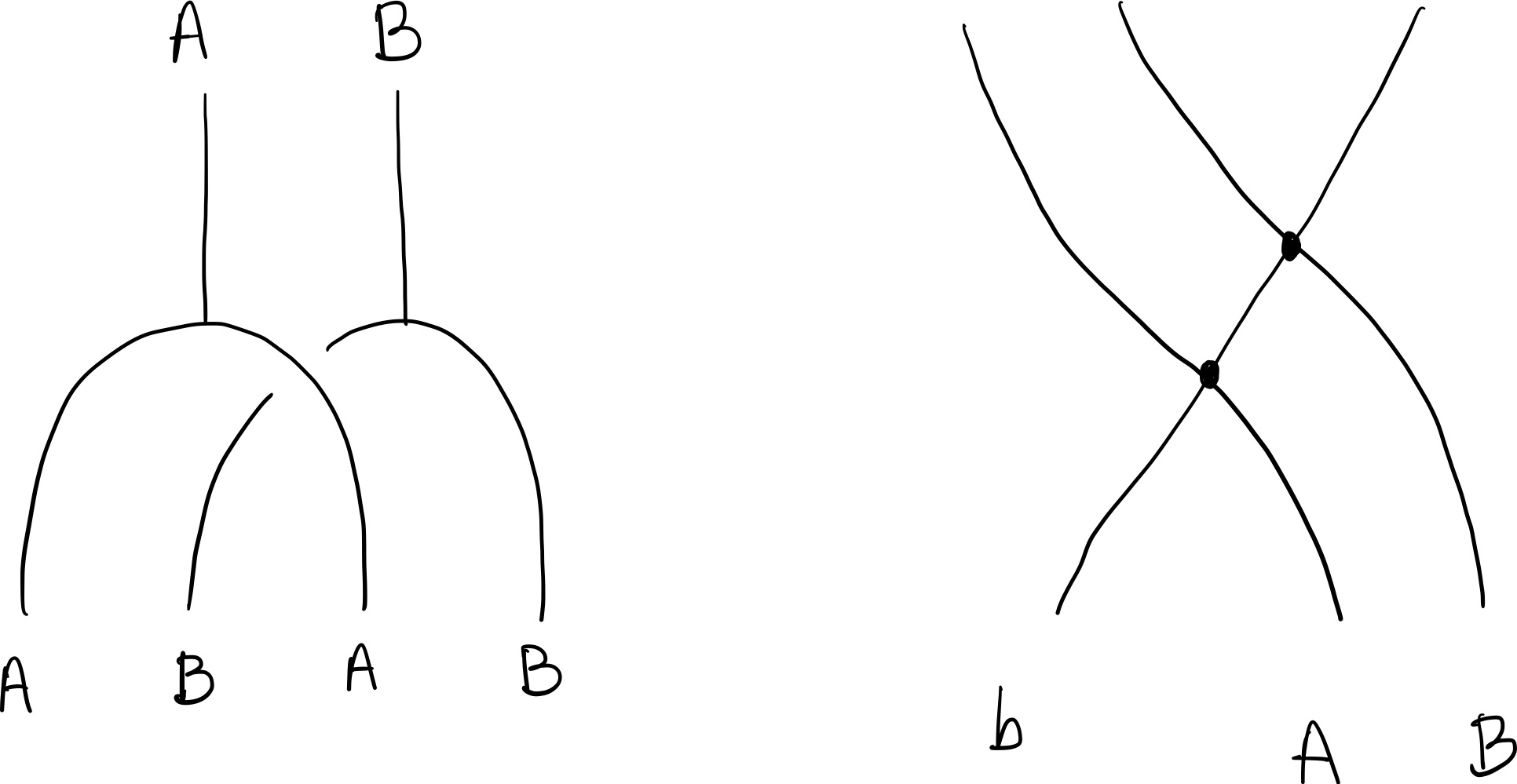}
        \caption{Product and half-braiding in $(A, \gamma)\boxtimes (B,\delta)$}
        \label{fig:halfbraiding}
\end{figure}

Therefore, both $\mathrm{sHBA}(\mathcal{B})$ and $\mathrm{BrMod}\text-\mathcal{B}$ can serve as models for $\mathcal{Z}(\Sigma(\mathcal{B}))$. The monoidal structure on $\mathrm{sHBA}(\mathcal{B})$ is given by:
\begin{enumerate}
    \item unit $(1, \mathrm{id})$,
    \item $(A, \gamma) \boxtimes (B,\delta) = (A\boxtimes B, \gamma\boxtimes\delta)$ with product and half-braiding as given in Figure \ref{fig:halfbraiding}.
\end{enumerate}

\paragraph{Lemma.}\textit{There is an isomorphism of half-braided algebras 
\begin{equation*}
    (A, \gamma)\boxtimes (B,\delta) \cong (B, \delta)\boxtimes (A, \gamma).
\end{equation*}}

\paragraph{Remark.} $\Omega\mathcal{Z}(\Sigma\mathcal{B}) \cong \mathcal{Z}_2(\mathcal{B})$ because any $b\in \mathcal{B}$ is a 1-1-bimodule, and this bimodule structure is half-braided iff $b$ is transparent.

\paragraph{Corollary.}\textit{$\mathcal{Z}(\mathcal{B})$ is a generating object of $\mathrm{BrMod}\text-\mathcal{B}$. In particular, there is an equivalence of 2-categories:
\begin{equation*}
    \Sigma\Omega_{\mathcal{Z}(\mathcal{B})}\mathrm{BrMod}\text-\mathcal{B} \cong \mathcal{Z}(\Sigma\mathcal{B}).
\end{equation*}}

\noindent\textit{Proof. }For $M\in \mathrm{BrMod}\text-\mathcal{B}$, we need to find a non-zero braided module functor $M\to \mathcal{Z}(\mathcal{B})$. Write $M = \mathrm{Mod}(A)$ with $A$ a half-braided algebra. As sketched in the previous talk, the half-braiding on $A$ induces a half-braiding on any $A$-module, so we obtain the horizontal functor in the following commutative diagram:
\begin{center}
\begin{tikzcd}
    \mathrm{Mod}(A) \arrow[r] & \mathcal{Z}(\mathcal{B})\\
    \{A\} \arrow[u] \arrow[ur]
\end{tikzcd}
\end{center}
This is indeed a braided module functor, which is evidently non-zero.

\begin{figure}
\centering
\begin{minipage}{.5\textwidth}
  \centering
  \includegraphics[width=3cm]{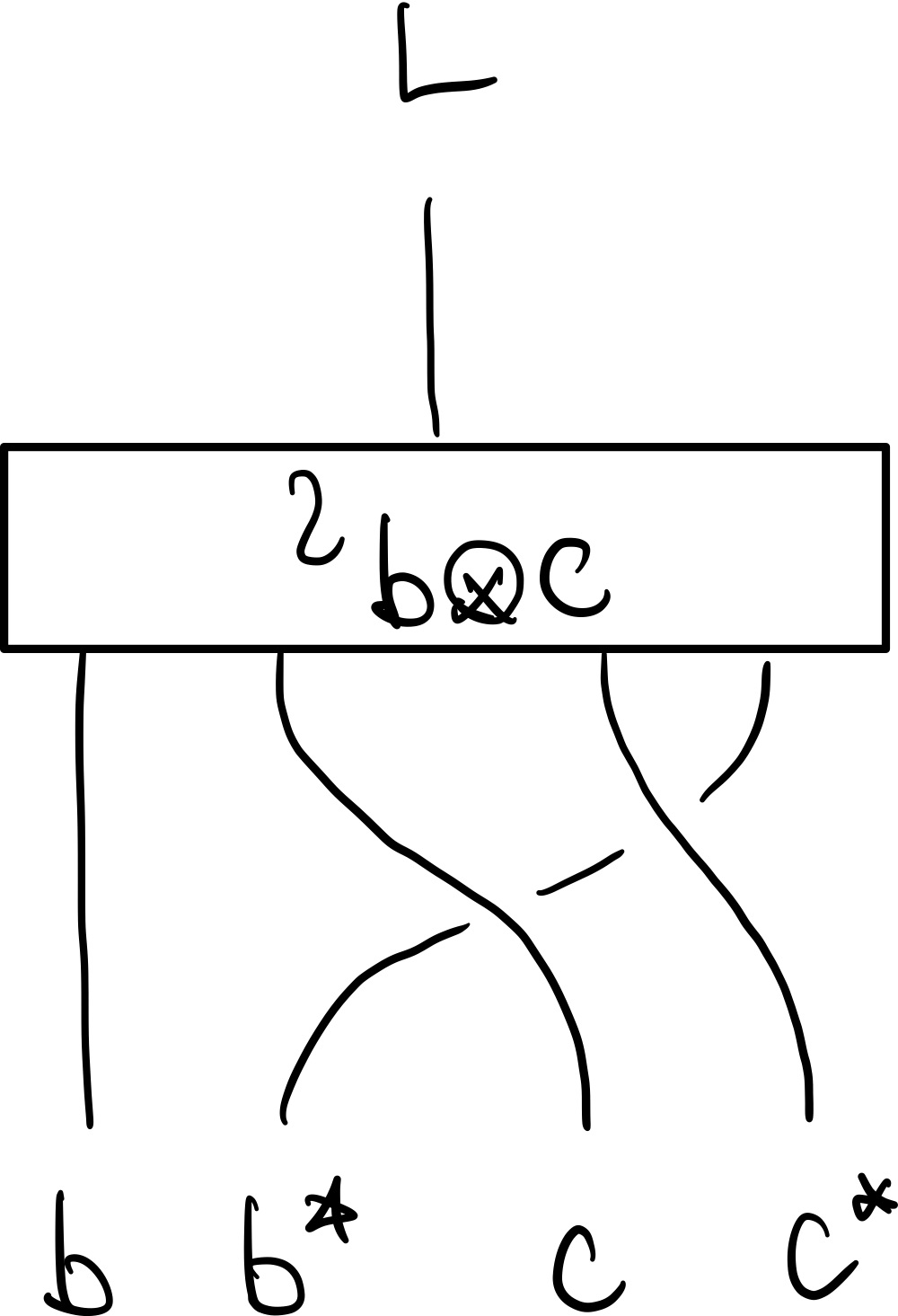}
  \caption{Multiplication in $L$}
  \label{fig:exmultiplication}
\end{minipage}%
\begin{minipage}{.5\textwidth}
  \centering
  \includegraphics[width=3.9cm]{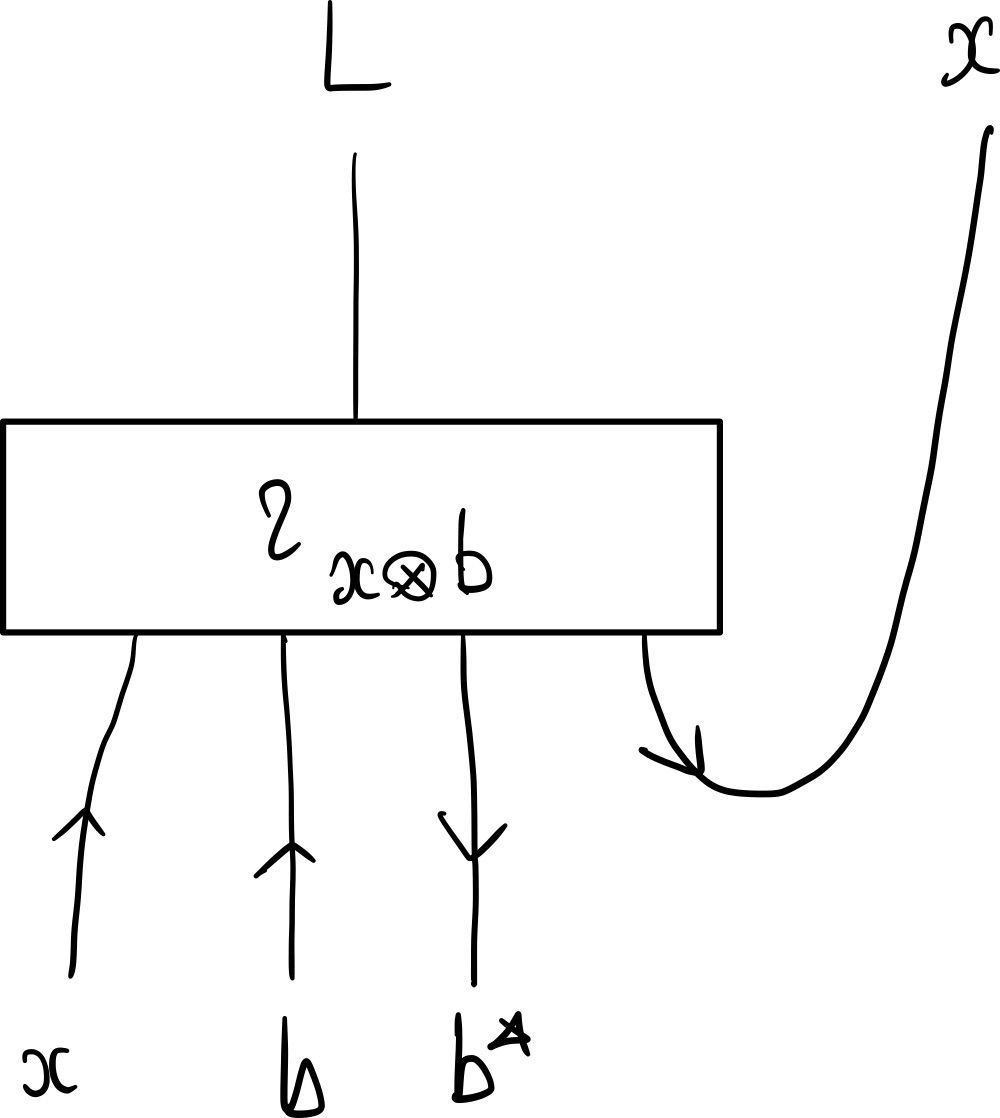}
  \caption{Half-braiding in $L$}
  \label{fig:exhalfbraiding}
\end{minipage}
\end{figure}

\paragraph{Example.} What is the half-braided algebra associated to $\mathcal{Z}(\mathcal{B}) \in \mathrm{BrMod}\text-\mathcal{B}$? Recall that if $F : \mathcal{B} \times \mathcal{B}^\mathrm{op} \to \mathcal{B}'$, then its coend is defined as 
\begin{center}
\begin{tikzcd}
    \displaystyle\int^\mathcal{B}F := \mathrm{coeq}(\bigoplus_{b\in \mathcal{B}}F(b,b) & \displaystyle\bigoplus_{f : b\to b'}F(b,b') \arrow[l, yshift=0.1cm] \arrow[l, yshift=-0.1cm])
\end{tikzcd}
\end{center}
and $\iota_b : F(b,b) \to \int^\mathcal{B}F$ denotes the inclusion. Notice that if $\mathcal{B}$ is semisimple, then
\begin{equation*}
    \int^\mathcal{B}F \cong \bigoplus_{[b] \in \pi_0\mathcal{B}}F(b,b).
\end{equation*}
$\displaystyle L:=\int^\mathcal{B}b\otimes b^*$ is an algebra with 
\begin{enumerate}
    \item unit $\iota_1 : 1 = 1\otimes 1^* \to L$,
    \item multiplication $L\otimes L \to L$ given by the dinatural (in $b$ and $c$) transformation in Figure \ref{fig:exmultiplication}
    \item half-braiding $x \otimes L \to L \otimes x$ for $x \in \mathcal{B}$ given by the dinatural (in $b$) transformation in Figure \ref{fig:exhalfbraiding}.
\end{enumerate}

\paragraph{Proposition. }\textit{The functor $L\text-\mathrm{Mod}_\mathcal{B} \to \mathcal{Z}(\mathcal{B})$ induced from the half-braiding on $L$ is an equivalence of braided module categories.}\\

\begin{figure}
    \centering
    \includegraphics[width=0.25\linewidth]{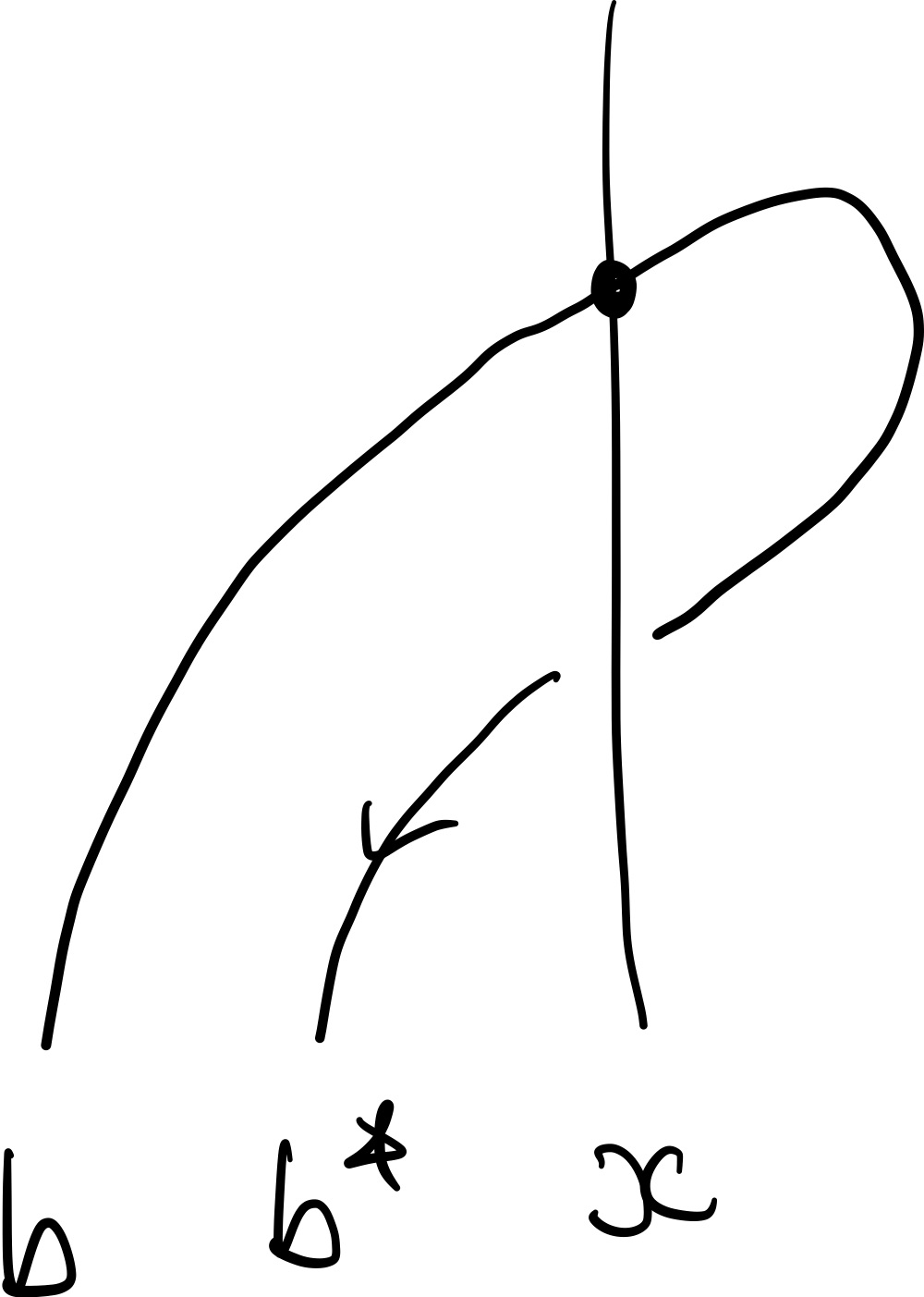}
    \caption{Left $L$ action on $x$}
    \label{fig:inverse}
\end{figure}

\noindent\textit{Proof. }The functor $L\text-\mathrm{Mod}_\mathcal{B} \to \mathcal{Z}(\mathcal{B})$ has inverse $(x,\gamma) \mapsto x$ with left $L$ action defined by the dinatural transformation in Figure \ref{fig:inverse}.

%% file: talks/2.3/main.tex
Talk by Adrià Marín Salvador, notes by Matthias Vancraeynest.

\section{Introduction}
The goal of this talk is to discuss how to control non-degenerate extensions of a braided fusion $1$-category $\mathcal{B}$ by means of Lagrangian algebras in $\mathcal{Z}(\Sigma \mathcal{B})$. One can understand this problem as the categorification of the problem of detecting Drinfeld centers in the $1$-categorical setting. 
We will introduce the notion of a Lagrangian algebra in a braided fusion $2$-category and its relation to
non-degenerate extensions. We start by recalling the definition of a Lagrangian algebra in the $1$-categorical setting, which will be starting point for the categorification procedure.
\begin{definition}
    Let $\mathcal{B}$ be a non-degenerate braided fusion category. A commutative algebra $A \in \mathcal{B}$ is called Lagrangian if it is 
    \begin{itemize}
        \item  connected, i.e. $\dim_{\mathbb{C}} \Hom_{\mathcal{B}}(I, A) = 1$,
        \item separable, i.e. the multiplication map admits a section as an $A$-bimodule morphism,
        \item the category of local (aka dyslectic) modules is trivial.
    \end{itemize}
\end{definition}
\section{Lagrangian algebras in braided fusion 2-categories}
In this section we will categorify the three notions appearing in the definition of a Lagrangian algebra.\\
\\
Let $\mathcal{C}$ be a monoidal $2$-category. Let $A \in \mathcal{C}$ be an algebra with multiplication $m$, unit $u$, associator $\alpha$ and unitors $\lambda$ and $\rho$, satisfying the usual coherence conditions. Additionally, if $\mathcal{C}$ is braided, then we say that the algebra $A$ is braided, if it comes equipped with a $2$-isomorphism $\beta: m \circ \br_{A,A} \Rightarrow m$, satisfying some coherence conditions.
\subsection{Separability}
The associator equips the multiplication $m: A \boxtimes A \to A$ with the structure of an $A-A$-bimodule morphism.
\begin{definition}
An algebra A is called rigid, if the multiplication $m$ has a right adjoint $\Delta: A \to A\boxtimes A$ as an $A-A$-bimodule $1$-morphism.
\end{definition}
\begin{definition}
    An algebra A is called separable, if it is rigid and $ev_m: m  \circ \Delta \Rightarrow id_A $ admits a section as an $A-A$-bimodule $2$-morphism.
\end{definition}
In characteristic zero, we expect these two notions to be equivalent, just as multifusion $1$-categories are separable in characteristic zero.

\begin{lemma}
    Let $\mathcal{C} :=$  Bimod-$\mathcal{X}$ be the $2$-category of finite semisimple bimodules over some fusion $1$-category $\mathcal{X}$. Algebra objects $A \in \mathcal{C}$ correspond to finite semisimple monoidal categories $\mathcal{Y}$ equipped with a monoidal functor $\mathcal{X} \to \mathcal{Y}$. Under this correspondence, the algebra $A$ is rigid if and only if it corresponding category $\mathcal{Y}$ is rigid.
\end{lemma}
\subsection{Triviality of local modules}
Let $\mathcal{C}$ be a braided monoidal $2$-category and $A$ a braided algebra in $\mathcal{C}$. Consider a $1$-morphism $f:I \to A$. We now consider the following endomorphism $ m \circ (f \boxtimes \id_A) : A \to A$. We consider two different $2$-isomorphisms $m\circ(f \boxtimes \id_A) \cong m\circ \br_{A,A} \circ \br_{A,A} \circ (f \boxtimes \id_A)$. One uses the braided monoidal structure of the ambient $2$-category $\mathcal{C}$, i.e.  $$f \boxtimes \id_A \cong \br_{A,A} \circ (\id_A \boxtimes f) \cong \br_{A,A} \circ \br_{A,A} \circ (f \boxtimes \id_A).$$
The other $2$-isomorphism, comes from the algebra itself being braided. It uses
$$ m \cong m \circ \br_{A,A} \cong m \circ \br_{A,A} \circ \br_{A,A}.$$

\begin{definition}
    A morphism $f:I \to A$ is called transparent if the two $2$-isomorphisms described above are equal. We define the Müger center of a braided algebra object $A$ as the full subcategory of transparent $1$-morphisms in $\Hom_\mathcal{C}(I,A)$. We denote it by $\mathcal{Z}_2(A)$.
\end{definition}
\begin{remark}
    When we take $\mathcal{C} = 2\Vect$, then a braided algebra $A$ is a braided monoidal category and $\mathcal{Z}_2(A)$ coincides with the usual definition of the Müger center of a braided monoidal category.
\end{remark}
We now describe a different way of obtaining $\mathcal{Z}_2(A)$, more resembling the local modules from the $1$-categorical case. We begin by introducing the $2$-category $\leftindex_{\mathcal{C}} \BrMod-A$, which serves as categorification of the category of local modules.
\begin{remark}
    Let $\leftindex_{\mathcal{C}} \BrMod-A$ be the braided monoidal $2$-category of braided right module objects over A, i.e. module
objects $M$ with act : $M \boxtimes A \to M$ in $\mathcal{C}$, equipped with a $2$-isomorphism $\sigma : \text{act} \circ \br_{A,M} \circ\br_{M,A} \cong \text{act}$,
satisfying some coherence conditions.\\
    We then have that $\mathcal{Z}_2(A) \cong \End_{\leftindex_{\mathcal{C}} \BrMod-A}(A,\beta \circ \beta)$.
\end{remark}
\subsection{Definition of Lagrangian algebra in BF2C}

\begin{definition}\label{Lagdef}
    A Lagrangian algebra $A$ in a braided fusion $2$-category $\mathcal{C}$ is a braided algebra object $A$ being
    \begin{itemize}
        \item  rigid,
        \item  nondegenerate, i.e. $\mathcal{Z}_2(A) \cong \Vect$,
        \item  strongly connected, meaning that the unit $u: I \to A$ is fully faithful or equivalently that $u:I \to A$ is a simple summand of $A$.
    \end{itemize}
\end{definition}
\begin{remark}
    In fact, such an object should only be called ``Lagrangian'', if the ambient $2$-category $\mathcal{C}$ is nondegenerate. In this case, the nondegeneracy condition in definition \ref{Lagdef} implies the triviality of the category of local modules.
\end{remark}
\begin{remark}
    The connectivity condition can be weakened to the requirement that $u: I \to A$ should be simple as a $1$-morphism. Doing so, results in a different notion of a Langrangian algebra.
\end{remark}
\begin{ex}
    Let $\mathcal{C}$ be $2\Vect$, then the Lagrangian algebras are the nondegenerated braided fusion $1$-categories. Note that in this case the strong connectivity condition and its weaker version are in fact equivalent.
\end{ex}
\section{Classification of extensions via Lagrangian algebras}
\begin{prop}
    Let $\mathcal{B}$ be a braided fusion $1$-category. There is a bijective correspondence between isomorphism classes of rigid braided algebras $A\in \mathcal{Z}(\Sigma \mathcal{B})$ and braided
multifusion categories $\mathcal{M}$  with a braided monoidal functor $i : \mathcal{B} \to \mathcal{M}$, up to equivalence. Moreover, \begin{enumerate}
    \item  $i : \mathcal{B} \to \mathcal{M}$ is fully faithful if and only if $u: I \to A$ is fully faithful,
    \item  $\mathcal{Z}_2(\mathcal{B} \to \mathcal{M}) \cong \Hom_{\mathcal{Z}(\Sigma \mathcal{B})}(I,A)$ and $\mathcal{Z}_2(\mathcal{M}) \cong \mathcal{Z}_2(A)$.
\end{enumerate}
\end{prop}
\begin{cor}
In particular, there is a bijective correspondence between the isomorphism classes of Lagrangian algebras in $\BrMod$-$\mathcal{B}$ and nondegenerate braided extensions $\mathcal{B} \subseteq \mathcal{M}$. 
\end{cor}
\begin{proof}[Sketch of proof]
    We use the following equivalence
    \begin{align*}
        \mathcal{Z}(\Sigma \mathcal{B}) &\cong \BrMod\text{-}\mathcal{B}\\
        A &\mapsto \Hom_{\mathcal{B}\text{-}\Mod}(\mathcal{B}, A).
    \end{align*}
    Hence, braided algebra objects correspond to finite semi-simple monoidal categories $\mathcal{M}$ with $i : \mathcal{B} \to \mathcal{M}$ a braided tensor functor. As seen before, the algebra $A$ is rigid, if and only if the corresponding $\mathcal{M}$ is rigid. Which leads to the one-to-one correspondence\\
    \\
    The fully faithfulness of $u: I \to A$ is equivalent with  $u: I \to A$ being a direct summand of $A$, which is, in turn, equivalent with $i: \mathcal{B} \to \mathcal{M}$ being the inclusion of an indecomposable module category. In other words, this last statement says that $i: \mathcal{B} \to \mathcal{M}$ is fully faithful, as was needed to be shown.
    \\
    \\
    For part $2$, we consider $N$, a braided module category over $\mathcal{B}$ with braiding $\{\sigma_{n,b}: n*b \to n*b\}$. A braided module functor $F: \mathcal{B} \to \mathcal{N}$ corresponds to an object $n := F(I_{\mathcal{B}})$, such that $\sigma_{n,b} = \id_b$ for all $b\in \mathcal{B}$. This allows us to write 
    \begin{align*}
        \Hom_{\BrMod\text{-}\mathcal{B}}(I,A) \cong \{m \in \mathcal{M} | \sigma_{m,b} = \id_b \forall b\} =: \mathcal{Z}_2(\mathcal{B} \to \mathcal{M}),
    \end{align*}
    where the last equality holds per definition.
    Furthermore, by a similar reasoning, you can convince yourself that $\mathcal{Z}_2(\mathcal{M}) \cong \mathcal{Z}_2(A)$.   
\end{proof}
As a corollary, we have the following theorem.
\begin{theorem}
    Let $\mathcal{B}$ and $\mathcal{C}$ be braided fusion $1$-categories. If $\mathcal{Z}(\Sigma \mathcal{B}) \cong \mathcal{Z}(\Sigma \mathcal{C})$, then there exists a nondegenerate extension $i:\mathcal{B} \to \mathcal{M}$ such that $\mathcal{Z}_2(\mathcal{B} \subset \mathcal{M})\cong \mathcal{C}^{\text{rev}}$.
\end{theorem}
\begin{proof}
    Let $\mathcal{L}$ be the lagrangian algebra in $\mathcal{Z}(\Sigma \mathcal{C})$ corresponding to the extension $\mathcal{C} \hookrightarrow \mathcal{Z}(\mathcal{C})$. Denote $F:\mathcal{Z}(\Sigma \mathcal{C}) \to \mathcal{Z}(\Sigma \mathcal{B})$ the equivalence provided as hypothesis. Then $F(\mathcal{L}) \in \mathcal{Z}(\Sigma \mathcal{B})$ is a Lagrangian algebra and hence has an associated nondegenerate extension $\mathcal{B} \to \mathcal{M}$. One then obtains,
    \begin{align*}
        \mathcal{Z}_2(\mathcal{B} \subset \mathcal{M}) \cong \Hom_{\mathcal{Z}(\Sigma \mathcal{B})}(I, F(\mathcal{L}) \cong \Hom_{\mathcal{Z}(\Sigma \mathcal{C})}(I, \mathcal{L}) \cong \mathcal{Z}_2(\mathcal{C} \subset \mathcal{Z}(\mathcal{C})) \cong \mathcal{C}^{\text{rev}}.
    \end{align*}
\end{proof}
\begin{remark}
    Why is this theorem helpful? Assume $\mathcal{Z}_2(\mathcal{B}) =\sVec$. In order to find a minimal nondegenerate extension of $\mathcal{B}$, it suffices to find an equivalence $\mathcal{Z}(\Sigma \mathcal{B})  \cong \mathcal{Z}(\Sigma \sVec)$. This provides $\mathcal{B} \to \mathcal{M}$, a nondegenerate extension with $\mathcal{Z}_2(\mathcal{B} \subset \mathcal{M}) \cong \sVec^{\text{rev}} \cong \mathcal{Z}_2(\mathcal{B)}$.
\end{remark}

%% file: talks/3.1/content.tex
Talk by Alea Hofstetter, notes by Julia Bierent.

\section{Introduction}

Let $\CC$ be a braided fusion category.
\begin{notation}
    \begin{itemize}
        \item Let $\sigma(\CC)$ be the set of simple objects in $\CC$, and
        \item let $d_{\pm}(X)=\Tr_\pm(\id_X)$. (See talk 1.1. by Theodoros Lagiotis on Braided fusion and Drinfeld center).
    \end{itemize}
\end{notation}
We can define the $S$-matrix $S:\sigma(\CC)\times\sigma(\CC)\to \Kb$. For $X,Y \in \sigma(\CC)$,
\begin{equation*}
    S_{X,Y}=m\frac{1}{d_-(X)d_+(Y)} 
    \begin{tikzpicture}[baseline=0]
        \coordinate (a) at (0,1);
        \coordinate (b) at (0,-1);
        \coordinate (l) at (-1.5,-0.25);
        \coordinate (r) at (1.5,0.25);
        \node (A) at (a) {};
        \node (B) at (b) {};
        \node at (-0.5,-1) {$X$};
        \node at (0.5,-0.9) {$Y$};
        \node (C) at (-1.75,0) {};
        \draw [-] (A) to[out=45,in=90] (r);
        \draw [-] (r) to[out=-90, in=-45] (b);
        \draw [-] (b) to[out=135, in=-135] (A);
        \draw [-] (B) to [out=-135,in=-90] (l);
        \draw [-] (l) to [out=90,in=135] (a);
        \draw [-] (a) to [out=-45,in=45] (B);
    \end{tikzpicture}
\end{equation*}
We want to build an $S$-matrix now for braided fusion $2$-categories.
\section{Braiding in braided fusion $2$-category}
Let $\CB$ be a braided fusion $2$-category.
\begin{notation}
    We will denote by $\pi_0 \CC$ the set of components in a $2$-category $\CC$, i.e., the set of simple objects in $\CC$ quotiented by the relation of being in the same component. Two simple objects $A,B \in \CC$ are in the same component if there exists a non-zero $1$-morphism $A\to B$.
\end{notation}
The entries of our $S$-matrix will now be in $\pi_0\CB$ and $\pi_0(\Omega\CB)=\pi_0(\Hom_\CB(I,I)) =$ set of simple objects in $\Omega\CB$ in this case.
So we want to define
\begin{equation*}
    \begin{aligned}
        S:\pi_0\CB\times &\pi_0\Omega\CB \to \Kb \\
        (A, &b)
    \end{aligned}
\end{equation*}
We have the braiding
\begin{equation*}
    \br_{A,I}:A\square I \to I\square A,
\end{equation*}
and the following diagram
\begin{equation*}
    \begin{tikzcd}
        A\square I \arrow[r,"A\square b"] \arrow[d, "\br_{A,I}"']
        &A\square I \arrow[d,"br_{A,I}"] \arrow[ld,Rightarrow,"\br_{A,b}"]\\
        I\square A \arrow[r,"b\square A"']
        &I\square A
    \end{tikzcd}
\end{equation*}
\begin{figure}[h]
    \centering
    \begin{tikzpicture}[baseline=-25]
        \node at (0,0) {$S(A,b)=\frac{1}{d_+(b)}\cdot$};
    \end{tikzpicture}
    \includegraphics[width=0.2\textwidth]{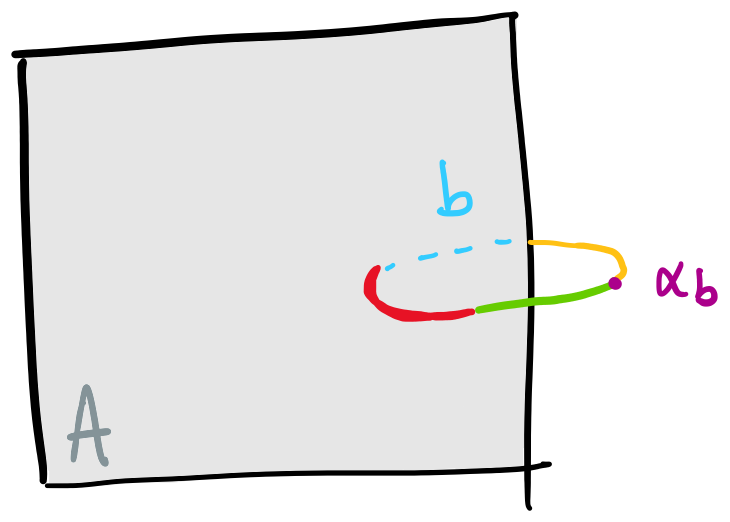}
    \caption{picture of the $S-$matrix, $b$ flying around the surface $A$.}
    \label{fig:baroundA}
\end{figure}
So we can define the $S-$matrix elements by
\colorlet{Myyellow}{green!10!orange}
\colorlet{Mygreen}{green!60!orange}
\begin{equation*}
     S(A,b)=\frac{1}{d_+(b)}(\id_A\square \color{Myyellow} \ev_b \color{black})(\color{red}(\br_{b,A}\cdot\br_{A,b})\color{black}\otimes \color{violet}\alpha_b\color{black})\cdot (\id_A \square \color{Mygreen}\coev_{^*b}\color{black}),
\end{equation*}
where $\alpha_b$ is some duality data, for a picture, see Figure \ref{fig:baroundA}.
\section{Application}
Let $\CC$ be a braided fusion $1$-category. Then $\CB:=\CZ(\Sigma \CC)$ is a braided fusion $2$-category. We are going to build an $S$-matrix on $\CB$
\begin{equation*}
    \CZ(\Sigma\CC)\cong \CZ(\Mod(\CC))\cong \BrMod(\CC) \cong \sHBA(\CC),
\end{equation*}
\begin{equation*}
    \Omega(\CZ(\Sigma\CC))\cong \CZ_2(\CC).
\end{equation*}
We are going to define $S:\pi_0\CB\times\pi_0\Omega\CB\to \Kb$ using two different descriptions of $\CB$.
\begin{enumerate}
    \item $S:\pi_0\BrMod(\CC)\times\pi_0\CZ_2(\CC)\to \Kb$.
    
    Let $[(\CM,\sigma)]\in\pi_0\BrMod(\CC), m\in\CM$ and $ [b] \in \pi_0$. Then,
    \begin{equation*}
        S_{[\CM,\sigma],[b]} = \frac{1}{d_+(b)}(m*\ev_b)\circ(\sigma_{m,b}*\alpha_b)\circ(m*\coev_{*b}),
    \end{equation*}
    which is well defined for different choices of object $m\in\CM$ and all different elements in the classes $[(\CM,\sigma)]$ and $[b]$.
    \item $S:\pi_0\sHBA(\CC)\times\pi_0\CZ_2(\CC)\to\Kb$

    Let $[(A,\gamma)]\in\pi_0\sHBA(\CC)$ and $[b]\in\pi_0\CZ_2(\CC)$. Then,
    \begin{equation*}
        \frac{1}{d_+(b)}(\id_A\otimes\ev_b)\circ((\gamma_{b,A}\cdot\br_{A,b})\otimes\alpha_b)\circ(\id_A\otimes\coev_{*b}),
    \end{equation*}
    which is well defined for all different elements in the classes $[(A,\gamma)]$ and $[b]$.
\end{enumerate}
\begin{example}
    Let $G$ be a finite abelian group and $\CC=\vect_G$ the category of $G$-graded $\Cb$-vector spaces. Write $U_g$ for a simple object in degree $g\in G$. Let $[\Psi,\Omega)]\in H^3_{ab}(G,\Cb^\times)$ with
    \begin{equation*}
        \Psi : G\times G\times G \to \Cb^\times,
    \end{equation*}
    giving the associator, and
    \begin{equation*}
        \Omega : G\times G \to \Cb^\times,
    \end{equation*}
    giving the braiding on $\CC$ by
    \begin{equation*}
        \br_{U_g,U_h}\cong\Omega(g,h)\id_{U_{g+h}},
    \end{equation*}
    for $g,h \in G$. Then,
    \begin{equation*}
        \begin{aligned}
            \pi_0\CZ_2(\CC) = \{U_g\mid \Omega(g,h)\Omega(h,g)=1 \forall h\in G\}, \\
            \CZ_2(G,\Omega) = \{g\in G \mid U_g \in \CZ_2(\CC)\}.
        \end{aligned}
    \end{equation*}
    Let us now look at $\Mod(\CC)$. Let $H\leq G$, $\mu:H\times H\to \Cb^\times$ and $d\mu = d\Phi$, then we can build $\CM_{H,\mu}\in\Mod(\CC)$. Let $M_k$ a simple in $\CM_{H,\mu}$ with $[k]\in G/H$. A braiding exists on $\CM_{H,\mu}$ if and only if $H\leq Z_2(G,\Omega)$. So let $H\leq Z_2(G,\Omega)$, then $(\CM_{H,\mu},\sigma)$ is a braided module category with braiding defined by
    \begin{equation*}
        \begin{aligned}
             \sigma_{M_k,U_g}&:M_k,U_g:M_k*U_g\to U_g*M_k\\
             &:=\Omega(k,g)\cdot\Omega(g,k)\cdot\chi(g)\cdot\id_{M_k*U_g},
        \end{aligned}
    \end{equation*}
    with $\chi$ a character on $G$.

    We still need to find the Schur equivalence classes. For all $(\CM_{H,\mu},\sigma_\chi)$, we can build a braiding $\hat{\sigma}_\chi$ also defined by the character $\chi$ but on $\CC_\CC$ the regular $\CC$-module category (take $H$ trivial) so that
    \begin{equation*}
        [(\CM_{H,\mu},\sigma_\chi)]=[(\CC_\CC,\hat{\sigma}_{\chi})],
    \end{equation*}
    Now $[(\CC_\CC,\sigma_\chi)]=[(\CC_\CC,\sigma_{\chi'}]$ if and only if $\chi_{|\CZ_2(G,\Omega)}=\chi_{|Z_2(G,\Omega)}$. This yields that $S$ is in fact the character table of $\CZ_2(G,\Omega)$.
    \begin{remark}
        If $\CC$ is symmetric, then $\CZ_2(\CC)$ is everything and $S$ gives us the character table of $G$.
    \end{remark}
\end{example}

%% file: talks/3.2/3.2.tex
\newcommand\arXiv[1]{\href{http://arxiv.org/abs/#1}{\nolinkurl{arXiv:#1}}}
\newcommand\DOI[1]{\href{http://dx.doi.org/#1}{\nolinkurl{DOI:#1}}}

\NewDocumentCommand {\cat}    {m  } {\mathcal{#1}}
\NewDocumentCommand {\ccat}   {m  } {\mathcal{#1}}
\NewDocumentCommand {\FVect}  {   } {\mathbf{Vec}}
\NewDocumentCommand {\VVect}  {   } {2\FVect}
\NewDocumentCommand {\field}  {   } {\Bbbk}
\NewDocumentCommand {\lin}    {m  } {\field[#1]}
\NewDocumentCommand {\rev}    {   } {\mathrm{rev}}
\NewDocumentCommand {\tomono} {   } {\hookrightarrow} 
\NewDocumentCommand {\mono}   {O{}} {\ar[hook,#1]} 

Talk by Leon Liu, notes by Jonathan Davies.

\section{Non-degenerate braided fusion 2-categories}

Recall that a braided fusion 1-category \(\cat{B}\) is defined to be non-degenerate if its M\"{u}ger centre \(\cat{Z}_2(\cat{B})\) is equivalent to the trivial symmetric fusion category \(\FVect\).
This holds if and only if its S-matrix is invertible.

\begin{definition}\label{def:ndmap}
  Let \(A,B\) be finite sets, and let \(\phi\) be a map \(A \times B \to \field\).
  Then \(\phi\) is called \textit{non-degenerate} if the induced map \(\lin{A} \to \lin{B}^{\ast}\) is an isomorphism.
\end{definition}
\begin{remark}
  Notice that \(\phi \colon A \times B \to \field\) can only be non-degenerate if \(|A| = |B|\).
\end{remark}

\begin{definition}
  A braided fusion 2-category \(\ccat{C}\) is called \textit{non-degenerate} if the S-matrix \(\pi_0 \ccat{C} \times \pi_0\Omega \ccat{C} \to \field\) is non-degenerate in the sense of Definition \ref{def:ndmap}.
\end{definition}

Since there are multiple equivalent definitions for the non-degeneracy of a braided fusion 1-category, we expect that the following should be\footnote{
  This is likely true, but depends on some higher Morita theory that may or may not exist.
  The equivalence (\ref{nd2} \(\iff\) \ref{nd3}) is in print assuming some TFT calculus (see \cite[Thm. 2]{j2022}).
  The equivalence (\ref{nd1} \(\iff\) \ref{nd2}) is easier, but not in print (and it may still require some TFT calculus).
}
equivalent for a braided fusion 2-category \(\ccat{C}\):
\begin{enumerate}
\item \(\ccat{C}\) is non-degenerate in the above sense of having a non-degenerate S-matrix \label{nd1}
\item the M\"{u}ger sylleptic centre \(\cat{Z}_2(\ccat{C})\) is equivalent to the trivial braided fusion 2-category \(\VVect := \Sigma\FVect\) \label{nd2}
\item \(\ccat{C}\) is \(\mathbb{E}_2\)-Morita invertible. \label{nd3}
\end{enumerate}

\section{Non-degeneracy of Drinfeld centers}

For any fusion 1-category\footnote{
  That is, any fusion 1-category over an algebraically closed field of characteristic 0.
  Over a general field we require the further assumption that \(\cat{C}\) is separable.
} \(\cat{C}\), recall that \(\cat{Z}(\cat{C})\) is a non-degenerate braided fusion 1-category.
For a braided fusion 1-category \(\cat{B}\), we have already seen that \(\cat{Z}(\Mod{\cat{B}}) \simeq \cat{Z}(\Sigma\cat{B})\) is a braided fusion 2-category equivalent to \(\BrMod{\cat{B}}\).
We therefore expect that \(\cat{Z}(\Mod{\cat{B}})\) is also non-degenerate; this is indeed true, and given in Theorem \ref{thm:main}.
We split Theorem \ref{thm:main} into two easier results which we prove first.


\begin{theorem}[2.52 of \cite{JFR}]\label{thm:A}
  Let \(\cat{B}\) be a braided fusion 1-category.
  Then we have
  \begin{equation*}
    |\pi_0\cat{Z}(\Mod{\cat{B}})|
    = |\pi_0\Omega\cat{Z}(\Mod{\cat{B}})|
    = |\pi_0\cat{Z}_2(\cat{B})| \,.
  \end{equation*}
  That is, the simple objects of \(\cat{Z}(\Mod{\cat{B}})\) are in bijection with the simple transparent objects of \(\cat{B}\) up to the equivalence relation
  \begin{equation*}
    x \sim y \iff \Hom_{\cat{C}}(x,y) \neq 0 \,.
  \end{equation*}
\end{theorem}
\begin{proof}
  We have discussed the second equality already.
  Recall that \(\cat{Z}(\cat{B})\) is in \(\BrMod{\cat{B}} \simeq \cat{Z}(\Mod{\cat{B}})\) as a generator.
  It follows by definition that the direct sum decomposition of \(\cat{Z}(\cat{B})\) into indecomposable braided \(\cat{B}\)-module categories contains at least one simple object from each component of \(\BrMod{\cat{B}}\).

  We have to check how many parts \(\cat{Z}(\cat{B})\) gets broken into.
  First, there exist braided monoidal inclusions of \(\cat{B}\) and \(\cat{B}^{\rev}\) into \(\cat{Z}(\cat{B})\) which give the commuting square
  \begin{equation*}
  \begin{tikzcd}
    & \cat{B} \mono[dr]
    \\
    \cat{Z}_2(\cat{B}) \mono[ur] \mono[dr]
    && \cat{Z}(\cat{B})
    \\
    & \cat{B}^{\rev} \mono[ur]
  \end{tikzcd}
  \end{equation*}
  Then there is a further induced monoidal subcategory \(\cat{B} \boxtimes_{\cat{Z}_2(\cat{B})} \cat{B}^{\rev}\) giving a commuting diagram
  \begin{equation*}
  \begin{tikzcd}
    & \cat{B} \mono[dr] \ar[d]
    \\
    \cat{Z}_2(\cat{B}) \mono[ur] \mono[dr]
    & \cat{B} \boxtimes_{\cat{Z}_2(\cat{B})} \cat{B}^{\rev} \mono[r]
    & \cat{Z}(\cat{B})
    \\
    & \cat{B}^{\rev} \mono[ur] \ar[u]
    &\,,
  \end{tikzcd}
  \end{equation*}
  where \(\cat{Z}(\cat{B} \tomono \cat{Z}(\cat{B})) \simeq \cat{B}^{\rev}\) and \(\cat{Z}(\cat{Z}_2(\cat{B}) \tomono \cat{Z}(\cat{B})) \simeq \cat{B} \boxtimes_{\cat{Z}_2\cat{B}} \cat{B}^{\rev}\) are centralisers.\footnote{
    We write \(\cat{Z}(\cat{C} \tomono \cat{E})\) for the \textit{centraliser} of an inclusion \(\cat{C} \tomono \cat{E}\) of braided fusion categories; that is, the full subcategory of objects in \(\cat{E}\) which are transparent to all objects in \(\cat{C}\).
    The M\"uger centre is then precisely the special case \(\cat{Z}_2(\cat{C}) := \cat{Z}(\cat{C} \tomono \cat{C})\).
    See \cite{dno2013}.
  }
  Define the following three equivalence relations on the indecomposable \(\cat{B}\)-module summands of \(\cat{Z}(\cat{B})\):
  \begin{enumerate}
  \item \(\cat{M} \sim_1 \cat{N}\) if they live in the same component of \(\BrMod{\cat{B}}\)
  \item \(\cat{M} \sim_2 \cat{N}\) if there exists a braided \(\cat{B}\)-module functor \(F \colon \cat{Z}(\cat{B}) \to \cat{Z}(\cat{B})\) which restricted to these summands is non-zero
  \item \(\cat{M} \sim_3 \cat{N}\) if they are are in the same \(\cat{B} \boxtimes_{\cat{Z}_2(\cat{B})} \cat{B}^{\rev}\)-module summand of \(\cat{Z}(\cat{B})\).
    By definition, this means there exist \(b \in \cat{B}^{\rev}, m \in \cat{M}, n \in \cat{N}\) such that \(\Hom_{\cat{Z}(\cat{B})}(b \otimes m, n) \neq 0\).
  \end{enumerate}

  Clearly \(\sim_1\) and \(\sim_2\) are equivalent, since any non-zero braided \(\cat{B}\)-module functor between \(\cat{M}\) and \(\cat{N}\) can be extended to a braided \(\cat{B}\)-module functor \(\cat{Z}(\cat{B}) \to \cat{Z}(\cat{B})\) with all other components vanishing.

  Suppose \(\cat{M} \sim_3 \cat{N}\).
  Given \(b \in \cat{B}^{\rev}\) with some \(\Hom_{\cat{Z}(\cat{B})}(b \otimes m,n)\) non-zero, we have \(\cat{M} \sim_2 \cat{N}\) by taking the braided \(\cat{B}\)-module functor  \(b \otimes - \colon \cat{Z}(\cat{B}) \to \cat{Z}(\cat{B})\).
  For the converse, suppose \(\cat{M} \sim_2 \cat{N}\).
  We claim that all braided \(\cat{B}\)-module endofunctors of \(\cat{Z}(\cat{B})\) are summands of functors of the form \(b \otimes -\), where \(b\) has \(\Hom_{\cat{Z}(\cat{B})}(b \otimes m,n)\) non-zero for some \(m \in \cat{M}, n \in \cat{N}\).
  This comes from the fact that \(\End_{\BrMod{\cat{B}}}(\cat{Z}(\cat{B}))\) is equivalent to \(\cat{Z}(\cat{B})\) under the annular tensor product.
  In particular, idempotent completion provides a map from the annular category \(\Ann(\cat{B})\).

  Since \(\cat{Z}(\cat{B})\) is a generator for \(\BrMod{\cat{B}}\), the set of equivalence classes of \(\sim_1\) is equal to \(\pi_0\BrMod{\cat{B}} \cong \pi_0\cat{Z}(\Mod{\cat{B}})\).
  Corollary 3.6 of \cite{dgno2010} gives the following:
  \begin{quote}
    Let \(\cat{C}\) be a fusion 1-category and \(\cat{E}\) a non-degenerate braided fusion 1-category with an embedding \(\cat{C} \tomono \cat{E}\).
    Then the number of \(\cat{Z}(\cat{C} \tomono \cat{E})\)-summands of \(\cat{E}\) (as a \(\cat{Z}(\cat{C} \tomono \cat{E})\)-module) is equal to \(|\pi_0\cat{C}|\).
  \end{quote}
  Identifying \(\cat{C} = \cat{Z}_2(\cat{B})\) and \(\cat{E} = \cat{Z}(\cat{B})\) so that \(\cat{Z}(\cat{C} \tomono \cat{E}) \simeq \cat{B} \boxtimes_{\cat{Z}_2} \cat{B}^{\rev}\), it follows that \(|\pi_0\cat{Z}_2(\cat{B})|\) is the number of equivalence classes of \(\sim_3\).
\end{proof}


Before proving the next theorem, we describe an explicit model for the S-matrix of \(\cat{Z}(\Mod{\cat{B}})\).
Let \(\cat{M}\) be a simple object of \(\BrMod{\cat{B}}\) and let \(m\) be in turn a simple object of \(\cat{M}\).
Let \(b \in \cat{Z}_2(\cat{B})\) be a simple transparent object of \(\cat{B}\).
Writing \(\bullet\) for a map \(b^{\ast\ast} \to b\) (which always exists in this case), we define \((R_b^{\cat{M}})_m \colon m \to m\) by the string diagram
\begin{equation*}
\begin{tikzpicture}
  \draw[thick,looseness=0.75] (-1,0) arc (180:0:1);

  \draw[white, line width=6pt] (-.5,-1.25) -- (-.5,1.25);
  \draw[thick] (-.5,-1.25) -- (-.5,1.25);

  \draw[white, line width=6pt, looseness=0.75] (-1,0) arc (180:360:1);
  \draw[thick, looseness=0.75] (-1,0) arc (180:360:1);

  \draw[thick] (-1.12,-0.1) -- (-1.0,0.1) -- (-0.88,-0.1); 
  \draw[thick] (-0.62,-0.1) -- (-0.5,0.1) -- (-0.38,-0.1); 

  \node at (1,0) {\(\bullet\)};

  \node[anchor=east, outer sep=3] at (-1,0) {\(b\)};
  \node[anchor=east] at (-0.5,-1.15) {\(m\)};
\end{tikzpicture}
\end{equation*}
We then take \((S_b^{\cat{M}})_m = \langle (R_b^{\cat{M}})_m \rangle\); that is, \((R_b^{\cat{M}})_m = (S_b^{\cat{M}})_m \id[m]\) with \((S_b^{\cat{M}})_m \in \field\).
Note that \((R_b^{\cat{M}})_m\) is natural in \(m\).
Provided \(b\) is transparent, \((R_b^{\cat{M}})_m\) is a \(\cat{B}\)-module natural transformation.
It follows that \((S_b^{\cat{M}})_m = (S_b^{\cat{M}})_{m'}\) for all \(m\) and \(m'\), so we simply write
\begin{equation*}
  \tilde{S}_{\cat{M},b} \colon \pi_0\cat{Z}(\cat{B}) \times \pi_0\cat{Z}_2(\cat{B}) \to \field \,.
\end{equation*}
for the \textit{framed S-matrix} of \(\cat{Z}(\Mod{\cat{B}})\) (treating \(\cat{M}\) as just an object in \(\cat{Z}(\cat{B})\)).
Looking at how this acts on rows, we find that it is a restriction of the normal S-matrix.
In particular, \(\rank S = \rank \tilde{S}_{\cat{M},b}\).
We have decategorified everything to questions about 1-morphisms.

\begin{theorem}[2.57 of \cite{JFR}]\label{thm:B}
  Let \(\cat{B}\) be a braided fusion 1-category.
  Then the S-matrix
  \begin{equation*}
    S \colon \pi_0 \cat{Z}(\Mod{\cat{B}}) \times \pi_0 \cat{Z}_2(\cat{B}) \to \field
  \end{equation*}
  is of full rank \(|\pi_0\cat{Z}_2(\cat{B})|\).
\end{theorem}
\begin{proof}
  This statement makes sense because Theorem \ref{thm:A} gives that \(S\) is a square matrix.
  Since \(\cat{Z}(\cat{B})\) is a non-degenerate braided fusion 1-category, we have \(\rank \tilde{S}_{\cat{M},b} = |\pi_0\cat{Z}_2(\cat{B})|\).
\end{proof}

\begin{remark}
  Let \(\cat{M}\) be an invertible simple object in a braided fusion 2-category \(\cat{C}\), let \(N\) be another simple object in \(\cat{C}\), and let \(b\) be a simple transparent object in \(\pi_0\Omega\cat{C}\).
  Then \(S_{\cat{M} \boxtimes \cat{N},b} = S_{\cat{M},b} S_{\cat{N},b} \in \field\), so that if \(\cat{M}\) is invertible then \(S\) becomes multiplicative.
  In particular, if every component has an invertible simple object then \(\pi_0\cat{C}\) has a group structure and \(S_{-,b} \colon \pi_0\cat{C} \to \field\) is a group homomorphism.
\end{remark}


\begin{theorem}\label{thm:main}
  If \(\cat{B}\) is a braided fusion 1-category then \(\cat{Z}(\Mod{\cat{B}})\) is non-degenerate.
\end{theorem}
\begin{proof}
  It will be shown in \cite{jrSMAT} that a braided fusion 2-category is non-degenerate if and only if its framed S-matrix is invertible.
  By Theorem \ref{thm:B}, this is the case for \(\cat{Z}(\Mod{\cat{B}})\).
\end{proof}

%% file: talks/4.1/STOAT_notes_talk_4_1.tex
Talk by Cameron Krulewski, notes by Michail Tagaris.

In this talk, we begin to classify the possible forms a slightly-degenerate braided fusion 2-category $\mathcal{B}$ can take. We use some known constraints on the form of $\mathcal{Z}(\Sigma \mathcal{B})$ and a homotopical computation to narrow down on two possibilities for $\mathcal{B}$: $\mathcal{S}$ or $\mathcal{T}$, both of which arise by linearization.
The next talk introduces the invariant by which we can tell these categories apart.

Before we start, we need to recall a fact about the $\Z/2$-cohomology of Eilenberg-MacLane spaces.

\begin{thm}[Serre \cite{serre_cohomologie_nodate}]
    Let $t_n \in H^n(K(\mathbb{Z}/2 , n ); \mathbb{Z}/2)$
be the element classified by the identity map
$\id\colon K(\mathbb{Z}/2 , n) \rightarrow K(\mathbb{Z}/2,n)$.
Then there is an isomorphism \[H^\bullet(K(\mathbb{Z}/2 , n); \mathbb{Z}/2)
\cong \mathbb{Z}/2 [\{\textnormal{Sq}^It_n~|~I~ \textnormal{admissible},~e(I)<n\}]\]
\end{thm}
A sequence is \emph{admissible} if $i_k \geq 2 i_{k+1}$ for each entry in $I$, and the \emph{excess} of a sequence is $e(I) = \sum_{k=1}^n i_k - 2i_{k+1}$.

This theorem is proven using Borel's transgression theorem and the Serre spectral sequence for the path space fibration. See \cite{serre_cohomologie_nodate} for the proof.

\section{Classifying Possible Forms of $\mathcal{B}$}

We have been studying the category $\mathcal{B}$, which is a slightly degenerate braided fusion 1-category (BF1C).
There are some constraints on the category $\mathcal{C}$ defined as the center $\mathcal{C} = \mathcal{Z}(\Sigma \mathcal{B}).$
We have the following facts:
\begin{enumerate}
    \item[$\circled{1}$] There's an equivalence of symmetric fusion 1-categories $\Omega \mathcal{C}\simeq \mathbf{sVec}$.
    
    \item[$\circled{2}$] $|\pi_0 \mathcal{C}|= 2$; i.e.\ $\mathcal{C}$ has two components, which we are going to call the $\id$ component and the magnetic component.
    
    \item[$\circled{3}$] For any $X$ in the non-identity magnetic component, and for $e\in\pi_0(\Omega\mathcal{C})$ the class of the odd line, the braiding must satisfy $br_{e,X} \cdot br_{X,e} = (-1)\id_{\id_X \boxtimes e}$.
\end{enumerate}

Fact $\circled{1}$ was explained in previous talks. Fact $\circled{2}$ follows from the fact that $S$ is non-degenerate, requiring $|\pi_0\mathcal{C}|=|\pi_0\Omega\mathcal{C}|$ (see \cite{JFR} Theorem 2.52), and that $\Omega \mathcal{C}\simeq \textbf{sVec}$, which has two components.
The point of this talk is to prove the following theorem:

\begin{thm}\label{main_thm}
    The properties of $\circled{1}$ + non-degeneracy are satisfied by precisely $2$ categories,  $\mathcal{S} \text{ and }\mathcal{T}$.
\end{thm}

Note that condition \circled{3} \emph{follows from} condition \circled{1}.
To show \circled{3}, we are going to use the following fact about the non-degeneracy of the $S$-matrix. (See \cite{JFR} Lemma 2.56 and Theorem 2.57.)

\begin{fact}
    If all simple objects in $\Omega \mathcal{C}$ are invertible, then ${S}$ induces a group homomorphism
    \begin{align*}
        {S}\colon \pi_0\mathcal{C} &\rightarrow \textnormal{Hom}_{\textnormal{Grp}}(\pi_0 \Omega \mathcal{C} , \Bbbk^{\times}) \\
        A &\mapsto (b \mapsto S_{A,b}).
    \end{align*}
\end{fact}
    
    We know that in our case of interest, the $S$ matrix takes the following form:
    \[
    S =\bordermatrix{&1 & e \cr
    1 & 1 & 1\cr
    \!\!\! \text{\footnotesize mag} \!\!& 1 & ?} .\]
    
    Since, by the fact, we have a group homomorphism from $\Z/2$, the last entry ($?$) must be $\pm 1$. Then because $S$ is non-degenerate, it has to be $-1$. This entry determines the braiding condition in \circled{3}.

Next, we use some results of Johnson-Freyd and Yu to learn more about the form of $\mathcal{C}$.

\begin{thm}[\cite{johnson-freyd_fusion_2021} Theorem B]
    Any fusion $2$-category satisfying \circled{1} has that all simple objects are invertible and each component has exactly 2 simple objects.
\end{thm}

We know that $\mathcal{C}$ has two components, and the previous theorem tells us that any object $M$ in the magnetic component is invertible. Tensoring with this object (taking $-\boxtimes M$) is an invertible functor, so we know that both components of $\mathcal{C}$ ``look the same."
Using this, we can draw a 
picture of $\mathcal{C}$:
\begin{center}
\begin{tikzpicture}[->,auto,node distance=1.5cm, thick,main node/.style=]

  \node[main node] (1) {\text{$I$}};
  \node[main node] (2) [right of=1] {\text{$C$}};
  \node[main node] (3) [right of=2] {\text{}};
  \node[main node] (4) [right of=3] {\text{$M$}};
  \node[main node] (5) [right of=4] {\text{$C\boxtimes M$}};

  \node at (0.75, 1.2) {\text{id component}};
  \node at (5.4, 1.2) {\text{magnetic component}};
  \node (6) at (0.1, -1.4) {\text{\textcolor{blue}{Identity}}};
  \node (7) at (1.7, -1.4) {\text{\textcolor{blue}{Clifford}}};
  
  \draw[thick] (0.75,0) ellipse (2cm and 1cm); 
  \draw[thick] (5.4,0) ellipse (2cm and 1cm);  

  \path[every node/.style={font=\sffamily\small}]
    (1) edge[bend left] node [left] {} (2)
    (2) edge[bend left] node [left] {} (1)
    (4) edge[bend left] node [left] {} (5)
    (5) edge[bend left] node [left] {} (4)
    (6) edge [right,blue] node {} (1)
    (7) edge [right, blue] node {} (2);
\end{tikzpicture}
\end{center}
In the identity component, we have $2$ simple objects, an Identity and a Clifford algebra object. For the magnetic component, we \emph{pick} a simple object that we call $M$. We know that it's invertible, so the other simple object in that component will be $C\boxtimes M$.

Next, we want to argue that $\mathcal{C}$ arises from linearization.
Here, it suffices to prove that $M^2=I$, since this gives a consistent way to multiply things.
Then $\mathcal{C}$ is recoverable from the full monoidal subcategory $\{I,M\}$.

So let's compute $M^2$.
The only two options for $M^2$ are $C$ or $I$ (because from the $S$-matrix you can show that $M^2$ has to be in the trivial component) so it suffices to prove that $M^2 \neq C$.
We can compare $\text{br}_{M,M}$ with $\text{br}_{M^2,M^2}$:

\begin{center}
\begin{tikzpicture}[thick, scale=1.5]

  \coordinate (A1) at (0, 0);  
  \coordinate (A2) at (2, 2);  
  \coordinate (B1) at (0, 2);  
  \coordinate (B2) at (2, 0);  
  \coordinate (C1) at (0.95, 1.05);  
  \coordinate (C2) at (1.05, 0.95);  
  
  \draw (B2) -- (1, 1); 
  \draw (B1) -- (B2); 
  \draw[white, line width=1.5pt] (C1) -- (C2); 
 \draw (A1) -- (A2); 
  
  \node at (-0.3, 0) {\text{M}};
  \node at (-0.3, 2) {\text{M}};
  \node at (2.3, 0) {\text{M}};
  \node at (2.3, 2) {\text{M}};

  \node at (3, 1) {$\text{= scalar} \cdot$};

  \coordinate (C1) at (4, 0);  
  \coordinate (C2) at (4, 2);  
  \coordinate (D1) at (5, 0);  
  \coordinate (D2) at (5, 2);  

  \draw (C1) -- (C2);  
  \draw (D1) -- (D2); 

  \node at (3.7, 0) {\text{M}};
  \node at (3.7, 2) {\text{M}};
  \node at (5.3, 0) {\text{M}};
  \node at (5.3, 2) {\text{M}};

  \node at (-1, 1) {\text{$\text{br}_{M,M}$:}};
  \node at (3, -0.2) {\text{$I$ or $e$}};
    
    \draw[decorate,decoration={snake,amplitude=0.5mm} , ->] (3, 0.0) -- (3.1, 0.8);
  
\end{tikzpicture}

\begin{tikzpicture}[thick, scale=1.5]

  \coordinate (A1) at (0, 2);   
  \coordinate (A2) at (2, 0);   
  \coordinate (B1) at (0.5, 2); 
  \coordinate (B2) at (2.5, 0); 
  
  \coordinate (C1) at (0.2, 0);   
  \coordinate (C2) at (2.2, 2);   
  \coordinate (D1) at (-0.5, 0); 
  \coordinate (D2) at (1.5, 2);  
  
  \draw (A1) -- (A2);
  \draw (B1) -- (B2);

  \draw[white, line width=5.5pt] (C1) -- (C2);
  \draw (C1) -- (C2);
  \draw[white, line width=5.5pt] (D1) -- (D2);
  \draw (D1) -- (D2);

  \node at (-1, 1) {\text{br$_{M^2,M^2}$:}};
  
  \node at (3, 1) {\text{= (1 or $e^4$)}};
  \node at (3.3, 0.75) {\text{\rotatebox{85}{$\,=$}}};
  \node at (3.25, 0.55) {\text{1}};
  
  \draw (4, 0) -- (4, 2);  
  \draw (4.5, 0) -- (4.5, 2);  
  \draw (5, 0) -- (5, 2);  
  \draw (5.5, 0) -- (5.5, 2); 
  
\end{tikzpicture}
\end{center}

We know that the self-braiding of $C$ is nontrivial: $\text{br}_{C,C} \cong e \cdot \id_{C\boxtimes C}$.
This means that $M^2$ cannot be $C$, and so we must have $M^2 = I$.

Next, what are the possible linearizations? The data involved are:
\begin{itemize}
    \item Choice of delooping $B^2 \mathcal{G}$
    \item Twist class $\alpha \in H^5(B^2\mathcal{G} ;\Bbbk^{\times})$.
\end{itemize}
From this information we can construct an extension as the pullback
\begin{center}
    \begin{tikzcd}[ampersand replacement=\&]
            K(\Bbbk^\times,4)\cdot B^2\mathcal{G} \ar[r] \ar[d] \arrow[dr, phantom, "\scalebox{2}{$\lrcorner$}" , very near start] \& * \ar[d] \\
            B^2\mathcal{G} \ar[r,"\alpha"] \& K(\Bbbk^\times,5)
        \end{tikzcd}
    \end{center}
and from there, linearization produces the category $2\textbf{Vec}^{\alpha}[\mathcal{G}]$.

Above, we recalled the form of the $\Z/2$-cohomology of Eilenberg-MacLane spaces. We can use the universal coefficient theorem to change from $\mathbb{Z}/2$ to $\Bbbk^{\times}$ coefficients; we map $\text{Sq}^I {t_n}$ in $H^\bullet(-;\mathbb{Z}/2)$ to $(-1)^{\text{Sq}^I {t_n}}$ in $H^\bullet(-;\Bbbk^{\times})$. Or we can look up the groups in \cite{serre_cohomologie_nodate} or \cite[Ex. 2.4 and Ex. 2.5]{DN21}.

To show Theorem \ref{main_thm}, we need to argue that there are only two options of linearization.\\

\noindent\underline{part 1}:
Note that there is a (non-monoidal) equivalence
\begin{center}
\begin{tikzpicture}
      \node at (0, 0)   {$\mathcal{C} \simeq 2\mathbf{Vec}[K(\mathbb{Z}/2 , 1 )\times K(\mathbb{Z}/2 , 0)].$};
       \node at (-0.3, -1) {$I,e$};
       \node at (1.3, -1) {$I,M$};

        \draw[decorate,decoration={snake,amplitude=0.5mm} , ->] (-0.3, -0.8) -- (-0.3, -0.2);
         \draw[decorate,decoration={snake,amplitude=0.5mm} , ->] (1.3, -0.8)-- (1.3, -0.2);
\end{tikzpicture}
\end{center}
The $K(\Z/2,1)$ piece comes from $\pi_0 \Omega \mathcal{C} = \pi_0 \textbf{sVec}$, while the $K(\Z/2,0)$ piece comes from $\pi_0\mathcal{C}$, which had a trivial and a magnetic component. To determine the monoidal structure, we need to know how these two pieces combine.

We will argue that we don't have a non-trivial extension; i.e., that
\[B^2 \mathcal{G} \simeq K(\mathbb{Z}/2,3) \times K(\mathbb{Z}/2,2).\]

Assume for the purposes of contradiction that $B^2 \mathcal{G}$ is the nontrivial extension, which we will denote by $X$. This extension is classified by the class $\text{Sq}^2 t_2 = t_2^2$. For now, we assume a lemma:
\begin{lemma}[\cite{JFR} Lemma 3.1]\label{loops_lemma}
    Let $\Omega$ be the desuspension map. There is an isomorphism
    $$ \Omega\colon H^6(BX;\Bbbk) \xrightarrow{\cong} H^5(X;\Bbbk). $$
\end{lemma}
This isomorphism ensures that any twist $\alpha\in H^5(X;\Bbbk)$ pulls back to a unique class in $H^6(BX;\Bbbk)$. Such a class guarantees a sylleptic structure on $2\textbf{Vec}^\alpha[\mathcal{G}]$, which would require that all objects have trivial braiding with 1-morphisms from $I$ to $I$. However, we know that this is not the case by condition \circled{3}. $\lightning$.
\\

\underline{part 2}:
Since $B^2\mathcal{G}$ is a product $\left( K(\mathbb{Z}/2,3) \times K(\mathbb{Z}/2,2)\right)$ we can use the K\"unneth formula to compute its cohomology:
\[
H^5(B^2\mathcal{G}; \Bbbk^\times ) \stackrel{\text{K\"unneth}}{\cong} (\mathbb{Z}/2)^3.
\]
We can pick a basis for this $(\Z/2)^3$ using the standard generators for the cohomology of each piece of the product:

\begin{center}
\begin{tikzpicture}
    \node (B) at (0, 0)  {$\{(-1)^{\text{Sq}^2t_3} , (-1)^{t_3t_2} , (-1)^{\text{Sq}^2\text{Sq}^1t_2}\}$};
    \node (A) at (3, 1.5) {encodes $S$-matrix};
    \draw[->] (A) -- (0.3,0.3);
    \node at (-4,-1) {from $ H^5(K(\mathbb{Z}/2 , 3))$};
    \node at (0,-1.5) {from $ H^3(K(\mathbb{Z}/2 , 3)) \otimes H^2(K(\mathbb{Z}/2 , 2)) $};
    \node at (4,-1) {from $ H_5(K(\mathbb{Z}/2 , 2))$};
    \draw[->] (-3,-0.8) -- (-2,-0.3);
    \draw[->] (0,-1.2) -- (0,-0.3);
    \draw[->] (3,-0.8) -- (2,-0.3);
\end{tikzpicture}
\end{center}

Next we want to narrow down the possible $\alpha$'s to just two options. Let's write $\alpha$ as a column vector in $(\Z/2)^3$.

Condition
\circled{1} implies that $  \alpha = \begin{pmatrix}
    1\\
    \cdot \\
    \cdot 
\end{pmatrix}$,
since the restriction $H^5(B^2\mathcal{G}) \rightarrow H^5(K(\mathbb{Z}/2 , 3))$\\
witnesses the non-trivial braiding $\textnormal{br}_{e,e} = -\id_{e,e}$ on $\Omega \mathcal{C} = \mathbf{sVec}$.\\\\
Meanwhile, condition $\circled{3}$ implies that $\alpha = \begin{pmatrix}
    \cdot\\
    1\\
    \cdot
\end{pmatrix}$, because we need the factor of ${t_3t_2}$ to be non-trivial in order to obstruct a sylleptic structure; i.e. to witness the nontrivial braiding
$\textnormal{br}_{e,M} \textnormal{br}_{M,e} = -\id.$

So we are left with $2$ possibilities left for $\alpha$ (depending on the choice of the $t_2$ term):
$$\sigma = \begin{pmatrix}
    1\\1\\0
\end{pmatrix} , \text{ and }  \tau = \begin{pmatrix}
    1\\1\\1
\end{pmatrix}.$$
This concludes our sketch proof of Theorem \ref{main_thm}.

\section{Comparing $\mathcal{S}$ and $\mathcal{T}$}

We will name the two categories produced by linearization for these two classes:
$\mathcal{S} \coloneqq 2\mathbf{Vec}t^\sigma [\mathcal{G}]$ and 
$\mathcal{T} \coloneqq 2\mathbf{Vec}t^\tau [\mathcal{G}].$
Observe that these two twists differ in their restrictions over the $K(\Z/2,2)$ factor:
$\sigma|_{K(\mathbb{Z}/2 , 2)} \simeq \mathrm{trivial}$, while $\tau|_{K({\mathbb{Z}/2 , 2})} =: \tau'$ is nontrivial.
We see that $\mathcal{S}$ and $\mathcal{T}$ admit the following subcategories:
\begin{itemize}
    \item $2\textbf{Vec}[\Z/2] \subset \mathcal{S}$
    \item $2\textbf{Vec}^{\tau'}[\Z/2] \subset \mathcal{T}$.
\end{itemize}

Each of these subcategories maps to $\Sigma \mathbf{sVec}$, since $\Sigma \mathbf{sVec} = 2\mathbf{Vec}^\omega [B\mathbb{Z}/2]$ for $\omega = (-1)^{\text{Sq}^2t_3}$. To define these maps, it suffices to give a map $f\colon K(\Z/2,2)\to K(\Z/2,3)$ that pulls back $\omega$ to the respective twist $\alpha$; i.e. $f^* \omega = \alpha$.
The trivial twist corresponds to the trivial map, and the nontrivial twist to the nontrivial map:
\begin{itemize}
    \item for $\sigma$, $f=0$
    \item for $\tau$, $f = \text{Sq}^1$.
\end{itemize}
Each map $f$ gives rise to a functor:
$$
f\colon K(\Z/2,2)\to K(\Z/2,3) \rightsquigarrow F\colon 2\mathbf{Vec}^\alpha[\mathbb{Z}/2]\rightarrow2\mathbf{Vec}^\omega[B\mathbb{Z}/2].
$$

Note that both $F$'s map $M\mapsto I$. The difference is in the monoidality isomorphisms:
\begin{center}
\begin{tikzpicture}[text width=8.8cm]
    \node at (0,0) {$I = I \boxtimes I = F(I) = F(M) \boxtimes F(M) \textcolor{blue}{\rightarrow} F(M \boxtimes M) = I$};
    \node at (6,-0.7) {\textcolor{blue}{ for $\mathcal{S}$, take the trivial map}};
    \node at (6,-1.4) {\textcolor{blue}{ for $\mathcal{T}$, take the map corresponding to $e\colon I^\vee \cong I$}};
\end{tikzpicture}
\end{center}

\section{Proof of the Lemma}

We now want to prove the lemma from earlier.
\begin{lemma}
For the non-trivial extension $X$
\[
K(\mathbb{Z}/2,3) \rightarrow X \rightarrow K(\mathbb{Z}/2 , 2)
\]
classified by $\text{Sq}^2t_2 = t^2_2$, the desuspension map
\[
\Omega\colon H^6(BX;\Bbbk^\times ) \rightarrow H^5(X;\Bbbk^\times )
\]
is an isomorphism.
\end{lemma}

It suffices to prove that we have a short exact sequence
\[
0 \rightarrow H^5(K(\mathbb{Z}/2,2);\Bbbk^\times) \rightarrow H^5(X;\Bbbk^\times) \rightarrow H^5(K(\mathbb{Z}/2,3);\Bbbk^\times) \rightarrow 0 
\]
and similarly for $BX$. Then, we can use the isomorphisms for each Eilenberg-MacLane space and conclude the desired statement from the Five Lemma.

To show this, we will use the Serre spectral sequence for computing the cohomology of $X$ from that of its parts.
It has the following type signature:
\[
E_2^{i,j} = H^i(K(\mathbb{Z}/2,2), H^j(K(\mathbb{Z}/2,3);\Bbbk^\times)) \Rightarrow H^{i+j}(X;\Bbbk^\times).
\]
This means that we fill in the $E_2$ page of the spectral sequence using the cohomology groups of $K(\mathbb{Z}/2,2)$ with coefficients in the cohomology of $K(\mathbb{Z}/2,3)$ with $\Bbbk^\times$ coefficients. After ``turning" a few pages of the spectral sequence by computing differentials, it will eventually stabilize to an $E_\infty$ page, whose entries will tell us the cohomology of $X$ (up to some possible extension questions).
We have the following $E_2$ page, using the cohomology groups we know from \cite{serre_cohomologie_nodate} or \cite{DN21}:

\begin{center}
\begin{tikzpicture}[scale=1.2]

  \draw[thick,->] (-0.5,-0.5) -- (7,-0.5) node[anchor=north] {$i$};
  \draw[thick,->] (-0.5,-0.5) -- (-0.5,6) node[anchor=east] {$j$};
  
  \node at (0,0) {$\Bbbk^\times$};
  \node (A) at (0,3) {$\mathbb{Z}/2$};
  \node at (2,0) {$\mathbb{Z}/2$};
  \node (C) at (2,3) {$\mathbb{Z}/2$};
  \node (E) at (3,3) {$\mathbb{Z}/2$};
  \node (B) at (6,0) {$\mathbb{Z}/2$};
  \node (D) at (5,0) {$\mathbb{Z}/2$};
  \foreach \x in {0,...,5} {
  \node at (-1,\x) {$\x$};
  }
  \foreach \x in {0,...,6} {
  \node at (\x , -1) {$\x$};
  }
  \node (F) at (0,5) {$\mathbb{Z}/2$};
  \node (G) at (4,0) {$\mathbb{Z}/4$};

    \foreach \x in {0,...,4} {
  \node at (\x,2) {$0$};
  }
  \foreach \x in {0,...,5} {
  \node at (\x,1) {$0$};
  }
  \foreach \x in {0,...,2} {
  \node at (\x,4) {$0$};
  }
  \node at (1,0) {$0$};
  \node at (3,0) {$0$};
  \node at (1,3) {$0$};
  \node at (1,5) {$0$};
\end{tikzpicture}
\end{center}

Differentials in the Serre spectral sequence have bidegree $d_n\colon E_n^{i,j}\to E_n^{i+n,j-n+1}$. For degree reasons, there are no possible $d_2$ differentials. But there is a possible $d_3$:

\begin{center}
\begin{tikzpicture}[scale=1.2]
  \draw[thick,->] (-0.5,-0.5) -- (7,-0.5) node[anchor=north] {$i$};
  \draw[thick,->] (-0.5,-0.5) -- (-0.5,6) node[anchor=east] {$j$};
  
  \node at (0,0) {$\Bbbk^\times$};
  \node (A) at (0,3) {$\mathbb{Z}/2$};
  \node at (2,0) {$\mathbb{Z}/2$};
  \node (C) at (2,3) {$\mathbb{Z}/2$};
  \node (E) at (3,3) {$\mathbb{Z}/2$};
  \node (B) at (6,0) {$\mathbb{Z}/2$};
  \node (D) at (5,0) {$\mathbb{Z}/2$};
  \foreach \x in {0,...,5} {
  \node at (-1,\x) {$\x$};
  }
  \foreach \x in {0,...,6} {
  \node at (\x , -1) {$\x$};
  }
  \node (F) at (0,5) {$\mathbb{Z}/2$};
  \node (G) at (4,0) {$\mathbb{Z}/4$};

  \draw[->, thick, blue] (F) -- (E); 
  \node at (3,4.5) [purple] {$d_3=0$ or not?};

    \foreach \x in {0,...,4} {
  \node at (\x,2) {$0$};
  }
  \foreach \x in {0,...,5} {
  \node at (\x,1) {$0$};
  }
  \foreach \x in {0,...,2} {
  \node at (\x,4) {$0$};
  }
  \node at (1,0) {$0$};
  \node at (3,0) {$0$};
  \node at (1,3) {$0$};
  \node at (1,5) {$0$};
\end{tikzpicture}
\end{center}

While this differential is possible, it doesn't actually happen. This is because we know the $\Z/2$ in degree $(0,5)$ has to survive (i.e.\ not be eliminated by the $d_3$): the edge map $H^5(X;\Bbbk^\times) \to E_2^{0,5}=H^5(K(\Z/2,2),\Bbbk^\times)$ is nontrivial. We witness this nontriviality with the category $\Sigma\textbf{sVec}$ discussed in the previous section, which maps to the (nontrivial) extension $X$ under the quotient:

\begin{center}
\begin{tikzcd}[column sep=normal, row sep=normal]
  H^5(X; \Bbbk^\times) \arrow[r] & 
  H^5(K(\mathbb{Z}/2,3); \Bbbk^\times) \\
  (\Sigma \mathbf{sVect})^\times \arrow[r, mapsto] \arrow[u, "\alpha" , draw=none] &
  \Omega(\Sigma \mathbf{sVect})^\times \\
  {K}(\Bbbk^\times,4) \cdot \textcolor{purple}{\underbrace{\textcolor{black}{{K}(\mathbb{Z}/2,3) \cdot {K}(\mathbb{Z}/2,2)}}_{X}}\arrow[u , "\text{\rotatebox[]{270}{$\simeq$}}" , draw = none ] \arrow[r, mapsto] &
  {K}(\Bbbk^\times,4) \cdot {K}(\mathbb{Z}/2,3)  \\
\end{tikzcd}
\begin{tikzpicture}[overlay]
    \node at (-1,0.2) {\rotatebox[]{315}{$\simeq$}};
     \node at (-0.45,-0.05) {\footnotesize $\mathbf{sVect}^\times$};
  \node at (-7.75, -2.35) [blue] {\text{nontrivial}};
  \node at (-7.75, -2.7) [blue] {\text{$\textnormal{br}_{c,c} = e \cdot \id$}};
   \draw[blue, decorate,decoration={snake,amplitude=0.5mm} , ->] (-7.75, -2.25) -- (-6.3, -0.76);
\end{tikzpicture}
\end{center}

\vspace{1cm}

In summary, we see explicitly that the edge map is an isomorphism $\Z/2\cong\Z/2$, and from this we know the nontrivial class in $E_2^{0,5}$ can't get killed by the $d_3$.

Next, we consider $d_4$ differentials. There is a differential (0,3) to (4,0) corresponding to the extension class of $X$: the image of this first transgressive differential is the extension class $t_2^2\in H^4(K(\Z/2,2);\Z/2)$.
We have $d_4((-1)^{t_3}) = (-1)^{t_2^2}$, hitting the order-two element in (4,0).

This behavior propagates by the Leibniz rule: the $d_4$ from (2,3) to (6,0) is given by
$$ d_4((-1)^{t_3t_2}) = d_4((-1)^{t_3}) (-1)^{t_2} + (-1)^{t_3}d_4((-1)^{t_2}) = (-1)^{t_2^3} + 0. $$

\begin{tikzpicture}[scale=1.2]
  \draw[thick,->] (-0.5,-0.5) -- (7,-0.5) node[anchor=north] {$i$};
  \draw[thick,->] (-0.5,-0.5) -- (-0.5,6) node[anchor=east] {$j$};
  
  \node at (0,0) {$\Bbbk^\times$};
  \node (A) at (0,3) {$\mathbb{Z}/2$};
  \node at (2,0) {$\mathbb{Z}/2$};
  \node (C) at (2,3) {$\mathbb{Z}/2$};
  \node (E) at (3,3) {$\mathbb{Z}/2$};
  \node (B) at (6,0) {$\mathbb{Z}/2$};
  \node (D) at (5,0) {$\mathbb{Z}/2$};
  \foreach \x in {0,...,5} {
  \node at (-1,\x) {$\x$};
  }
  \foreach \x in {0,...,6} {
  \node at (\x , -1) {$\x$};
  }
  \node (F) at (0,5) {$\mathbb{Z}/2$};
  \node (G) at (4,0) {$\mathbb{Z}/4$};

  \draw[->, thick, blue] (A) -- (G);
  \draw[->, thick, blue] (C) -- (B); 

  \node at (2.55,1.7) [blue] {differential};

  \node[below, purple] (I) at (4.4,-1.2) {\footnotesize order $2$};
  \node[purple] at (4.4,-1.8) {\footnotesize $(-1)^{t_2^2}$};
  \node[below] (H) at (-2,3) [purple] {\footnotesize $(-1)^{t_3}$};
  \draw[->, thick, purple] (H) to [out=350,in=210](A); 
  \draw[->, thick, purple] (I) to [out=180,in=210](G);

    \foreach \x in {0,...,4} {
  \node at (\x,2) {$0$};
  }
  \foreach \x in {0,...,5} {
  \node at (\x,1) {$0$};
  }
  \foreach \x in {0,...,2} {
  \node at (\x,4) {$0$};
  }
  \node at (1,0) {$0$};
  \node at (3,0) {$0$};
  \node at (1,3) {$0$};
  \node at (1,5) {$0$};
\end{tikzpicture}

There are no $d_5$'s or higher differentials for degree reasons, so the spectral sequence collapses on page 5. (We really did have to check the $d_4$'s though; if they had been trivial then there could have been a $d_6$ from (0,5) to (6,0).)
Let's draw the $E_5=E_\infty$ page. We are interested in the classes in total degree five since we are computing $H^5(X;\Bbbk^\times)$:

\begin{center}
\begin{tikzpicture}[scale=1.2]
  \draw[thick,->] (-0.5,-0.5) -- (7,-0.5) node[anchor=north] {$i$};
  \draw[thick,->] (-0.5,-0.5) -- (-0.5,6) node[anchor=east] {$j$};
  
  \node at (0,0) {$\Bbbk^\times$};
  \node (A) at (0,3) {$0$};
  \node at (2,0) {$\mathbb{Z}/2$};
  \node (C) at (2,3) {$0$};
  \node (E) at (3,3) {$\mathbb{Z}/2$};
  \node (B) at (6,0) {$0$};
  \node[draw, rectangle] (D) at (5,0) {$\mathbb{Z}/2$};
  \foreach \x in {0,...,5} {
  \node at (-1,\x) {$\x$};
  }
  \foreach \x in {0,...,6} {
  \node at (\x , -1) {$\x$};
  }
  \node[draw, rectangle] (F) at (0,5) {$\mathbb{Z}/2$};
  \node (G) at (4,0) {$\mathbb{Z}/2$};

    \foreach \x in {0,...,4} {
  \node at (\x,2) {$0$};
  }
  \foreach \x in {0,...,5} {
  \node at (\x,1) {$0$};
  }
  \foreach \x in {0,...,2} {
  \node at (\x,4) {$0$};
  }
  \node at (1,0) {$0$};
  \node at (3,0) {$0$};
  \node at (1,3) {$0$};
  \node at (1,5) {$0$};
\end{tikzpicture}
\end{center}

We reconstruct $H^5(X;\Bbbk^\times)$ from the groups in total degree five, which make up its filtered pieces. There is a short exact sequence
$$
0 \to
  E_{\infty}^{5,0} \to
  H^5(X; \Bbbk^\times) \to 
  E_{\infty}^{0,5} \to 
  0.
$$
We argued above that the edge map from $H^5(X;\Bbbk^\times)$ to $E_2^{5,0}$ was nontrivial, and we know the other edge map is also nontrivial since the entry (5,0) survives the sequence.
We have an identification of the entry $E^{5,0}_\infty$ with the entry $E^{5,0}_2$ since there were no differentials modifying that group, and thus we can identify $E^{5,0}_\infty$ with the cohomology group $H^5(K(\Z/2,2),\Bbbk^\times)$. A similar statement is true for $E^{0,5}_\infty$.

Finally, an analogous argument for $BX$ gives a similar short exact sequence. The desuspension maps $\Omega\colon H^i(K(G,n))\to H^{i-1}(K(G,n-1))$ are always isomorphisms, so we can apply the Five Lemma to conclude the lemma:
\begin{center}
\begin{tikzcd}
  0 \arrow[r] & 
  E_{\infty}^{5,0} \arrow[r] & 
  H^5(X; \Bbbk^\times) \arrow[r] & 
  E_{\infty}^{0,5} \arrow[r] & 
  0 \\
  & 
  E_2^{5,0} \arrow[d, "\text{\rotatebox{90}{$\,=$}}" description, draw=none]
  \arrow[u, "\text{\rotatebox{90}{$\,=$}}" description, draw=none] & &  E_2^{0,5} 
  \arrow[u, "\text{\rotatebox{90}{$\,=$}}" description, draw=none]
  \\
 & H^5({K}(\mathbb{Z}/2,2); \Bbbk^\times) 
 & 
  & 
  H^5({K}(\mathbb{Z}/3,2); \Bbbk^\times) \arrow[u, "\text{\rotatebox{90}{$\,=$}}" description, draw=none]\\
  0 \arrow[r] & 
  H^6({K}(\mathbb{Z}/2,3)) \arrow[u  , "\simeq", purple] \arrow[r] & 
  H^6(BX) \arrow[r] & 
  H^6({K}(\mathbb{Z}/2,4)) \arrow[u,"\simeq", purple] \arrow[r] & 
  0
\end{tikzcd}

\begin{tikzpicture}[overlay]
  \draw[->, thick , purple] (0,0.6) -- (0,4) node[midway, right] {\footnotesize \textcolor{purple}{$\cong$ }}; 
  \node at (0.7 , 2.7) {\footnotesize \textcolor{purple}{$5$ lemma}};
\end{tikzpicture}
\end{center}

%% file: talks/4.2/main.tex
Talk by Nivedita, notes by Chetan Vuppulury.

\section{The 1-categorical Quantum Dimension}
    Let $\mathcal E$ be a braided monoidal $1$-category. For any right-dualisable object $x\in\mathcal E$, define
    \begin{center}
        \tikzset{every picture/.style={line width=0.75pt}}
        \begin{tikzpicture}[x=0.75pt,y=0.75pt,yscale=-1,xscale=1]
            \draw    (232.6,237.07) .. controls (232.63,262.86) and (256.59,262.65) .. (256.56,237.07) ;
            \draw    (232.6,205.5) .. controls (232.63,178.23) and (256.59,178.87) .. (256.56,205.5) ;
            \draw [shift={(249.33,186.82)}, rotate = 194.87] [fill={rgb, 255:red, 0; green, 0; blue, 0 }  ][line width=0.08]  [draw opacity=0] (8.93,-4.29) -- (0,0) -- (8.93,4.29) -- cycle    ;
            \draw    (232.6,204.85) .. controls (232.63,217.32) and (256.59,224.63) .. (256.56,237.07) ;
            \draw    (232.6,237.07) .. controls (232.63,230.64) and (239.23,224.91) .. (242.21,222.48) ;
            \draw    (245.81,219.26) .. controls (248.92,216.06) and (256.56,210.01) .. (256.56,205.5) ;
            \draw (257.56,237.07) node [anchor=north west][inner sep=0.75pt]  [font=\scriptsize]  {$x$};
            \draw (70,212) node [anchor=north west][inner sep=0.75pt]    {$\eta\left(x\right)=\mathrm{ev}_\circ\beta_{x^*,x}\circ\mathrm{coev}_x=$};
        \end{tikzpicture}
    \end{center}
    Using the naturality of the braiding and the fact that there is a unique morphism between any two choices of right-dual data, we see that $\eta\left(x\right)$ does not depend on the choice of duality data. Further, it only depends on the isomorphism class of $x$, and thus descends to a map
    \begin{equation*}
        \eta\colon\pi_0\mathcal E^\mathrm{rd}\to\Omega\mathcal E
    \end{equation*}
    where $\mathcal E^\mathrm{rd}$ is the full subcategory of $\mathcal E$ on right-dualisable objects. From the diagram bellow, we can also see that when $x$ has a left dual, we have $\eta\left({}^*x\right)=\eta\left(x\right)$.
    \begin{center}
        \tikzset{every picture/.style={line width=0.75pt}}
        \begin{tikzpicture}[x=0.75pt,y=0.75pt,yscale=-1,xscale=1]
            \draw    (254.6,237.07) .. controls (254.63,262.86) and (278.59,262.65) .. (278.56,237.07) ;
            \draw    (254.6,205.5) .. controls (254.63,178.23) and (278.59,178.87) .. (278.56,205.5) ;
            \draw [shift={(271.33,186.82)}, rotate = 194.87] [fill={rgb, 255:red, 0; green, 0; blue, 0 }  ][line width=0.08]  [draw opacity=0] (8.93,-4.29) -- (0,0) -- (8.93,4.29) -- cycle    ;
            \draw    (254.6,204.85) .. controls (254.63,217.32) and (278.59,224.63) .. (278.56,237.07) ;                \draw    (254.6,237.07) .. controls (254.63,230.64) and (261.23,224.91) .. (264.21,222.48) ;
            \draw    (267.81,219.26) .. controls (270.92,216.06) and (278.56,210.01) .. (278.56,205.5) ;
            \draw    (327.6,237.07) .. controls (327.63,262.86) and (376.59,262.65) .. (376.56,237.07) ;
            \draw    (302.6,205.5) .. controls (302.33,179) and (351.59,178.87) .. (351.56,205.5) ;
            \draw [shift={(332.12,185.96)}, rotate = 183.66] [fill={rgb, 255:red, 0; green, 0; blue, 0 }  ][line width=0.08]  [draw opacity=0] (8.93,-4.29) -- (0,0) -- (8.93,4.29) -- cycle    ;
            \draw    (327.6,204.85) .. controls (327.63,217.32) and (351.59,224.63) .. (351.56,237.07) ;
            \draw    (327.6,237.07) .. controls (327.63,230.64) and (334.23,224.91) .. (337.21,222.48) ;
            \draw    (340.81,219.26) .. controls (343.92,216.06) and (351.56,210.01) .. (351.56,205.5) ;
            \draw    (302.6,205.5) .. controls (302.53,229.51) and (325.33,185) .. (327.6,204.85) ;
            \draw    (351.56,237.07) .. controls (351.49,261.08) and (376.33,215) .. (376.56,237.43) ;
            \draw (279.56,237.07) node [anchor=north west][inner sep=0.75pt]  [font=\scriptsize]  {$x$};
            \draw (348,221) node [anchor=north west][inner sep=0.75pt]  [font=\scriptsize]  {$x$};
            \draw (285,217) node [anchor=north west][inner sep=0.75pt]    {$=$};
        \end{tikzpicture}
    \end{center}
    \begin{remark}
        In the above diagrams, framing data has been suppressed but it is important for distinguishing between left and right duals. However, in the case of a symmetric monoidal category this becomes trivial.
    \end{remark}
    \section{A review of duality in 2-categories}
    Let $\mathcal C$ be a 2-category. We represent objects $X$, $1$-morphisms $f\colon X\to Y$, and $2$-morphisms $\alpha\colon f\to g$ by
    \begin{center}
        \tikzset{every picture/.style={line width=0.75pt}} \begin{tikzpicture}[x=0.75pt,y=0.75pt,yscale=-1,xscale=1]
            \draw  [fill={rgb, 255:red, 65; green, 117; blue, 5 }  ,fill opacity=1 ] (100,110.09) -- (130.4,110.09) -- (130.4,170.09) -- (100,170.09) -- cycle ;
            \draw  [fill={rgb, 255:red, 208; green, 2; blue, 27 }  ,fill opacity=1 ] (130,110.09) -- (160.4,110.09) -- (160.4,170.09) -- (130,170.09) -- cycle ;
            \draw    (130.4,110.09) -- (130,170.09) ;
            \draw [shift={(130.24,133.59)}, rotate = 90.38] [fill={rgb, 255:red, 0; green, 0; blue, 0 }  ][line width=0.08]  [draw opacity=0] (8.93,-4.29) -- (0,0) -- (8.93,4.29) -- cycle    ;
            \draw  [fill={rgb, 255:red, 65; green, 117; blue, 5 }  ,fill opacity=1 ] (59,110.09) -- (89.4,110.09) -- (89.4,170.09) -- (59,170.09) -- cycle ;
            \draw  [fill={rgb, 255:red, 65; green, 117; blue, 5 }  ,fill opacity=1 ] (170,110.09) -- (200.4,110.09) -- (200.4,170.09) -- (170,170.09) -- cycle ;
            \draw  [fill={rgb, 255:red, 208; green, 2; blue, 27 }  ,fill opacity=1 ] (200,110.09) -- (230.4,110.09) -- (230.4,170.09) -- (200,170.09) -- cycle ;
            \draw    (200,139.33) -- (200,170.09) ;
            \draw [shift={(200,148.21)}, rotate = 90] [fill={rgb, 255:red, 0; green, 0; blue, 0 }  ][line width=0.08]  [draw opacity=0] (8.93,-4.29) -- (0,0) -- (8.93,4.29) -- cycle    ;
            \draw    (200,110.09) -- (200,140.85) ;
            \draw [shift={(200,118.97)}, rotate = 90] [fill={rgb, 255:red, 0; green, 0; blue, 0 }  ][line width=0.08]  [draw opacity=0] (8.93,-4.29) -- (0,0) -- (8.93,4.29) -- cycle    ;
            \draw  [fill={rgb, 255:red, 255; green, 255; blue, 255 }  ,fill opacity=1 ] (197.2,140.09) .. controls (197.2,138.42) and (198.55,137.07) .. (200.22,137.07) .. controls (201.89,137.07) and (203.25,138.42) .. (203.25,140.09) .. controls (203.25,141.76) and (201.89,143.12) .. (200.22,143.12) .. controls (198.55,143.12) and (197.2,141.76) .. (197.2,140.09) -- cycle ;
            \draw (133,134) node [anchor=north west][inner sep=0.75pt]  [font=\scriptsize]  {$\textcolor[rgb]{1,1,1}{f}$};
            \draw (99,155) node [anchor=north west][inner sep=0.75pt]  [font=\scriptsize,color={rgb, 255:red, 255; green, 255; blue, 255 }  ,opacity=1 ]  {$X$};
            \draw (147,155) node [anchor=north west][inner sep=0.75pt]  [font=\scriptsize,color={rgb, 255:red, 255; green, 255; blue, 255 }  ,opacity=1 ]  {$Y$};
            \draw (57,155) node [anchor=north west][inner sep=0.75pt]  [font=\scriptsize,color={rgb, 255:red, 255; green, 255; blue, 255 }  ,opacity=1 ]  {$X$};
            \draw (169,155) node [anchor=north west][inner sep=0.75pt]  [font=\scriptsize,color={rgb, 255:red, 255; green, 255; blue, 255 }  ,opacity=1 ]  {$X$};
            \draw (217,155) node [anchor=north west][inner sep=0.75pt]  [font=\scriptsize,color={rgb, 255:red, 255; green, 255; blue, 255 }  ,opacity=1 ]  {$Y$};
            \draw (203,148) node [anchor=north west][inner sep=0.75pt]  [font=\scriptsize]  {$\textcolor[rgb]{1,1,1}{f}$};
            \draw (203,120) node [anchor=north west][inner sep=0.75pt]  [font=\scriptsize]  {$\textcolor[rgb]{1,1,1}{g}$};
            \draw (187,136) node [anchor=north west][inner sep=0.75pt]  [font=\scriptsize]  {$\textcolor[rgb]{1,1,1}{\alpha }$};
        \end{tikzpicture}
    \end{center}
    \begin{definition}
        A right adjoint to a $1$-morphism $f\colon X\to Y$ in $\mathcal C$ is the data of a $1$-morphism $f^*\colon Y\to X$ and two $2$-morphisms $\mathrm{ev}_f$ and $\mathrm{coev}_f$. Diagrammatically, we represent them as
        \begin{center}
            \tikzset{every picture/.style={line width=0.75pt}}
            \begin{tikzpicture}[x=0.75pt,y=0.75pt,yscale=-1,xscale=1]
                \draw  [fill={rgb, 255:red, 208; green, 2; blue, 27 }  ,fill opacity=1 ] (100,191.09) -- (130.4,191.09) -- (130.4,251.09) -- (100,251.09) -- cycle ;
                \draw  [fill={rgb, 255:red, 65; green, 117; blue, 5 }  ,fill opacity=1 ] (130,191.09) -- (160.4,191.09) -- (160.4,251.09) -- (130,251.09) -- cycle ;
                \draw    (130.4,251.09) -- (130.4,191.09) ;
                \draw [shift={(130.4,227.59)}, rotate = 270] [fill={rgb, 255:red, 0; green, 0; blue, 0 }  ][line width=0.08]  [draw opacity=0] (8.93,-4.29) -- (0,0) -- (8.93,4.29) -- cycle    ;
                \draw  [fill={rgb, 255:red, 208; green, 2; blue, 27 }  ,fill opacity=1 ] (180.38,190.13) -- (240.09,190.13) -- (240.09,250.21) -- (180.38,250.21) -- cycle ;
                \draw [fill={rgb, 255:red, 65; green, 117; blue, 5 }  ,fill opacity=1 ]   (201.07,249.39) .. controls (201.07,219.85) and (222.87,218.86) .. (222.55,249.72) ;
                \draw [shift={(206.62,229.59)}, rotate = 331.84] [fill={rgb, 255:red, 0; green, 0; blue, 0 }  ][line width=0.08]  [draw opacity=0] (8.93,-4.29) -- (0,0) -- (8.93,4.29) -- cycle    ;
                \draw  [fill={rgb, 255:red, 65; green, 117; blue, 5 }  ,fill opacity=1 ] (322.09,250.02) -- (261.14,250.04) -- (261.13,190.14) -- (322.08,190.13) -- cycle ;
                \draw [fill={rgb, 255:red, 208; green, 2; blue, 27 }  ,fill opacity=1 ]   (301.92,190.96) .. controls (301.93,220.41) and (279.68,221.4) .. (280,190.63) ;
                \draw [shift={(286.03,211.41)}, rotate = 18.63] [fill={rgb, 255:red, 0; green, 0; blue, 0 }  ][line width=0.08]  [draw opacity=0] (8.93,-4.29) -- (0,0) -- (8.93,4.29) -- cycle    ;
                \draw (133,213) node [anchor=north west][inner sep=0.75pt]  [font=\scriptsize]  {$\textcolor[rgb]{1,1,1}{f^{*}}$};
                \draw (99,236) node [anchor=north west][inner sep=0.75pt]  [font=\scriptsize,color={rgb, 255:red, 255; green, 255; blue, 255 }  ,opacity=1 ]  {$Y$};
                \draw (147,236) node [anchor=north west][inner sep=0.75pt]  [font=\scriptsize,color={rgb, 255:red, 255; green, 255; blue, 255 }  ,opacity=1 ]  {$X$};
                \draw (179.6,235.11) node [anchor=north west][inner sep=0.75pt]  [font=\scriptsize,color={rgb, 255:red, 255; green, 255; blue, 255 }  ,opacity=1 ]  {$Y$};
                \draw (201.4,235.11) node [anchor=north west][inner sep=0.75pt]  [font=\scriptsize,color={rgb, 255:red, 255; green, 255; blue, 255 }  ,opacity=1 ]  {$X$};
                \draw (211.71,215.16) node [anchor=north west][inner sep=0.75pt]  [font=\tiny]  {$\textcolor[rgb]{1,1,1}{f}$};
                \draw (290.61,211.32) node [anchor=north west][inner sep=0.75pt]  [font=\tiny]  {$\textcolor[rgb]{1,1,1}{f}$};
                \draw (259.77,235.13) node [anchor=north west][inner sep=0.75pt]  [font=\scriptsize,color={rgb, 255:red, 255; green, 255; blue, 255 }  ,opacity=1 ]  {$X$};
                \draw (279.73,190.49) node [anchor=north west][inner sep=0.75pt]  [font=\scriptsize,color={rgb, 255:red, 255; green, 255; blue, 255 }  ,opacity=1 ]  {$Y$};
            \end{tikzpicture}
        \end{center}
        We require them to satisfy two conditions represented by
        \begin{center}
            \tikzset{every picture/.style={line width=0.75pt}}
            \begin{tikzpicture}[x=0.75pt,y=0.75pt,yscale=-1,xscale=1]
                \draw  [fill={rgb, 255:red, 208; green, 2; blue, 27 }  ,fill opacity=1 ] (92.38,289.13) -- (152.09,289.13) -- (152.09,349.21) -- (92.38,349.21) -- cycle ;
                \draw [draw opacity=0][fill={rgb, 255:red, 65; green, 117; blue, 5 }  ,fill opacity=1 ]   (113.07,348.39) .. controls (113.07,334.15) and (117.16,305.92) .. (123.35,325.94) .. controls (129.53,345.96) and (136.33,338) .. (136.67,337.67) ;
                \draw [shift={(114.38,331.71)}, rotate = 278.91] [line width=0.08]  [draw opacity=0] (5.36,-2.57) -- (0,0) -- (5.36,2.57) -- (3.56,0) -- cycle    ;
                \draw [shift={(126.33,333.55)}, rotate = 61.58] [line width=0.08]  [draw opacity=0] (5.36,-2.57) -- (0,0) -- (5.36,2.57) -- (3.56,0) -- cycle    ;
                \draw [draw opacity=0][fill={rgb, 255:red, 65; green, 117; blue, 5 }  ,fill opacity=1 ]   (135.67,337.67) .. controls (140.33,338) and (138.33,284) .. (139.33,289) .. controls (140.33,294) and (152.33,289) .. (152.09,349.21) ;
                \draw  [draw opacity=0][fill={rgb, 255:red, 65; green, 117; blue, 5 }  ,fill opacity=1 ] (139.33,289) -- (151.33,289) -- (151.33,349.33) -- (139.33,349.33) -- cycle ;
                \draw [draw opacity=0]   (136.67,337.67) .. controls (140,327) and (139,319) .. (139,319) ;
                \draw  [draw opacity=0][fill={rgb, 255:red, 65; green, 117; blue, 5 }  ,fill opacity=1 ] (113.07,338.33) -- (147.07,338.33) -- (147.07,349.33) -- (113.07,349.33) -- cycle ;
                \draw [fill={rgb, 255:red, 65; green, 117; blue, 5 }  ,fill opacity=1 ]   (113.07,349.33) .. controls (112.21,349.58) and (116.33,290) .. (128.33,337) ;
                \draw    (128.33,337) .. controls (139,344) and (138,327) .. (138.04,332.08) ;
                \draw    (139.33,289) -- (138.04,332.08) ;
                \draw [shift={(138.54,315.54)}, rotate = 271.72] [fill={rgb, 255:red, 0; green, 0; blue, 0 }  ][line width=0.08]  [draw opacity=0] (8.93,-4.29) -- (0,0) -- (8.93,4.29) -- cycle    ;
                \draw   (92.38,289.13) -- (152.09,289.13) -- (152.09,349.21) -- (92.38,349.21) -- cycle ;
                \draw  [fill={rgb, 255:red, 208; green, 2; blue, 27 }  ,fill opacity=1 ] (175,290.09) -- (205.4,290.09) -- (205.4,350.09) -- (175,350.09) -- cycle ;
                \draw  [fill={rgb, 255:red, 65; green, 117; blue, 5 }  ,fill opacity=1 ] (205,290.09) -- (235.4,290.09) -- (235.4,350.09) -- (205,350.09) -- cycle ;
                \draw    (205.4,350.09) -- (205.4,290.09) ;
                \draw [shift={(205.4,324.19)}, rotate = 270] [fill={rgb, 255:red, 0; green, 0; blue, 0 }  ][line width=0.08]  [draw opacity=0] (7.14,-3.43) -- (0,0) -- (7.14,3.43) -- (4.74,0) -- cycle    ;
                \draw (91.6,334.11) node [anchor=north west][inner sep=0.75pt]  [font=\scriptsize,color={rgb, 255:red, 255; green, 255; blue, 255 }  ,opacity=1 ]  {$Y$};
                \draw (139.04,335.08) node [anchor=north west][inner sep=0.75pt]  [font=\scriptsize,color={rgb, 255:red, 255; green, 255; blue, 255 }  ,opacity=1 ]  {$X$};
                \draw (208,312) node [anchor=north west][inner sep=0.75pt]  [font=\scriptsize]  {$\textcolor[rgb]{1,1,1}{f}\textcolor[rgb]{1,1,1}{^{*}}$};
                \draw (174,335) node [anchor=north west][inner sep=0.75pt]  [font=\scriptsize,color={rgb, 255:red, 255; green, 255; blue, 255 }  ,opacity=1 ]  {$Y$};
                \draw (222,335) node [anchor=north west][inner sep=0.75pt]  [font=\scriptsize,color={rgb, 255:red, 255; green, 255; blue, 255 }  ,opacity=1 ]  {$X$};
                \draw (157,313) node [anchor=north west][inner sep=0.75pt]    {$=$};
                \draw (122,299) node [anchor=north west][inner sep=0.75pt]  [font=\scriptsize]  {$\textcolor[rgb]{1,1,1}{f}\textcolor[rgb]{1,1,1}{^{*}}$};
            \end{tikzpicture}
        \end{center}
        and another similar diagram.
    \end{definition}
    Now let $\mathcal C$ have a monoidal structure, we now will need an additional direction to represent the direction of the tensor product. So objects will now be represented as sheets in $3$d, and we need to remember which side of the sheet is facing front (i.e., the framing which we have been suppressing).
    \begin{definition}
        For an object $X\in\mathcal C$, a right dual is an object $X^\vee$ with $1$-morphisms $\mathrm{ev}_X\colon X\boxtimes X^\vee\to I$ and $\mathrm{coev}_X\colon I\to X^\vee\boxtimes X$ and an invertible $2$-morphism $\mathrm{cusp}_X\colon\mathrm{id}_X\Rightarrow\left(\mathrm{ev}_X\boxtimes\mathrm{id}_X\right)\circ\left(\mathrm{id}_X\boxtimes\mathrm{coev}_X\right)$. We require the existence of an invertible $2$-morphism between $\mathrm{id}_{X^\vee}$ and $\left(\mathrm{id}_{X^\vee}\boxtimes\mathrm{ev}_X\right)\circ\left(\mathrm{coev}_X\boxtimes\mathrm{id}_{X^\vee}\right)$ but it is not part of the data of the dual.
    \end{definition}
    \begin{remark}
        We require what would have been `$\mathrm{cusp}^X$' to exist but not ask it to be part of the data as if a right dual exists, the $2$-groupoids of right duality data for $X$ is contractible, making the choice of $\mathrm{cusp}^X$ essentially uniquely determined by the data we have asked for.
    \end{remark}
    We visualise the sheets like origami sheets which have colours on both sides, one light and one dark. For example, if the first sheet is an object $Y$, the second sheet would be the object $Y^\vee$, and the diagram itself would represent the tensor product $Y\boxtimes Y^\vee$.
    \begin{center}
        \tikzset{every picture/.style={line width=0.75pt}}
        \begin{tikzpicture}[x=0.75pt,y=0.75pt,yscale=-1,xscale=1]
            \draw  [fill={rgb, 255:red, 208; green, 2; blue, 27 }  ,fill opacity=1 ] (360.24,40.02) -- (359.87,104) -- (282.5,83.04) -- (282.87,19.06) -- cycle ;
            \draw  [fill={rgb, 255:red, 255; green, 173; blue, 183 }  ,fill opacity=1 ] (448.24,43.02) -- (447.87,107) -- (370.5,86.04) -- (370.87,22.06) -- cycle ;
        \end{tikzpicture}
    \end{center}
    $\mathrm{coev}_Y$ and $\mathrm{cusp}_Y$ are represented by
    \begin{center}
        \tikzset{every picture/.style={line width=0.75pt}}
        \begin{tikzpicture}[x=0.75pt,y=0.75pt,yscale=-1,xscale=1]
            \draw  [fill={rgb, 255:red, 208; green, 2; blue, 27 }  ,fill opacity=1 ] (475.24,50.02) -- (474.87,114) -- (397.5,93.04) -- (397.87,29.06) -- cycle ;
            \draw  [fill={rgb, 255:red, 255; green, 173; blue, 183 }  ,fill opacity=1 ] (461.24,77.02) -- (460.87,141) -- (383.5,120.04) -- (383.87,56.06) -- cycle ;
            \draw  [dash pattern={on 0.84pt off 2.51pt}] (475.24,50.02) -- (474.87,114) -- (397.5,93.04) -- (397.87,29.06) -- cycle ;
            \draw    (383.87,56.06) .. controls (366,51) and (378.51,24.35) .. (397.87,29.06) ;
            \draw    (383.5,120.04) .. controls (371.58,117.3) and (376.46,102.56) .. (374.87,107) ;
            \draw    (375.24,43.02) -- (374.87,107) ;
            \draw  [draw opacity=0][fill={rgb, 255:red, 208; green, 2; blue, 27 }  ,fill opacity=1 ] (376.94,37.5) .. controls (380.28,30.83) and (396.79,24.96) .. (397.87,30.06) .. controls (398.94,35.17) and (414.94,64.83) .. (393.28,58.83) .. controls (371.61,52.83) and (380.28,54.83) .. (377.61,53.5) .. controls (374.94,52.17) and (373.61,44.17) .. (376.94,37.5) -- cycle ;
            \draw    (383.87,56.06) .. controls (366,51) and (378.51,24.35) .. (397.87,29.06) ;
            \draw  [draw opacity=0][fill={rgb, 255:red, 255; green, 173; blue, 183 }  ,fill opacity=1 ] (374.93,50.6) .. controls (373.92,53.03) and (416.53,70.6) .. (396.53,90.6) .. controls (376.53,110.6) and (387.24,117.57) .. (383.5,120.04) .. controls (379.76,122.5) and (373.45,112.5) .. (374.53,111) .. controls (375.61,109.5) and (375.94,48.17) .. (374.93,50.6) -- cycle ;
            \draw  [draw opacity=0][fill={rgb, 255:red, 255; green, 173; blue, 183 }  ,fill opacity=1 ] (383.49,57.23) -- (415.59,64.81) -- (415.21,78.77) -- (383.11,71.19) -- cycle ;
            \draw  [draw opacity=0][fill={rgb, 255:red, 255; green, 173; blue, 183 }  ,fill opacity=1 ] (379.37,103.18) -- (409.16,110.22) -- (411.53,126.67) -- (381.75,119.63) -- cycle ;
            \draw  [dash pattern={on 0.84pt off 2.51pt}]  (374.87,107) .. controls (375.46,101.56) and (384.46,92.56) .. (397.5,93.04) ;
            \draw    (381.75,119.63) -- (460.87,141) ;
            \draw    (383.87,56.06) .. controls (361.81,41.57) and (389.24,26.14) .. (397.87,29.06) ;
            \draw    (383.87,56.06) -- (461.24,77.02) ;
            \draw  [draw opacity=0][fill={rgb, 255:red, 255; green, 173; blue, 183 }  ,fill opacity=1 ] (379.14,80.14) -- (430.67,80.14) -- (430.67,91.29) -- (379.14,91.29) -- cycle ;
            \draw    (374.87,107) .. controls (374.45,105.91) and (373.45,118.91) .. (383.5,120.04) ;
            \draw  [draw opacity=0][fill={rgb, 255:red, 255; green, 173; blue, 183 }  ,fill opacity=1 ] (559,50.33) -- (609,50.33) -- (609,70) -- (559,70) -- cycle ;
            \draw  [draw opacity=0][fill={rgb, 255:red, 208; green, 2; blue, 27 }  ,fill opacity=1 ] (530.67,67.33) -- (620.32,67.33) -- (620.32,127) -- (530.67,127) -- cycle ;
            \draw  [draw opacity=0][fill={rgb, 255:red, 208; green, 2; blue, 27 }  ,fill opacity=1 ] (499.3,47.33) -- (560,47.33) -- (560,127) -- (499.3,127) -- cycle ;
            \draw [fill={rgb, 255:red, 208; green, 2; blue, 27 }  ,fill opacity=1 ]   (560,47.33) .. controls (591,47.33) and (588.47,57.75) .. (560,57.33) ;
            \draw [fill={rgb, 255:red, 255; green, 173; blue, 183 }  ,fill opacity=1 ]   (560,57.33) .. controls (530,57.33) and (530,67.33) .. (560,67.33) ;
            \draw    (499.33,47) -- (560,47.33) ;
            \draw    (559,67.33) -- (619.67,67.67) ;
            \draw    (500.33,127) -- (620.98,126.78) ;
            \draw    (500.33,127) -- (499.33,47) ;
            \draw    (620.98,126.78) -- (619.67,67.67) ;
            \draw    (537.58,87.15) -- (537.58,62.15) ;
            \draw  [dash pattern={on 0.84pt off 2.51pt}]  (582.58,86.15) -- (582.58,53.15) ;
            \draw  [dash pattern={on 0.84pt off 2.51pt}] (537.58,86.15) .. controls (537.58,83.39) and (547.65,81.15) .. (560.08,81.15) .. controls (572.5,81.15) and (582.58,83.39) .. (582.58,86.15) .. controls (582.58,88.91) and (572.5,91.15) .. (560.08,91.15) .. controls (547.65,91.15) and (537.58,88.91) .. (537.58,86.15) -- cycle ;
            \draw  [dash pattern={on 0.84pt off 2.51pt}]  (537.58,87.15) -- (559.71,101.65) ;
            \draw  [dash pattern={on 0.84pt off 2.51pt}]  (559.71,101.65) -- (582.58,86.15) ;
            \draw  [draw opacity=0][fill={rgb, 255:red, 255; green, 255; blue, 255 }  ,fill opacity=1 ] (583,50.33) -- (609,50.33) -- (609,67.17) -- (583,67.17) -- cycle ;
            \draw    (583,66.83) -- (582,50.33) ;
        \end{tikzpicture}
    \end{center}
    respectively and $\mathrm{ev}_Y$ is similar to $\mathrm{coev}_Y$.
    \begin{definition}
        Given a $1$-morphism $f\colon X\to Y$ and $X$ and $Y$ have right duals $X^\vee$ and $Y^\vee$, we define its right mate $f^\vee\colon Y^\vee\to X^\vee$ as a morphism equipped with two invertible $2$-morphisms
        \begin{equation*}
            \begin{split}
                \mathrm{rot}_f\colon\left(\mathrm{id}_{X^\vee}\boxtimes f\right)\circ\mathrm{coev}_X\Rightarrow&\left(f^\vee\boxtimes\mathrm{id}_Y\right)\circ\mathrm{coev}_Y\\
                \mathrm{rot}^f\colon\mathrm{ev}_Y\circ\left(f\boxtimes\mathrm{id}_{Y^\vee}\right)\Rightarrow&\;\mathrm{ev}_X\circ\left(\mathrm{id}_X\boxtimes f^\vee\right)
            \end{split}
        \end{equation*}
        which satisfy obvious coherences with the $\mathrm{cusp}$ maps of $X$ and $Y$.
    \end{definition}
    We represent $\mathrm{rot}_f$ graphically as 
    \begin{center}
        \tikzset{every picture/.style={line width=0.75pt}}
        \begin{tikzpicture}[x=0.75pt,y=0.75pt,yscale=-1,xscale=1]
            \draw  [fill={rgb, 255:red, 208; green, 2; blue, 27 }  ,fill opacity=1 ] (253.24,250.02) -- (252.87,314) -- (175.5,293.04) -- (175.87,229.06) -- cycle ;
            \draw  [fill={rgb, 255:red, 157; green, 194; blue, 127 }  ,fill opacity=1 ] (239.24,277.02) -- (238.87,341) -- (161.5,320.04) -- (161.87,256.06) -- cycle ;
            \draw  [dash pattern={on 0.84pt off 2.51pt}] (253.24,250.02) -- (252.87,314) -- (175.5,293.04) -- (175.87,229.06) -- cycle ;
            \draw    (161.87,256.06) .. controls (144,251) and (156.51,224.35) .. (175.87,229.06) ;
            \draw    (161.5,320.04) .. controls (149.58,317.3) and (154.46,302.56) .. (152.87,307) ;
            \draw    (153.24,243.02) -- (152.87,307) ;
            \draw  [draw opacity=0][fill={rgb, 255:red, 208; green, 2; blue, 27 }  ,fill opacity=1 ] (154.94,236.5) .. controls (158.28,229.83) and (174.79,223.96) .. (175.87,229.06) .. controls (176.94,234.17) and (192.94,263.83) .. (171.28,257.83) .. controls (149.61,251.83) and (158.28,253.83) .. (155.61,252.5) .. controls (152.94,251.17) and (151.61,243.17) .. (154.94,236.5) -- cycle ;
            \draw    (161.87,256.06) .. controls (144,251) and (156.51,224.35) .. (175.87,229.06) ;
            \draw  [draw opacity=0][fill={rgb, 255:red, 157; green, 194; blue, 127 }  ,fill opacity=1 ] (152.93,250.6) .. controls (151.92,253.03) and (194.53,270.6) .. (174.53,290.6) .. controls (154.53,310.6) and (165.24,317.57) .. (161.5,320.04) .. controls (157.76,322.5) and (151.45,312.5) .. (152.53,311) .. controls (153.61,309.5) and (153.94,248.17) .. (152.93,250.6) -- cycle ;
            \draw  [draw opacity=0][fill={rgb, 255:red, 157; green, 194; blue, 127 }  ,fill opacity=1 ] (161.49,257.23) -- (193.59,264.81) -- (193.21,278.77) -- (161.11,271.19) -- cycle ;
            \draw  [draw opacity=0][fill={rgb, 255:red, 157; green, 194; blue, 127 }  ,fill opacity=1 ] (159.13,303.59) -- (188.91,310.62) -- (191.28,327.07) -- (161.5,320.04) -- cycle ;
            \draw  [dash pattern={on 0.84pt off 2.51pt}]  (152.87,307) .. controls (153.46,301.56) and (162.46,292.56) .. (175.5,293.04) ;
            \draw    (159.75,319.63) -- (238.87,341) ;
            \draw    (161.87,256.06) .. controls (139.81,241.57) and (167.24,226.14) .. (175.87,229.06) ;
            \draw    (161.87,256.06) -- (239.24,277.02) ;
            \draw  [draw opacity=0][fill={rgb, 255:red, 157; green, 194; blue, 127 }  ,fill opacity=1 ] (157.14,279.14) -- (208.67,279.14) -- (208.67,290.29) -- (157.14,290.29) -- cycle ;
            \draw    (152.87,307) .. controls (152.45,305.91) and (151.45,318.91) .. (161.5,320.04) ;
            \draw  [dash pattern={on 0.84pt off 2.51pt}]  (153.05,275.01) .. controls (156.77,287.71) and (188.77,275.71) .. (187.77,296.71) ;
            \draw [shift={(166.92,282.06)}, rotate = 7.59] [fill={rgb, 255:red, 0; green, 0; blue, 0 }  ][line width=0.08]  [draw opacity=0] (8.93,-4.29) -- (0,0) -- (8.93,4.29) -- cycle    ;
            \draw [fill={rgb, 255:red, 255; green, 173; blue, 183 }  ,fill opacity=1 ]   (153.05,275.01) .. controls (159.52,268.02) and (187.52,285.02) .. (187.52,264.02) ;
            \draw  [draw opacity=0][fill={rgb, 255:red, 255; green, 173; blue, 183 }  ,fill opacity=1 ] (152.93,250.6) .. controls (151.28,253.46) and (187.71,263.37) .. (187.52,264.02) .. controls (187.33,264.67) and (186.33,265.67) .. (172.59,270.74) .. controls (158.84,275.81) and (152.77,272.36) .. (153.05,275.01) .. controls (153.33,277.67) and (154.59,247.74) .. (152.93,250.6) -- cycle ;
            \draw [fill={rgb, 255:red, 255; green, 173; blue, 183 }  ,fill opacity=1 ]   (153.05,275.01) .. controls (159.52,268.02) and (187.52,285.02) .. (187.52,264.02) ;
            \draw [shift={(178.53,274.86)}, rotate = 181.11] [fill={rgb, 255:red, 0; green, 0; blue, 0 }  ][line width=0.08]  [draw opacity=0] (8.93,-4.29) -- (0,0) -- (8.93,4.29) -- cycle    ;
            \draw    (161.87,256.06) -- (239.24,277.02) ;
            \draw    (161.87,256.06) .. controls (139.81,241.57) and (167.24,226.14) .. (175.87,229.06) ;
            \draw (230,250) node [anchor=north west][inner sep=0.75pt]    {$\textcolor[rgb]{1,1,1}{Y}$};
            \draw (214,311) node [anchor=north west][inner sep=0.75pt]    {$\textcolor[rgb]{1,1,1}{X^{\lor }}$};
            \draw (185,269) node [anchor=north west][inner sep=0.75pt]  [font=\tiny,color={rgb, 255:red, 255; green, 255; blue, 255 }  ,opacity=1 ]  {$f$};
        \end{tikzpicture}
    \end{center}
    and $\mathrm{rot}^f$ similarly.
    \begin{remark}
        Left adjoints, duals, and mates can be defined in a similar way. Adjoints, duals, and mates exist up to a contractible choice whenever they do exist.
    \end{remark}
    \begin{definition}
        We define the category $\mathcal C^\mathrm{coh}$ to be the category whose objects are objects of $\mathcal C$ equipped with right duality data, $1$-morphisms are $1$-morphisms of $\mathcal C$ with right adjunction data and right mate data, and the same $2$-morphisms as those in $\mathcal C$.

        We can similarly define the category ${}^\mathrm{coh}\mathcal C^\mathrm{coh}$ having the objects of $\mathcal C$ equipped with right and left duality data, and so on.
    \end{definition}
    The forgetful functor $\mathcal C^\mathrm{coh}\to\mathcal C$ is an equivalence if all objects and morphisms have right duals and adjoints, as then the fibres of the forgetful functor are all contractible. The same is true of the forgetful functor $\mathrm h_1{}^\mathrm{coh}\mathcal C^\mathrm{coh}\to\mathcal C$ if all objects have left and right duals and all morphisms have left and right  adjoints.
    \begin{remark}
        We note that since adjoints and mates are unique up to invertible $2$-morphisms, the category $\mathrm h_1\mathcal C^\mathrm{coh}$ has the same objects of $\mathcal C^\mathrm{coh}$ but the $1$-morphisms are isomorphism classes of $1$-morphisms in $\mathcal C$.
    \end{remark}
    We can view $\mathrm h_1{}^\mathrm{coh}\mathcal C^\mathrm{coh}$ as a (strict) homotopy pullback
    \begin{equation*}
        \begin{tikzcd}
            \mathrm h_1{}^\mathrm{coh}\mathcal C^\mathrm{coh}\arrow[rrr]\arrow[d]&&&\mathrm h_1\mathcal C^\mathrm{coh}\arrow[d,"\left(Y{,}\dots\right)\mapsto Y{,}g\mapsto g"]\\
            \mathrm h_1\mathcal C^\mathrm{coh}\arrow[rrr,"\left(Y{,}\dots\right)\mapsto Y^\vee{,}g\mapsto\left(g^*\right)^\vee"']&&&\mathrm h_1\mathcal C
        \end{tikzcd}
    \end{equation*}
    \section{The 2-categorical Quantum Dimension}
    Let $\mathcal C$ be a braided monoidal category now.
    \begin{definition}
        We define
        \begin{equation*}
            \eta\colon\mathrm h_1\mathcal C^\mathrm{coh}\to\Omega\mathcal C
        \end{equation*}
        as
        \begin{itemize}
            \item for an object $\left(X,X^\vee,\mathrm{ev}_X,\mathrm{coev}_X,\mathrm{cusp}_X\right)$, define
            \begin{equation*}
                \eta\left(x\right)=\mathrm{ev}_X\circ b_{X^\vee,X}\circ\mathrm{coev}_X
            \end{equation*}
            \item for a $1$-morphism in $\mathrm h_1\mathcal C^\mathrm{coh}$ from $\left(X,X^\vee,\mathrm{ev}_X,\mathrm{coev}_X,\mathrm{cusp}_X\right)\to\left(Y,Y^\vee,\mathrm{ev}_Y,\mathrm{coev}_Y,\mathrm{cusp}_Y\right)$, pick a representative $f\colon X\to Y$ and pick a right adjoint and right mate, and define $\eta\left(\left[f\right]\right)$
            \begin{center}
                \tikzset{every picture/.style={line width=0.75pt}}
                \begin{tikzpicture}[x=0.75pt,y=0.75pt,yscale=-1,xscale=1]
                    \draw [color={rgb, 255:red, 208; green, 2; blue, 27 }  ,draw opacity=1 ]   (11.33,41.33) .. controls (11.33,10.33) and (61.33,10.33) .. (61.33,41.33) ;
                    \draw [shift={(41.55,18.54)}, rotate = 184.28] [fill={rgb, 255:red, 208; green, 2; blue, 27 }  ,fill opacity=1 ][line width=0.08]  [draw opacity=0] (8.93,-4.29) -- (0,0) -- (8.93,4.29) -- cycle    ;
                    \draw [color={rgb, 255:red, 208; green, 2; blue, 27 }  ,draw opacity=1 ]   (11.33,80.33) .. controls (10.67,109.33) and (60.67,110.33) .. (61.33,80.33) ;
                    \draw [color={rgb, 255:red, 208; green, 2; blue, 27 }  ,draw opacity=1 ]   (11.33,41.33) .. controls (11.33,51) and (61.33,71) .. (61.33,80.33) ;
                    \draw [color={rgb, 255:red, 208; green, 2; blue, 27 }  ,draw opacity=1 ]   (11.33,80.33) .. controls (10.67,71) and (23.67,68) .. (32.67,63) ;
                    \draw [color={rgb, 255:red, 208; green, 2; blue, 27 }  ,draw opacity=1 ]   (61.33,41.33) .. controls (61.06,51.49) and (40.06,58.49) .. (40.06,58.49) ;
                    \draw [color={rgb, 255:red, 208; green, 2; blue, 27 }  ,draw opacity=1 ]   (102.33,41.33) .. controls (102.33,10.33) and (152.33,10.33) .. (152.33,41.33) ;
                    \draw [shift={(132.55,18.54)}, rotate = 184.28] [fill={rgb, 255:red, 208; green, 2; blue, 27 }  ,fill opacity=1 ][line width=0.08]  [draw opacity=0] (8.93,-4.29) -- (0,0) -- (8.93,4.29) -- cycle    ;
                    \draw [color={rgb, 255:red, 208; green, 2; blue, 27 }  ,draw opacity=1 ]   (102.33,80.33) .. controls (101.67,109.33) and (133.67,103) .. (139.67,98) ;
                    \draw [color={rgb, 255:red, 208; green, 2; blue, 27 }  ,draw opacity=1 ]   (102.33,41.33) .. controls (102.33,51) and (132.67,61) .. (138.67,69) ;
                    \draw [color={rgb, 255:red, 208; green, 2; blue, 27 }  ,draw opacity=1 ]   (102.33,80.33) .. controls (101.67,71) and (114.67,68) .. (123.67,63) ;
                    \draw [color={rgb, 255:red, 208; green, 2; blue, 27 }  ,draw opacity=1 ]   (152.33,41.33) .. controls (152.06,51.49) and (131.06,58.49) .. (131.06,58.49) ;
                    \draw [color={rgb, 255:red, 65; green, 117; blue, 5 }  ,draw opacity=1 ]   (138.67,69) .. controls (147.67,77) and (151.67,91) .. (139.67,98) ;
                    \draw [shift={(145.12,77.47)}, rotate = 72.34] [fill={rgb, 255:red, 65; green, 117; blue, 5 }  ,fill opacity=1 ][line width=0.08]  [draw opacity=0] (8.93,-4.29) -- (0,0) -- (8.93,4.29) -- cycle    ;
                    \draw  [fill={rgb, 255:red, 0; green, 0; blue, 0 }  ,fill opacity=1 ] (137.67,97) .. controls (137.67,95.7) and (138.72,94.65) .. (140.01,94.65) .. controls (141.31,94.65) and (142.36,95.7) .. (142.36,97) .. controls (142.36,98.3) and (141.31,99.35) .. (140.01,99.35) .. controls (138.72,99.35) and (137.67,98.3) .. (137.67,97) -- cycle ;
                    \draw  [fill={rgb, 255:red, 0; green, 0; blue, 0 }  ,fill opacity=1 ] (136.67,70) .. controls (136.67,68.7) and (137.72,67.65) .. (139.01,67.65) .. controls (140.31,67.65) and (141.36,68.7) .. (141.36,70) .. controls (141.36,71.3) and (140.31,72.35) .. (139.01,72.35) .. controls (137.72,72.35) and (136.67,71.3) .. (136.67,70) -- cycle ;
                    \draw [color={rgb, 255:red, 208; green, 2; blue, 27 }  ,draw opacity=1 ]   (193.33,41.33) .. controls (193.33,10.33) and (243.33,10.33) .. (243.33,41.33) ;
                    \draw [shift={(223.55,18.54)}, rotate = 184.28] [fill={rgb, 255:red, 208; green, 2; blue, 27 }  ,fill opacity=1 ][line width=0.08]  [draw opacity=0] (8.93,-4.29) -- (0,0) -- (8.93,4.29) -- cycle    ;
                    \draw [color={rgb, 255:red, 208; green, 2; blue, 27 }  ,draw opacity=1 ]   (193.33,41.33) .. controls (193.33,51) and (230.3,65.9) .. (236.3,73.9) ;
                    \draw [color={rgb, 255:red, 208; green, 2; blue, 27 }  ,draw opacity=1 ]   (193.33,80.33) .. controls (192.67,71) and (205.67,68) .. (214.67,63) ;
                    \draw [color={rgb, 255:red, 208; green, 2; blue, 27 }  ,draw opacity=1 ]   (243.33,41.33) .. controls (243.06,51.49) and (222.06,58.49) .. (222.06,58.49) ;
                    \draw [color={rgb, 255:red, 65; green, 117; blue, 5 }  ,draw opacity=1 ]   (236.3,73.9) .. controls (247.67,99) and (200.67,117) .. (193.33,80.33) ;
                    \draw [shift={(225.44,99)}, rotate = 162.85] [fill={rgb, 255:red, 65; green, 117; blue, 5 }  ,fill opacity=1 ][line width=0.08]  [draw opacity=0] (8.93,-4.29) -- (0,0) -- (8.93,4.29) -- cycle    ;
                    \draw  [fill={rgb, 255:red, 0; green, 0; blue, 0 }  ,fill opacity=1 ] (234.96,76.24) .. controls (234.96,74.95) and (236.01,73.9) .. (237.3,73.9) .. controls (238.6,73.9) and (239.65,74.95) .. (239.65,76.24) .. controls (239.65,77.54) and (238.6,78.59) .. (237.3,78.59) .. controls (236.01,78.59) and (234.96,77.54) .. (234.96,76.24) -- cycle ;
                    \draw  [fill={rgb, 255:red, 0; green, 0; blue, 0 }  ,fill opacity=1 ] (191.33,80.33) .. controls (191.33,79.04) and (192.38,77.99) .. (193.68,77.99) .. controls (194.98,77.99) and (196.03,79.04) .. (196.03,80.33) .. controls (196.03,81.63) and (194.98,82.68) .. (193.68,82.68) .. controls (192.38,82.68) and (191.33,81.63) .. (191.33,80.33) -- cycle ;
                    \draw [color={rgb, 255:red, 208; green, 2; blue, 27 }  ,draw opacity=1 ]   (300.33,41.33) .. controls (300.33,10.33) and (350.33,10.33) .. (350.33,41.33) ;
                    \draw [shift={(330.55,18.54)}, rotate = 184.28] [fill={rgb, 255:red, 208; green, 2; blue, 27 }  ,fill opacity=1 ][line width=0.08]  [draw opacity=0] (8.93,-4.29) -- (0,0) -- (8.93,4.29) -- cycle    ;
                    \draw [color={rgb, 255:red, 65; green, 117; blue, 5 }  ,draw opacity=1 ]   (300.33,80.33) .. controls (299.67,109.33) and (349.67,110.33) .. (350.33,80.33) ;
                    \draw [shift={(330.43,102.02)}, rotate = 175.81] [fill={rgb, 255:red, 65; green, 117; blue, 5 }  ,fill opacity=1 ][line width=0.08]  [draw opacity=0] (8.93,-4.29) -- (0,0) -- (8.93,4.29) -- cycle    ;
                    \draw [color={rgb, 255:red, 65; green, 117; blue, 5 }  ,draw opacity=1 ]   (300.33,41.33) .. controls (300.33,51) and (350.33,71) .. (350.33,80.33) ;
                    \draw [color={rgb, 255:red, 65; green, 117; blue, 5 }  ,draw opacity=1 ]   (300.33,80.33) .. controls (299.67,71) and (312.67,68) .. (321.67,63) ;
                    \draw [color={rgb, 255:red, 65; green, 117; blue, 5 }  ,draw opacity=1 ]   (350.33,41.33) .. controls (350.06,51.49) and (329.06,58.49) .. (329.06,58.49) ;
                    \draw  [fill={rgb, 255:red, 0; green, 0; blue, 0 }  ,fill opacity=1 ] (298.99,42.68) .. controls (298.99,41.38) and (300.04,40.33) .. (301.33,40.33) .. controls (302.63,40.33) and (303.68,41.38) .. (303.68,42.68) .. controls (303.68,43.98) and (302.63,45.03) .. (301.33,45.03) .. controls (300.04,45.03) and (298.99,43.98) .. (298.99,42.68) -- cycle ;
                    \draw  [fill={rgb, 255:red, 0; green, 0; blue, 0 }  ,fill opacity=1 ] (347.64,41.33) .. controls (347.64,40.04) and (348.69,38.99) .. (349.99,38.99) .. controls (351.28,38.99) and (352.33,40.04) .. (352.33,41.33) .. controls (352.33,42.63) and (351.28,43.68) .. (349.99,43.68) .. controls (348.69,43.68) and (347.64,42.63) .. (347.64,41.33) -- cycle ;
                    \draw [color={rgb, 255:red, 208; green, 2; blue, 27 }  ,draw opacity=1 ]   (415,50) .. controls (399,43) and (411,29) .. (414,27) ;
                    \draw [shift={(408.26,34.33)}, rotate = 103.46] [fill={rgb, 255:red, 208; green, 2; blue, 27 }  ,fill opacity=1 ][line width=0.08]  [draw opacity=0] (8.93,-4.29) -- (0,0) -- (8.93,4.29) -- cycle    ;
                    \draw [color={rgb, 255:red, 65; green, 117; blue, 5 }  ,draw opacity=1 ]   (407.33,81.33) .. controls (406.67,110.33) and (456.67,111.33) .. (457.33,81.33) ;
                    \draw [shift={(437.43,103.02)}, rotate = 175.81] [fill={rgb, 255:red, 65; green, 117; blue, 5 }  ,fill opacity=1 ][line width=0.08]  [draw opacity=0] (8.93,-4.29) -- (0,0) -- (8.93,4.29) -- cycle    ;
                    \draw [color={rgb, 255:red, 65; green, 117; blue, 5 }  ,draw opacity=1 ]   (415,50) .. controls (418,55) and (457.33,72) .. (457.33,81.33) ;
                    \draw [color={rgb, 255:red, 65; green, 117; blue, 5 }  ,draw opacity=1 ]   (407.33,81.33) .. controls (406.67,72) and (419.67,69) .. (428.67,64) ;
                    \draw [color={rgb, 255:red, 65; green, 117; blue, 5 }  ,draw opacity=1 ]   (457.33,42.33) .. controls (457.06,52.49) and (436.06,59.49) .. (436.06,59.49) ;
                    \draw  [fill={rgb, 255:red, 0; green, 0; blue, 0 }  ,fill opacity=1 ] (413,50) .. controls (413,48.7) and (414.05,47.65) .. (415.35,47.65) .. controls (416.64,47.65) and (417.69,48.7) .. (417.69,50) .. controls (417.69,51.3) and (416.64,52.35) .. (415.35,52.35) .. controls (414.05,52.35) and (413,51.3) .. (413,50) -- cycle ;
                    \draw  [fill={rgb, 255:red, 0; green, 0; blue, 0 }  ,fill opacity=1 ] (409.31,28) .. controls (409.31,26.7) and (410.36,25.65) .. (411.65,25.65) .. controls (412.95,25.65) and (414,26.7) .. (414,28) .. controls (414,29.3) and (412.95,30.35) .. (411.65,30.35) .. controls (410.36,30.35) and (409.31,29.3) .. (409.31,28) -- cycle ;
                    \draw [color={rgb, 255:red, 65; green, 117; blue, 5 }  ,draw opacity=1 ]   (414,27) .. controls (440,2.33) and (458,36.33) .. (457.33,42.33) ;
                    \draw [color={rgb, 255:red, 65; green, 117; blue, 5 }  ,draw opacity=1 ]   (494.33,42.33) .. controls (494.33,11.33) and (544.33,11.33) .. (544.33,42.33) ;
                    \draw [shift={(524.55,19.54)}, rotate = 184.28] [fill={rgb, 255:red, 65; green, 117; blue, 5 }  ,fill opacity=1 ][line width=0.08]  [draw opacity=0] (8.93,-4.29) -- (0,0) -- (8.93,4.29) -- cycle    ;
                    \draw [color={rgb, 255:red, 65; green, 117; blue, 5 }  ,draw opacity=1 ]   (494.33,81.33) .. controls (493.67,110.33) and (543.67,111.33) .. (544.33,81.33) ;
                    \draw [color={rgb, 255:red, 65; green, 117; blue, 5 }  ,draw opacity=1 ]   (494.33,42.33) .. controls (494.33,52) and (544.33,72) .. (544.33,81.33) ;
                    \draw [color={rgb, 255:red, 65; green, 117; blue, 5 }  ,draw opacity=1 ]   (494.33,81.33) .. controls (493.67,72) and (506.67,69) .. (515.67,64) ;
                    \draw [color={rgb, 255:red, 65; green, 117; blue, 5 }  ,draw opacity=1 ]   (544.33,42.33) .. controls (544.06,52.49) and (523.06,59.49) .. (523.06,59.49) ;
                    \draw (60,46) node [anchor=north west][inner sep=0.75pt]    {$\xrightarrow{\mathrm{coev}_f}$};
                    \draw (152.33,46) node [anchor=north west][inner sep=0.75pt]    {$\xrightarrow{\mathrm{rot}_f}$};
                    \draw (240,46) node [anchor=north west][inner sep=0.75pt]    {$\xrightarrow[\text{ of }\beta]{\text{naturality}}$};
                    \draw (355,42.33) node [anchor=north west][inner sep=0.75pt]    {$\xrightarrow{\left(\mathrm{rot}^f\right)^{-1}}$};
                    \draw (460,48) node [anchor=north west][inner sep=0.75pt]    {$\xrightarrow{\mathrm{ev}_f}$};
                    \draw (56,92) node [anchor=north west][inner sep=0.75pt]  [font=\scriptsize]  {$\textcolor[rgb]{0.82,0.01,0.11}{x}$};
                    \draw (97,91) node [anchor=north west][inner sep=0.75pt]  [font=\scriptsize]  {$\textcolor[rgb]{0.82,0.01,0.11}{x}$};
                    \draw (239,16) node [anchor=north west][inner sep=0.75pt]  [font=\scriptsize]  {$\textcolor[rgb]{0.82,0.01,0.11}{x}$};
                    \draw (346,16) node [anchor=north west][inner sep=0.75pt]  [font=\scriptsize]  {$\textcolor[rgb]{0.82,0.01,0.11}{x}$};
                    \draw (394,29) node [anchor=north west][inner sep=0.75pt]  [font=\scriptsize]  {$\textcolor[rgb]{0.82,0.01,0.11}{x}$};
                    \draw (148,73) node [anchor=north west][inner sep=0.75pt]  [font=\scriptsize,color={rgb, 255:red, 65; green, 117; blue, 5 }  ,opacity=1 ]  {$\textcolor[rgb]{0.25,0.46,0.02}{y}$};
                    \draw (202,100) node [anchor=north west][inner sep=0.75pt]  [font=\scriptsize,color={rgb, 255:red, 65; green, 117; blue, 5 }  ,opacity=1 ]  {$\textcolor[rgb]{0.25,0.46,0.02}{y}$};
                    \draw (330,101) node [anchor=north west][inner sep=0.75pt]  [font=\scriptsize,color={rgb, 255:red, 65; green, 117; blue, 5 }  ,opacity=1 ]  {$\textcolor[rgb]{0.25,0.46,0.02}{y}$};
                    \draw (437,102) node [anchor=north west][inner sep=0.75pt]  [font=\scriptsize,color={rgb, 255:red, 65; green, 117; blue, 5 }  ,opacity=1 ]  {$\textcolor[rgb]{0.25,0.46,0.02}{y}$};
                    \draw (518,103) node [anchor=north west][inner sep=0.75pt]  [font=\scriptsize,color={rgb, 255:red, 65; green, 117; blue, 5 }  ,opacity=1 ]  {$\textcolor[rgb]{0.25,0.46,0.02}{y}$};
                    \draw (140.01,94.65) node [anchor=north west][inner sep=0.75pt]  [font=\scriptsize]  {$f$};
                    \draw (143,63) node [anchor=north west][inner sep=0.75pt]  [font=\scriptsize]  {$f^{*}$};
                    \draw (238,70) node [anchor=north west][inner sep=0.75pt]  [font=\scriptsize]  {$f^{*}$};
                    \draw (177.01,71.65) node [anchor=north west][inner sep=0.75pt]  [font=\scriptsize]  {$f^{\vee }$};
                    \draw (351.01,33.65) node [anchor=north west][inner sep=0.75pt]  [font=\scriptsize]  {$f^{\vee }$};
                    \draw (287.67,32) node [anchor=north west][inner sep=0.75pt]  [font=\scriptsize]  {$f^{*}$};
                    \draw (416.67,39) node [anchor=north west][inner sep=0.75pt]  [font=\scriptsize]  {$f^{*}$};
                    \draw (400.06,17.49) node [anchor=north west][inner sep=0.75pt]  [font=\scriptsize]  {$f$};
                \end{tikzpicture}
            \end{center}
            We view this as slices of a $2$d dimensional surface embedded in $4$d representing a $2$-morphism.
        \end{itemize}
    \end{definition}
    As before, $\eta\left(f\right)$ is independent of the choices of representatives, mates, and adjoints. Similarly, the isomorphism class of $\eta\left(X\right)$ does not depend on the choice of duality data, and equivalent objects have isomorphic quantum dimensions.

    $\eta$ is \emph{not} monoidal. For example, picking an invertible object $Y$ and an invertible $1$-morphism $x\colon I\to I$, we have
    \begin{equation}
        \label{eq:2etanotmonoidal}
        \eta\left(x\boxtimes\mathrm{id}_Y\right)=\left\langle\beta_{\mathrm{id}_Y,x}\cdot\beta_{\mathrm{x,id}_Y}\right\rangle\eta\left(x\right)\mathrm{id}_{\eta\left(Y\right)}
    \end{equation}
    \begin{remark}
        Let $x\colon I\to I$ and equip $I$ with the trivial duality data ($I^\vee=I$ and the duality morphisms are all given by unitors), then $\eta\left(x\right)$ computed as a morphism in $\mathcal C$ or as an object in $\Omega\mathcal C$ agree.
    \end{remark}
    \begin{definition}
        Let $X$ be an object of $\mathcal C$ with a left dual. We define a natural isomorphism
        \begin{equation*}
            \mathrm{fl}_X\colon\eta\left(x\right)\to\eta\left({}^\vee X\right)
        \end{equation*}
        as
        \begin{equation*}
            \begin{split}
                \eta\left({}^\vee X\right)=&\eta\left({}^\vee X\right)^\vee=\left(\mathrm{ev}_{{}^\vee X}\circ\beta_{\left({}^\vee X\right)^\vee,{}^\vee X}\circ\mathrm{coev}_{{}^\vee X}\right)^\vee\simeq\left(\mathrm{ev}_{{}^\vee X}\circ\beta_{X,{}^\vee X}\circ\mathrm{coev}_{{}^\vee X}\right)^\vee\\
                \simeq&\left(\mathrm{coev}_{{}^\vee X}\right)^\vee\circ\left(\beta_{X,{}^\vee X}\right)^\vee\circ\left(\mathrm{ev}_{{}^\vee X}\right)^\vee\simeq\mathrm{ev}_X\circ\beta_{X^\vee,X}\circ\mathrm{coev}_X=\eta\left(X\right)
            \end{split}
        \end{equation*}
        The first equality follows from the fact that $\eta\left({}^\vee X\right)$ is an endomorphism of $I$ and $I^\vee=I$, so all endomorphisms are their own mates. The fifth equivalence follows from the diagrammatic algebra.

        This in fact defines a natural isomorphism
        \begin{equation*}
            \begin{tikzcd}
                &\mathrm h_1\mathcal C^\mathrm{coh}\arrow[dr,"\eta"]\arrow[dd,Rightarrow,"\mathrm{fl}","\sim"',shorten=5mm]\\
                \mathrm h_1{}^\mathrm{coh}\mathcal C^\mathrm{coh}\arrow[ur,"X\mapsto X"]\arrow[dr,"X\mapsto{}^\vee X"']&&\Omega\mathcal C\\
                &\mathrm h_1\mathcal C^\mathrm{coh}\arrow[ur,"\eta"']
            \end{tikzcd}
        \end{equation*}
    \end{definition}
    \section{Computing the quantum dimension in some examples}
    Consider the case when $\mathcal C$ is $\mathcal S$ or $\mathcal T$. Their simple objects are $I,C,M,C\boxtimes M$. We compute
    \begin{itemize}
        \item$\eta\left(I\right)=\mathrm{id}_I$
        \item$\eta\left(C\right)=\mathrm{id}_e$ as $C^\vee=C$ and $\beta_{C,C}=e\boxtimes\mathrm{id}_{C\boxtimes C}$.
        \item$\eta\left(M\right)=\mathrm{id}_I$ as $M^\vee=M$ and $\beta_{M,M}=\mathrm{id}_{M\boxtimes M}$.
    \end{itemize}
    Now we compute $\eta\left(C\boxtimes M\right)$. Using Equation~\ref{eq:2etanotmonoidal}, we have
    \begin{equation*}
        \begin{split}
            \eta\left(e\boxtimes\mathrm{id}_C\right)=&\left\langle\beta_{\mathrm{id}_C,e}\cdot\beta_{e,\mathrm{id}_C}\right\rangle\eta\left(e\right)\mathrm{id}_{\eta\left(C\right)}=\left(1\right)\left(-1\right)\mathrm{id}_{\eta\left(C\right)}=-\mathrm{id}_{\eta\left(C\right)}\\
            \eta\left(e\boxtimes\mathrm{id}_M\right)=&\left\langle\beta_{\mathrm{id}_M,e}\cdot\beta_{e,\mathrm{id}_M}\right\rangle\eta\left(e\right)\mathrm{id}_{\eta\left(M\right)}=\left(-1\right)\left(-1\right)\mathrm{id}_{\eta\left(M\right)}=\mathrm{id}_{\eta\left(M\right)}
        \end{split}
    \end{equation*}
    We have $v_1\colon I\to C$ and $v_2\colon C\to I$ the simple Clifford modules, then $v_2\circ v_1$ is $C$ viewed as an element of $\Omega\mathcal S=\Omega\mathcal T=\mathrm{sVect}$, and thus $v_2\cdot v_1\simeq\mathrm{id}_I\oplus\left(e\boxtimes\mathrm{id}_I\right)$. Consider the maps
    \begin{equation*}
        f_1=v_1\boxtimes\mathrm{id}_M\colon M\to C\boxtimes M\quad\quad f_2=v_2\boxtimes\mathrm{id}_M\colon C\boxtimes M\to M 
    \end{equation*}
    Then $f_2\cdot f_1\simeq\left(\mathrm{id}_I\boxtimes e\right)\oplus\mathrm{id}_M$. Thus,
    \begin{equation*}
        \eta\left(f_2\right)\cdot\eta\left(f_1\right)=\eta\left(f_2\cdot f_1\right)\simeq\eta\left(\left(\mathrm{id}_I\boxtimes e\right)\oplus\mathrm{id}_M\right)\simeq2\mathrm{id}_M
    \end{equation*}
    Thus $\eta\left(f_2\right)\cdot\eta\left(f_1\right)$ is an isomorphism, and similarly, so is $\eta\left(f_1\right)\cdot\eta\left(f_2\right)$. Thus $\eta\left(f_1\right)$ is an isomorphism, and thus $\eta\left(C\boxtimes M\right)\simeq I$.

%% file: talks/4.3/main.tex
\input{lol.tikzstyles}
\input{style.tikzstyles}

\newcommand{\ldual}[1]{^\vee#1}
\newcommand{\uK}{\underline{\kappa}}
\newcommand{\ue}{\underline{\eta}}

\title{ The Klein invariant }
Talk by Lorenzo Riva, notes by by Tudor Caba.

\section{Recap}

Let $\mathcal{B}$ and $\mathcal{C}$ be braided fusion 1-categories. In Adria's talk we saw that there is a correspondence 

\[\begin{tikzcd}
	\begin{array}{c} \left \{ \; \begin{matrix} \text{Nondegenerate extensions} \\ \iota: \mathcal{B}\hookrightarrow  \mathcal{M} \\ \text{with } \mathcal{C}^\text{rev} \simeq \mathcal{Z}_2(\iota: \mathcal{B}\hookrightarrow \mathcal{M})  \end{matrix} \; \right \} \end{array} & \begin{array}{c} \left\{ \; \begin{matrix} \text{Equivalences} \\ \mathcal{Z}(\Sigma \mathcal{B}) \simeq \mathcal{Z}(\Sigma \mathcal{C}) \\ \text{of braided fusion 2-cats} \end{matrix} \; \right\} \end{array}
	\arrow["{\huge \simeq}", shift right=2, draw=none, from=1-1, to=1-2]
\end{tikzcd}\]

Assume now that $\mathcal{B}$ is slightly degenerate, i.e. $\mathcal{Z}_2 (\mathcal{B}) \simeq \text{sVect}$. If we could show that $\mathcal{Z}(\Sigma \mathcal{B})$ is independent of $\mathcal{B}$, then in particular we'd have $\mathcal{Z}(\Sigma \mathcal{B}) \simeq \mathcal{Z}(\Sigma \text{sVect} )$ and the above correspondence would provide us with a nondegenerate extension $\iota: \mathcal{B}\xhookrightarrow{} \mathcal{M}$ with $\text{sVect}^\text{rev} \simeq \mathcal{Z}_2(\iota)$. It would follow that $\mathcal{Z}_2(\mathcal{B}) \simeq \mathcal{Z}_2(\iota: \mathcal{B}\to \mathcal{M})$ and so $\mathcal{M}$ would be a minimally nondegenerate extension of $\mathcal{B}$. \\

We know from Cameron's talk that if $\mathcal{B}$ is slightly degenerate then $\mathcal{Z}(\Sigma \mathcal{B})$ can only be one of two 2-categories $\mathcal{S}$ and $\mathcal{T}$ (up to braided monoidal equivalence) which are both of the form $2\text{Vec}^\alpha[\mathcal{G}]$ for some twisting cocycle $\alpha$ and (the same) groupoid $\mathcal{G}$. Of course, if it happens that $\mathcal{S}\simeq \mathcal{T}$ then we'd be done. Sadly, things are not that easy, and in this talk we'll construct an invariant called the \emph{Klein invariant} $k$ which 
\begin{enumerate}
    \item Allows us to show that $\mathcal{S}$ and $\mathcal{T}$ are not equivalent, and
    \item Will eventually be used to disallow one of $\mathcal{S}$ and $\mathcal{T}$ from occurring. 
\end{enumerate}

\section{The Klein invariant}
Fix a braided monoidal $2$-category $\mathcal{C}$. Nivedita constructed the eta invariant functor 
$$\eta(-): h_1 \mathcal{C} \to \Omega \mathcal{C}$$

which comes with the natural flip isomorphism 
$$\operatorname{fl}_x: \eta(\ldual{x}) \xrightarrow[]{\sim} \eta(x).$$

\begin{definition} Let $x$ be a fully dualizable object in $\mathcal{C}$ equipped with an isomorphism $r: x \xrightarrow{\simeq} \ldual{x}$ (called \emph{duality datum}). The \emph{Klein invariant} of this data is 
$$\kappa(\mathcal{C}, X, r) := \operatorname{tr}(\operatorname{fl}_x \cdot 
 \eta(r)) \in \Omega^2 \mathcal{C}$$
\end{definition}

Note that $\operatorname{fl}_x \cdot \eta(r) = ( \eta(x) \xrightarrow[]{\sim} \eta(\ldual{x}) \xrightarrow[]{\sim} \eta(x) )$ is an automorphism of $\eta(x)$.
\begin{remark}
    The Klein invariant does not depend on the choice of $\ldual{x}$ by naturality of $\operatorname{fl}_x$ and functoriality of $\eta$.
\end{remark}

\begin{example}
    \begin{itemize}
        \item If $\mathcal{C}$ is any braided monoidal 2-category, $I$ has a canonical duality datum $I \xrightarrow[]{\operatorname{id}} I$. The associated Klein invariant is $\kappa(\mathcal{C},I,\operatorname{id}) = 1 \in k$ in the base field. 
        \item If $\mathcal{C} = \Sigma \operatorname{sVect}, I$ has a second duality datum $e: I \to \ldual{I}$ for which $\kappa(\operatorname{sVect}, I, e) = -1$, where $e$ is a (choice of) the non-identity simple element in $\operatorname{sVec}$. 
    \end{itemize}
\end{example}

\subsection{The Tangle Hypothesis interpretation}

The name of the Klein invariant comes from the fact that it can be interpreted via the Tangle Hypothesis as the value on the Klein bottle of a certain TFT arising from the fully dualizable object $x \in \mathcal{C}$ with duality datum $r$. We recall the Tangle Hypothesis here. 

\begin{definition}
    Let $k \leq n$ be positive integers. The tangle category $\operatorname{Tang}_{k,n}$ is an $\mathbb{E}_{n-k}$-monoidal $(\infty,k)$-category whose $d$-morphisms are embedded $d$-dimensional submanifolds of $D^{n-k} \times [0,1]^d$ together with appropriate transversality boundary conditions on $D^{n-k} \times \{0,1\}^d$.
\end{definition}

If $\beta$ is a map $BG \to \operatorname{Gr}(n,k)$, the Grassmanian of $k$-planes in $n$-space, we also define $\operatorname{Tang}_{k,n}^\beta$ to be the same category as above, but where every manifold is equipped with a compatible $G$-structure on the tangent bundle. 

\begin{theorem}{(Tangle Hypothesis)}
    Let $\mathscr{C}$ be an $\mathbb{E}_{n-k}$-monoidal $(\infty,k)$-category. Evaluation at the point induces an equivalence of $(\infty,k)$-categories
   $$\operatorname{Fun}^\otimes(\operatorname{Tang}_{k,n}^\beta, \mathscr{C}) \simeq \left(\mathscr{C}^{\operatorname{f.d.}}\right)^{hG}$$
   where the right hand side is formed by taking homotopy fixed points of the $G$-action on the fully dualizable objects of $\mathscr{C}$.
\end{theorem}

If $\mathcal{C}$ is a braided monoidal 2-category, then $\mathcal{C}$ is an $\mathbb{E}_2$-algebra in $\operatorname{Cat}_2$. A fully dualizable object $x \in C$ with duality datum $r$ should give rise via the tangle hypothesis to a TFT 
$$Z_X: \operatorname{Tang}_{2,4}^\beta \to \mathcal{C}^{\operatorname{f.d.}} + \text{ duality datum } \simeq \left(\mathcal{C}^{\operatorname{f.d}}\right)^{h\mathbb{Z}}$$
for some tangential structure $\beta: B\mathbb{Z} \to \operatorname{Gr}(4,2)$ called a \emph{projective framing}, whose homotopy fixed points encode the duality datum. [T. Johnson-Freyd points out that he does not remember what circle in the Grassmanian this map picks out.] We can evaluate this TFT on various things: \\\\
\begin{tblr}{X[c,m] X[c,m] X[c,m] X[l,m]}
 $Z_X(\text{point})$ & ${\small \hspace{0.8cm} \input{circle.tikz}}$ \;\; & \hspace{0.8cm}$\longmapsto$&$ (x,r)$ \\\\
 {$Z_X(S^1)$ \\ \phantom{} \\ \phantom{} } &  ${\small \input{band.tikz}}$ \;\; & {\hspace{0.8cm} $\longmapsto$ \\ \phantom{} \\ \phantom{} }& {$\operatorname{tr}(\operatorname{id}_x) = \eta(x) \in \Omega \mathcal{C} $ \\ \phantom{} \\ \phantom{}}\\\\
   {$Z_X \left(\parbox{1.9cm}{\centering \footnotesize antipode map cobordism}\right)  $ \\ \phantom{} \\ \phantom{}}& ${\small \input{antipode.tikz}}$ \;\; & {\hspace{0.8cm}$\longmapsto$ \\ \phantom{} \\ \phantom{}} & {$\operatorname{fl}_x \cdot \eta(r)$ \\ \phantom{} \\ \phantom{}} \\\\
   {\phantom{} \\ $Z_X(\text{Klein})$ } & ${\small \input{klein.tikz}}$ &  { \phantom{} \\ \ \hspace{0.8cm} $\longmapsto$  } & { \phantom{} \\ $ \operatorname{tr}(\operatorname{fl}_x \cdot \eta(r)) \in \Omega^2\mathcal{C}$ \\} 
\end{tblr}

\section{Distinguishing $\mathcal{S}$ and $\mathcal{T}$}

We can now use the Klein invariant to distinguish $\mathcal{S}$ and $\mathcal{T}$. We will use the fact that $\mathcal{S}$ and $\mathcal{T}$ admit certain braided fusion sub-2-categories.
\begin{itemize}
    \item $2 \operatorname{Vect}^\text{triv}[\mathbb{Z}/2]$ is a sub-2-category of $\mathcal{S}$ via $(M,r) \mapsto (I, \operatorname{id})$ (here $M$ is a choice of invertible element in the magnetic component). We have that 
    $$\kappa(\mathcal{S},M,r) = \kappa(\Sigma \operatorname{sVect}, I, \operatorname{id}) = 1$$
    since $\kappa(\mathcal{S},M,r)$ only depends on where $M$ lives. 
    \item $2\operatorname{Vect}^\tau[\mathbb{Z}/2]$ includes in $\mathcal{T}$ via $(M,r) \mapsto (I,e)$. Here $\tau$ is the twisting cocycle computed in Cameron's talk. [T. Johnson-Freyd points out that picking a duality datum $r$ here is equivalent to picking a product decomposition $B^2 \mathcal{G} \simeq \kappa(\mathbb{Z}/2,2) \times \kappa(\mathbb{Z}/2,3)$]. We have that 
    $$\kappa(\mathcal{T}, M, r) = \kappa(\Sigma \text{sVect}, I, e) = -1.$$
\end{itemize}

But this does not prove that $\mathcal{S}$ and $\mathcal{T}$ are different - for that, we need to show that the above values are independent of the choice of magnetic object $M$ and duality datum $r$. Up to isomorphism, the only other choice of duality datum is $e \boxtimes r$. We compute
$$\kappa(\mathcal{S}\text{ or } \mathcal{T}, M, e \boxtimes r) = \operatorname{tr}(\operatorname{fl} \cdot \eta(e \boxtimes r)) = \operatorname{tr}(\operatorname{fl} \cdot \eta(r) \cdot \eta(e \boxtimes \operatorname{id}_M)) = \kappa(\mathcal{S}\text{ or } \mathcal{T}, M, r)$$
The last equality follows from a cancellation of two $-1$s. So there is no dependence on $r$. The only other magnetic object is $C \boxtimes M$ where $C$ is the Clifford algebra. To show that there is no dependence on this choice, note that by Nivedita's talk there exists a morphism  $a: M \to C \boxtimes M$ such that $\eta(a)$ is an isomorphism. ($a$ is the unique simple 1-morphism). We have a commutative diagram
\[\begin{tikzcd}
	{\eta(M) } && {\eta(^\vee M)} && {\eta(M)} \\
	&& {} \\
	{\eta(C \boxtimes M) } && {\eta(^\vee (C \boxtimes M)) } && {\eta(C \boxtimes M)}
	\arrow["{\eta(r)}", Rightarrow, from=1-1, to=1-3]
	\arrow["{\eta(a)}"', Rightarrow, from=1-1, to=3-1]
	\arrow["{\operatorname{fl}_M}", Rightarrow, from=1-3, to=1-5]
	\arrow["{\eta(^{*\vee}a)}"', Rightarrow, from=1-3, to=3-3]
	\arrow["{\eta(a)}", Rightarrow, from=1-5, to=3-5]
	\arrow["{\eta(s)}"', Rightarrow, from=3-1, to=3-3]
	\arrow["{\operatorname{fl}_{C \boxtimes M}}"', Rightarrow, from=3-3, to=3-5]
\end{tikzcd}\]

Here $r$ and $s$ are duality data. The left square commutes since $s \circ a $ and $^{*\vee}a \circ r$ are isomorphic since they are both simple, combined with functoriality of $\eta$. The right square commutes by naturality of $\operatorname{fl}$. It follows that the top is conjugate to the bottom, and so the traces of both agree, which shows that 
$$\kappa(\mathcal{S}\text{ or } \mathcal{T}, M) = \kappa(\mathcal{S}\text{ or } \mathcal{T}, C \boxtimes M)$$
and so there is no dependence on $M$.

\begin{cor}
    $\mathcal{S}$ is not equivalent to $\mathcal{T}$.
\end{cor}
\begin{proof}
    If they were, an equivalence would have to send $M$ to either $M$ or $C \boxtimes M$. But in both cases the Klein invariant disagrees. 
\end{proof}

This in fact makes our life harder, since now it's possible that $\mathcal{Z}(\Sigma \mathcal{B})$ is not independent of $\mathcal{B}$. We will start working towards ruling out one of $\mathcal{S}$ and $\mathcal{T}$ as a possibility for $\mathcal{Z}(\Sigma \mathcal{B})$.

\begin{prop}
    Let $X \in \mathcal{R} = (\mathcal{S}\text{ or } \mathcal{T})$ be purely magnetic, i.e. one of the following equivalent conditions hold: 
    \begin{itemize}
        \item The full braiding of $X$ with $e$ is $- \operatorname{id},$ or
        \item $X = \sum_{i} X_i, $ with each $x_i$ simple and magnetic.
    \end{itemize}
  Let $r: X \to^\vee X$ be a self duality datum. Then: 
  \begin{itemize}
      \item If $\mathcal{R} = \mathcal{S}$ then $\kappa(\mathcal{S}, X, r) \in \mathbb{Z}_{\geq 0},$
      \item If $\mathcal{R} = \mathcal{T}$ then $\kappa(\mathcal{T}, X, r) \in \mathbb{Z}_{\leq 0}.$
  \end{itemize}
\end{prop}

\begin{proof}
    Pick auxiliary self duality datums $r_i: x_i \to^\vee x_i$ for all $x_i$. (Note that all $x_i$ are either $M$ or $C \boxtimes M$). These sum together to another duality datum 
    $$q := \sum_i r_i: X \xrightarrow[]{\sim}^\vee X.$$
    We obtain an automorphism $q^{-1}r$ of $X$ which when written at the level of the sum decomposition of $X$ turns out to be of the form 
    $$q^{-1}r =  \begin{bmatrix}
    \lambda_{1} & & \\
    & \ddots & \\
    & & \lambda_{n}
  \end{bmatrix} P : \sum_{i} X_i \to \sum_{i} X_i$$
  where $P$ is a permutation matrix of the summands, and $\lambda_i$ are invertible morphisms. [D. Reutter explains that this is specific to semisimple 2-categories, and you would not encounter this phenomenon in semisimple 1-categories. T. Johnson-Freyd adds that the formula for the inverse of a matrix has some minus signs in it, and the only way to avoid them is if you're a permutation matrix.] Thus
  $$\eta(q^{-1})\eta(r) =   \begin{bmatrix}
    \eta(\lambda_{1}) & & \\
    & \ddots & \\
    & & \eta(\lambda_{n})
  \end{bmatrix} \eta(P) \simeq \eta(P)$$
 since $\eta(\text{automorphism of simple})$ is trivial in $\mathcal{R}$. By functoriality of $\eta$ we obtain that $\eta(r) = \eta(qP)$ and so 
 \begin{align*}
     \kappa(\mathcal{R}, X,r) &= \operatorname{tr}(\operatorname{fl}_X \cdot \eta(r)) \\
            &= \operatorname{tr}(\operatorname{fl}_X \cdot \eta(qP)) \\
            &= \sum_{j \text{ fixed by } P } \operatorname{tr}(\operatorname{fl}_{X_j} \eta(r_j)) \\
            &= \sum_{j \text{ fixed by } P } \kappa(\mathcal{R}, X_j, r_j)
 \end{align*}
 Since $X_i$ is always $M$ or $C \boxtimes M$, and since the Klein invariant is either always positive or always negative for those two objects depending on whether $\mathcal{R} = \mathcal{S}$ or $\mathcal{R} = \mathcal{T}$, the result is immediate from the above description. 
\end{proof}

\section{An extension of the Klein invariant}

There exists an useful extension of the Klein invariant, denoted by $\uK$. To define it, we first define an extension of $\eta$, denoted $$\ue(x) = \operatorname{Hom}_{\Omega \mathcal{C} }(\eta(x), I).$$ We also denote by $\underline{\operatorname{fl}}_X: \ue(x) \to \ue(^\vee x )$ the map $\phi \mapsto \phi \circ \operatorname{fl}_X$. If $r$ is a self duality datum of $x$, we define $\ue(r) : \ue(^\vee x) \to \ue(x)$ by $\phi \mapsto \phi \cdot \eta(r)$.

\begin{definition}
    The extended Klein invariant $\uK$ is given by 
    $$\uK(\mathcal{C}, x, r) := \operatorname{tr}(\ue(r) \underline{\operatorname{fl}}_x).$$
\end{definition}

\begin{prop}
    If $X$ is a magnetic object in $\mathcal{C}$ with self duality datum $r$, then $\kappa(\mathcal{C}, x, r) = \uK(\mathcal{C},x,r).$
\end{prop}

[T. Johnson-Freyd remarks that $\eta(\text{Clifford})$ is an electron, whereas $\ue(\text{Clifford})$ is $0$. In particular, $\kappa(\text{Clifford})$ is not well-defined, whereas $\uK(\text{Clifford})$ does not care, so we defined $\uK$ to avoid the Clifford object. D. Reutter further remarks that this allows one to do computations purely with magnetic objects.]

%% file: talks/4.3/lol.tikzstyles

\tikzstyle{new style 0}=[fill=white, draw=black, shape=circle, tikzit draw=black, tikzit shape=circle]
\tikzstyle{new style 1}=[fill=black, draw=black, shape=circle, tikzit fill=black, tikzit draw=black, tikzit shape=circle]

\tikzstyle{new edge style 0}=[-]

%% file: talks/4.3/style.tikzstyles


\tikzstyle{new edge style 0}=[-, draw={rgb,255: red,128; green,128; blue,128}, dashed, dash pattern=on 2mm off 2mm, fill=white, tikzit draw={rgb,255: red,128; green,128; blue,128}]
\tikzstyle{new edge style 1}=[-, draw={rgb,255: red,128; green,128; blue,128}, fill=none, tikzit draw={rgb,255: red,128; green,128; blue,128}]
\tikzstyle{new edge style 2}=[-, tikzit draw=black, fill=white, dashed, dash pattern=on 0.5mm off 0.5mm]

%% file: talks/4.3/circle.tikz
\begin{tikzpicture}
	\begin{pgfonlayer}{nodelayer}
		\node [style=none] (0) at (0, -2) {};
		\node [style=none] (1) at (-2, 0) {};
		\node [style=none] (2) at (0, 2) {};
		\node [style=none] (3) at (2, 0) {};
		\node [style=new style 1] (4) at (0, 0) {};
		\node [style=none] (5) at (-2.75, 0) {$D^2$};
	\end{pgfonlayer}
	\begin{pgfonlayer}{edgelayer}
		\draw (0.center)
			 to [bend left=45] (1.center)
			 to [bend left=45] (2.center)
			 to [bend left=45] (3.center)
			 to [bend left=45] cycle;
	\end{pgfonlayer}
\end{tikzpicture}

%% file: talks/4.3/band.tikz
\begin{tikzpicture}
	\begin{pgfonlayer}{nodelayer}
		\node [style=none] (0) at (0, 1) {};
		\node [style=none] (1) at (0, -2) {};
		\node [style=none] (2) at (6, 1) {};
		\node [style=none] (3) at (6, -2) {};
		\node [style=none] (4) at (1.25, -0.5) {};
		\node [style=none] (5) at (3, -0.5) {};
		\node [style=none] (6) at (4.75, -0.5) {};
		\node [style=none] (7) at (2.75, -1) {};
		\node [style=none] (8) at (3.25, 0) {};
		\node [style=none] (9) at (-1.25, -0.5) {};
		\node [style=none] (10) at (-1.25, -0.5) {};
		\node [style=none] (11) at (-1.5, -0.5) {$D^2$};
		\node [style=none] (12) at (3, -2.5) {$[0,1]$};
	\end{pgfonlayer}
	\begin{pgfonlayer}{edgelayer}
		\draw [style=new edge style 1, bend left=90] (0.center) to (1.center);
		\draw [style=new edge style 1, bend right=90] (0.center) to (1.center);
		\draw [style=new edge style 1, bend left=90] (2.center) to (3.center);
		\draw [style=new edge style 0, bend left=270] (2.center) to (3.center);
		\draw [style=new edge style 1] (0.center) to (2.center);
		\draw [style=new edge style 1] (1.center) to (3.center);
		\draw [bend left=90, looseness=1.75] (4.center) to (5.center);
		\draw [bend right=90, looseness=1.75] (5.center) to (6.center);
		\draw [bend right=75, looseness=1.25] (4.center) to (7.center);
		\draw [bend left=75, looseness=1.25] (8.center) to (6.center);
	\end{pgfonlayer}
\end{tikzpicture}

%% file: talks/4.3/antipode.tikz
\begin{tikzpicture}
	\begin{pgfonlayer}{nodelayer}
		\node [style=none] (0) at (0, 1) {};
		\node [style=none] (1) at (0, -2) {};
		\node [style=none] (2) at (6, 1) {};
		\node [style=none] (3) at (6, -2) {};
		\node [style=none] (9) at (-1.25, -0.5) {};
		\node [style=none] (11) at (-1.5, -0.5) {$D^2$};
		\node [style=none] (12) at (3, -2.5) {$[0,1]^2$};
		\node [style=none] (14) at (2, 0.5) {};
		\node [style=none] (15) at (1.5, -0.25) {};
		\node [style=none] (16) at (1.75, -1.5) {};
		\node [style=none] (17) at (2, -0.75) {};
		\node [style=none] (18) at (1.75, -0.5) {};
		\node [style=none] (19) at (4.25, 0) {};
		\node [style=none] (20) at (4.25, -1) {};
	\end{pgfonlayer}
	\begin{pgfonlayer}{edgelayer}
		\draw [style=new edge style 1, bend left=90] (0.center) to (1.center);
		\draw [style=new edge style 1, bend right=90] (0.center) to (1.center);
		\draw [style=new edge style 1, bend left=90] (2.center) to (3.center);
		\draw [style=new edge style 0, bend left=270] (2.center) to (3.center);
		\draw [style=new edge style 1] (0.center) to (2.center);
		\draw [style=new edge style 1] (1.center) to (3.center);
		\draw [bend left=75, looseness=1.25] (15.center) to (14.center);
		\draw [bend left=75, looseness=1.25] (17.center) to (16.center);
		\draw [bend left=90, looseness=1.75] (16.center) to (18.center);
		\draw [bend right=90, looseness=1.75] (18.center) to (14.center);
		\draw [bend left=90, looseness=1.25] (19.center) to (20.center);
		\draw [bend left=270, looseness=1.25] (19.center) to (20.center);
		\draw (14.center) to (19.center);
		\draw (16.center) to (20.center);
	\end{pgfonlayer}
\end{tikzpicture}

%% file: talks/4.3/klein.tikz
\begin{tikzpicture}
	\begin{pgfonlayer}{nodelayer}
		\node [style=none] (0) at (0, 1) {};
		\node [style=none] (1) at (0, -2) {};
		\node [style=none] (2) at (6, 1) {};
		\node [style=none] (3) at (6, -2) {};
		\node [style=none] (9) at (-1.25, -0.5) {};
		\node [style=none] (11) at (-1.5, -0.5) {$D^2$};
		\node [style=none] (12) at (3, -2.5) {$[0,1]^2$};
		\node [style=none] (13) at (1.75, 0) {};
		\node [style=none] (14) at (1.75, -1.25) {};
		\node [style=none] (15) at (4.75, 0.25) {};
		\node [style=none] (16) at (3, -0.25) {};
		\node [style=none] (17) at (4.25, 0) {};
		\node [style=none] (18) at (2.25, -1) {};
		\node [style=none] (19) at (2.75, -0.75) {};
		\node [style=none] (20) at (2, -0.25) {};
		\node [style=none] (21) at (1.5, -1) {};
		\node [style=none] (23) at (2.5, -0.5) {};
		\node [style=none] (24) at (3, -0.5) {};
	\end{pgfonlayer}
	\begin{pgfonlayer}{edgelayer}
		\draw [style=new edge style 1, bend left=90] (0.center) to (1.center);
		\draw [style=new edge style 1, bend right=90] (0.center) to (1.center);
		\draw [style=new edge style 1, bend left=90] (2.center) to (3.center);
		\draw [style=new edge style 0, bend left=270] (2.center) to (3.center);
		\draw [style=new edge style 1] (0.center) to (2.center);
		\draw [style=new edge style 1] (1.center) to (3.center);
		\draw [bend right=60] (14.center) to (15.center);
		\draw [in=120, out=30, looseness=1.25] (13.center) to (16.center);
		\draw [in=90, out=45, looseness=1.25] (19.center) to (17.center);
		\draw [style=new edge style 2] (19.center)
			 to [in=0, out=-135] (18.center)
			 to [in=45, out=-180, looseness=1.25] (14.center);
		\draw [in=45, out=90, looseness=1.25] (15.center) to (16.center);
		\draw [style=new edge style 2, in=135, out=-30, looseness=1.50] (13.center) to (20.center);
		\draw [in=135, out=-150] (13.center) to (14.center);
		\draw [bend left] (13.center) to (14.center);
		\draw [bend right=15, looseness=1.25] (23.center) to (16.center);
		\draw [style=new edge style 2, in=-180, out=-45] (20.center) to (23.center);
		\draw [in=-90, out=-90, looseness=1.75] (24.center) to (17.center);
	\end{pgfonlayer}
\end{tikzpicture}

%% file: talks/5.1/main.tex
\def\dist{1.0}

Talk by Daniel Teixeira, notes by Adrien DeLazzer Munier.

In this talk, we put everything together to finish the proof of the Main Theorem. We'll do a few computations in the graphical calculus, so we start by recalling some definitions and showing their graphical representations. If $\mathcal{B}$ is our braided fusion 1-category, then we denote its braiding as $\beta_{x,y} \colon x \otimes y \rightarrow y \otimes x$. Graphically, we'll represent the braiding and its inverse by
\[
\begin{tikzpicture}

    \node[draw] {
    \begin{tikzpicture}
    
        \node[draw] (left) at (0,0) {
        \includegraphics[scale=.25]{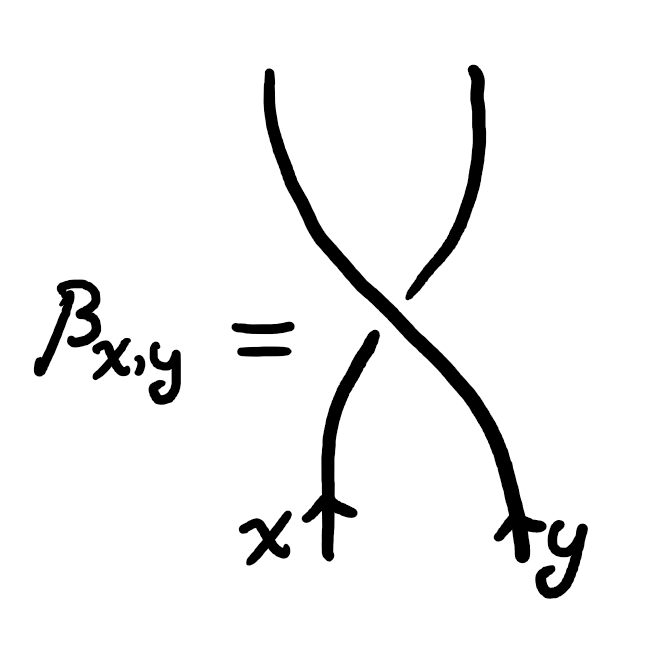}
        };
        
        \node[draw, anchor=west] at ($(left.east)+(\dist,0)$) {
        \includegraphics[scale=.25]{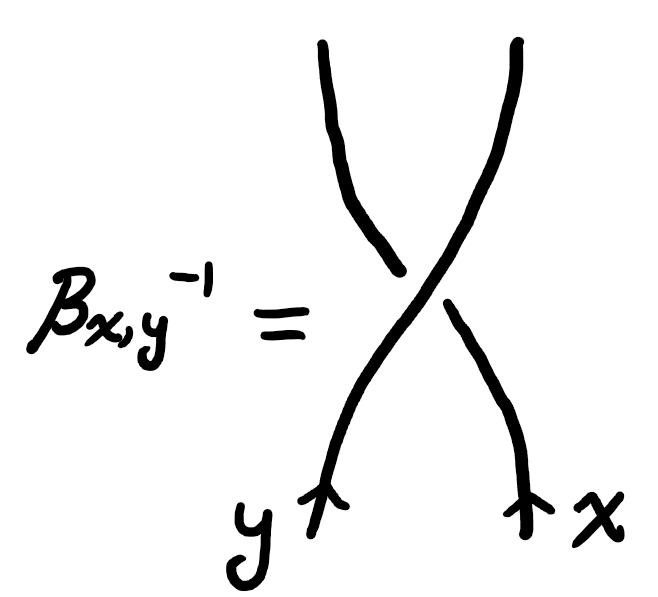}
        };
    
    \end{tikzpicture}
    };
  
\end{tikzpicture}
\]
These satisfy the hexagon axioms, which are graphically represented as
\[
\begin{tikzpicture}

    \node[draw] {
    \begin{tikzpicture}
    
        \node[draw] (left) at (0,0) {
        \includegraphics[scale=.25]{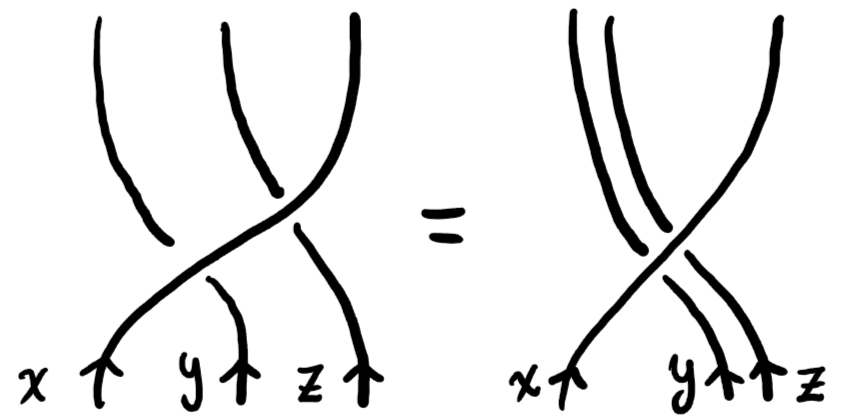}
        };
        
        \node[draw, anchor=west] at ($(left.east)+(\dist,0)$) {
        \includegraphics[scale=.25]{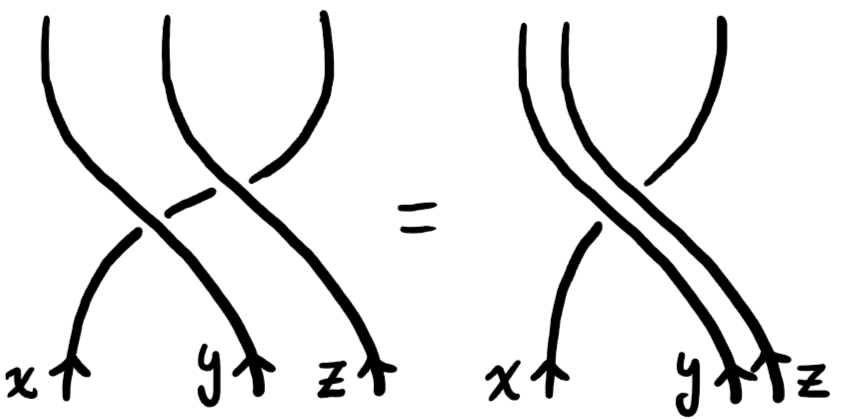}
        };
    
    \end{tikzpicture}
    };
  
\end{tikzpicture}
\]
In order to make the graphical calculus more legible, we will often implicitly make use of the hexagon axioms to group braided strands together. Let $(A,\gamma) \in \operatorname{HBA}(\mathcal{B})$ be a half-braided algebra. We'll draw the multiplication without any additional marking, while the half-braid will be marked with a red dot.
\[
\begin{tikzpicture}

    \node[draw] {
    \begin{tikzpicture}
    
        \node[draw] (left) at (0,0) {
        \includegraphics[scale=.25]{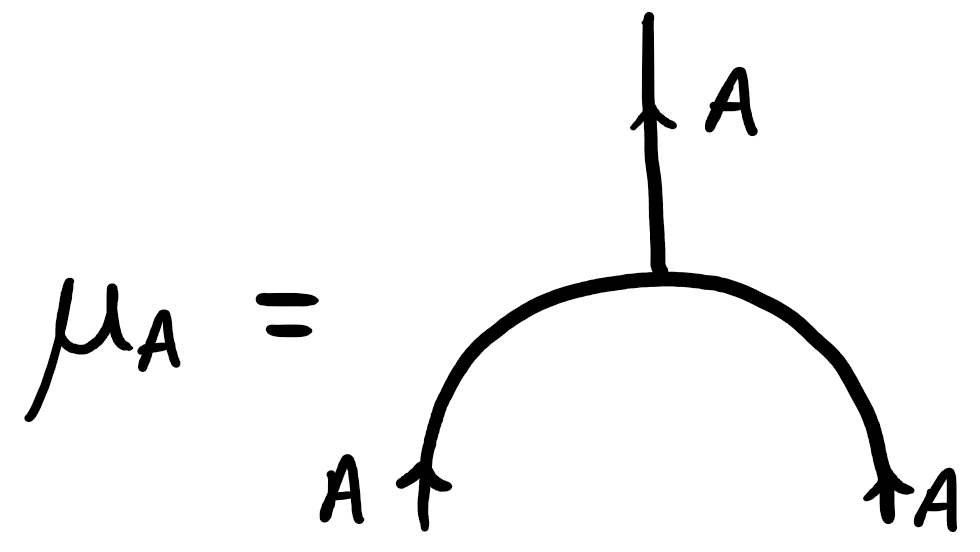}
        };
        
        \node[draw, anchor=west] at ($(left.east)+(\dist,0)$) {
        \includegraphics[scale=.25]{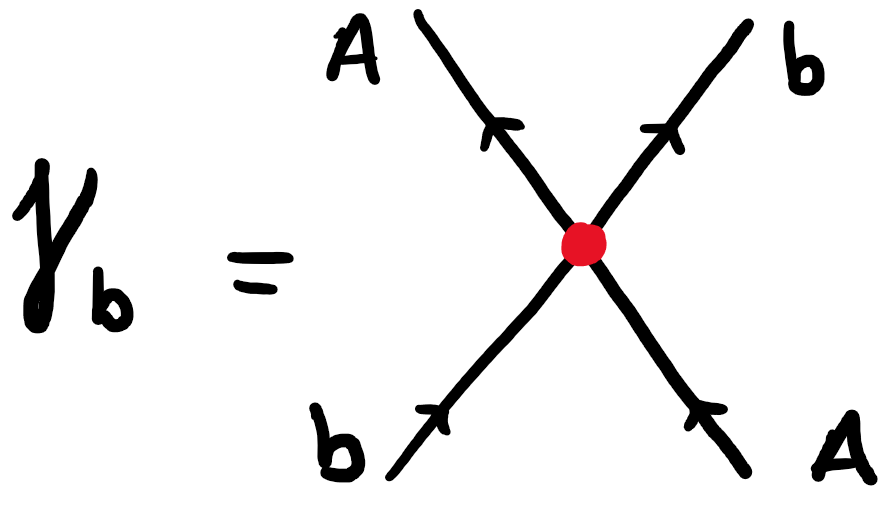}
        };
    
    \end{tikzpicture}
    };
  
\end{tikzpicture}
\]
The map $\gamma_b \colon b \otimes A \rightarrow A \otimes b$ is natural in $b$, and satisfies an analogous hexagon axiom. The naturality means that, for any morphism $f \colon b \rightarrow c$ in $\mathcal{B}$, we have
\[
\begin{tikzpicture}
    \node[draw] {
    \includegraphics[scale=0.25]{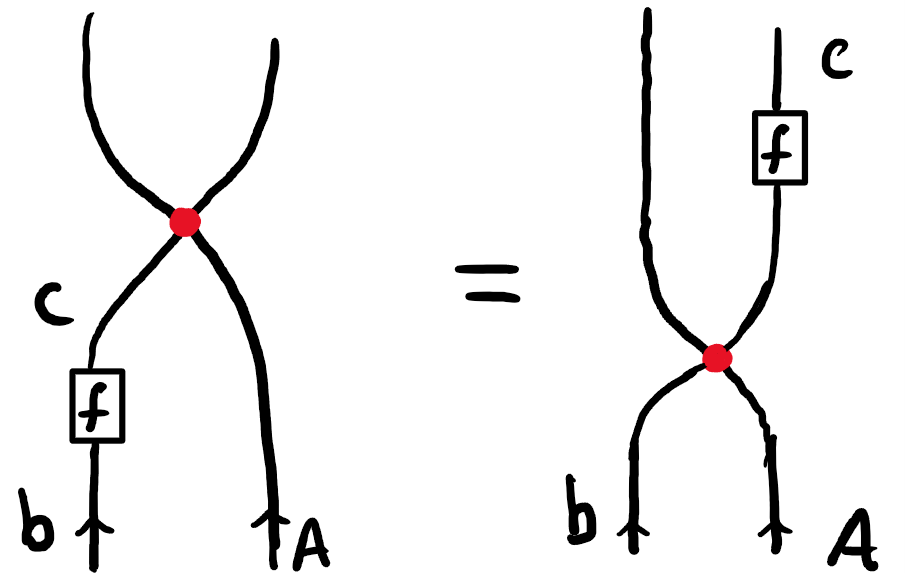}
    };
\end{tikzpicture}
\]
And finally, there is the compatibility condition between $\gamma, \mu_A,$ and $\beta$ that makes $(A, \mu_A, \gamma)$ into an actual half-braided algebra, namely
\[
\begin{tikzpicture}
    \node[draw] {
    \includegraphics[scale=0.25]{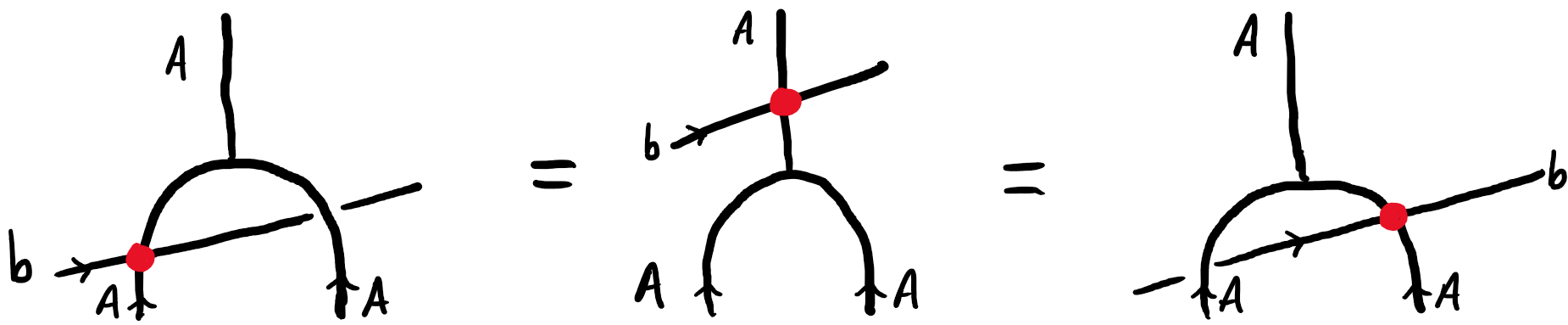}
    };
\end{tikzpicture}
\]
Any half-braided algebra discussed is assumed to have a unit, though we often suppress it from the notation. Associated to each half-braided algebra $(A,\gamma)$ is a morphism $\theta_{(A,\gamma)} \colon A \rightarrow A$ in $\mathcal{B}$, called its \textit{twist}, that is graphically defined as follows.
\[
\begin{tikzpicture}
    \node[draw] {
    \includegraphics[scale=0.25]{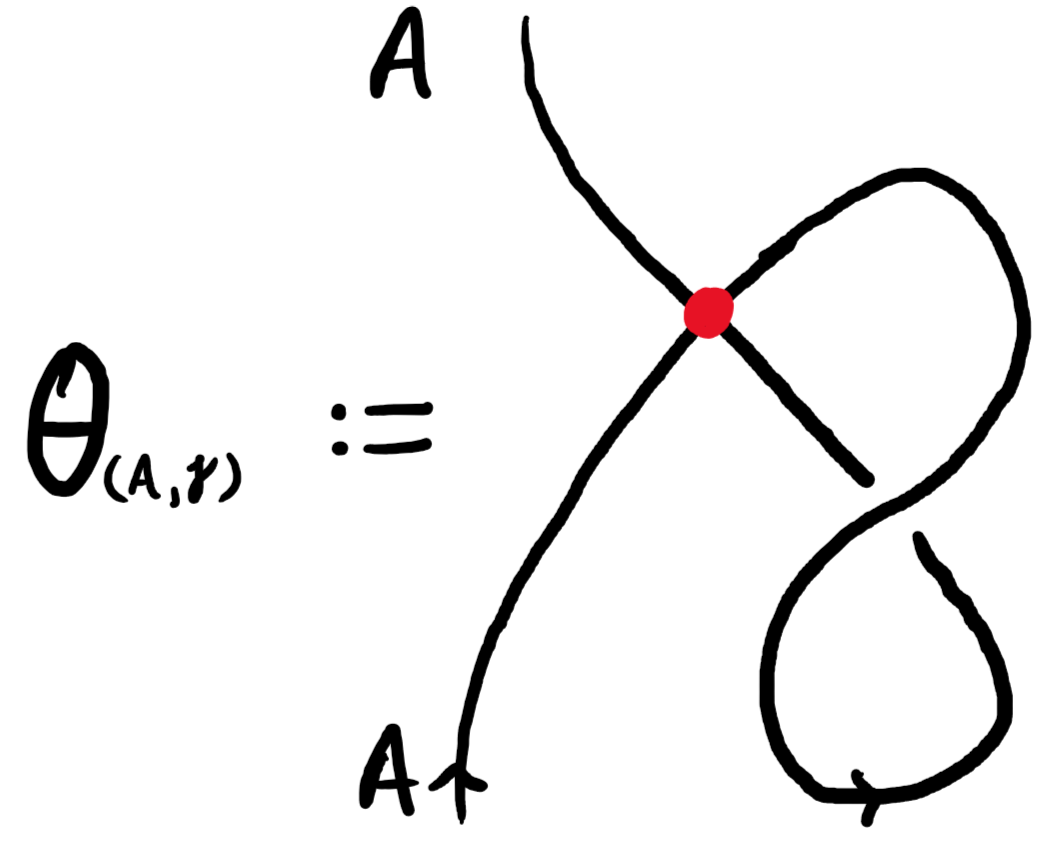}
    };
\end{tikzpicture}
\]
We use this twist to define a notion of trace of half-braided algebras. A morphism $\varepsilon \colon A \rightarrow \mathds{1}$ is called a \textit{twisted trace} if the following equation holds in $\mathcal{B}$.
\[
\begin{tikzpicture}
    \node[draw] {
    \includegraphics[scale=0.25]{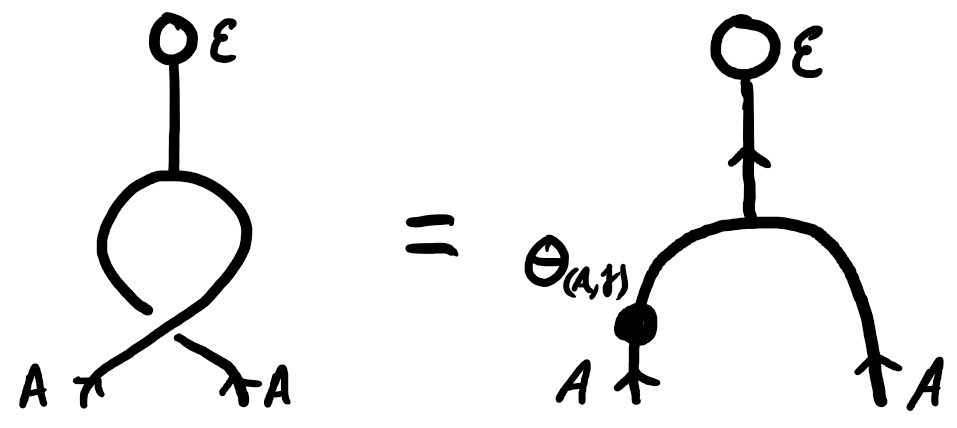}
    };
\end{tikzpicture}
\]
Let $\mathcal{E}_A$ be the coequalizer of the two morphisms that get composed with $\varepsilon$ in this definition, i.e. $\mu_A \cdot \beta_{A,A}^{-1}$ and $(\theta_{(A,\gamma)} \otimes \operatorname{id}_A)$. By the coequalizer universal property, a twisted trace $\varepsilon \colon A \rightarrow \mathds{1}$ is then equivalent to a map $\mathcal{E}_A \rightarrow \mathds{1}$. This is why the authors in \cite{JFR} call $\mathcal{E}_A$ the ``universal trace". We can rewrite the right hand side of the above equation in a more convenient way as follows.
\[
\begin{tikzpicture}
    \node[draw] {
    \includegraphics[scale=0.25]{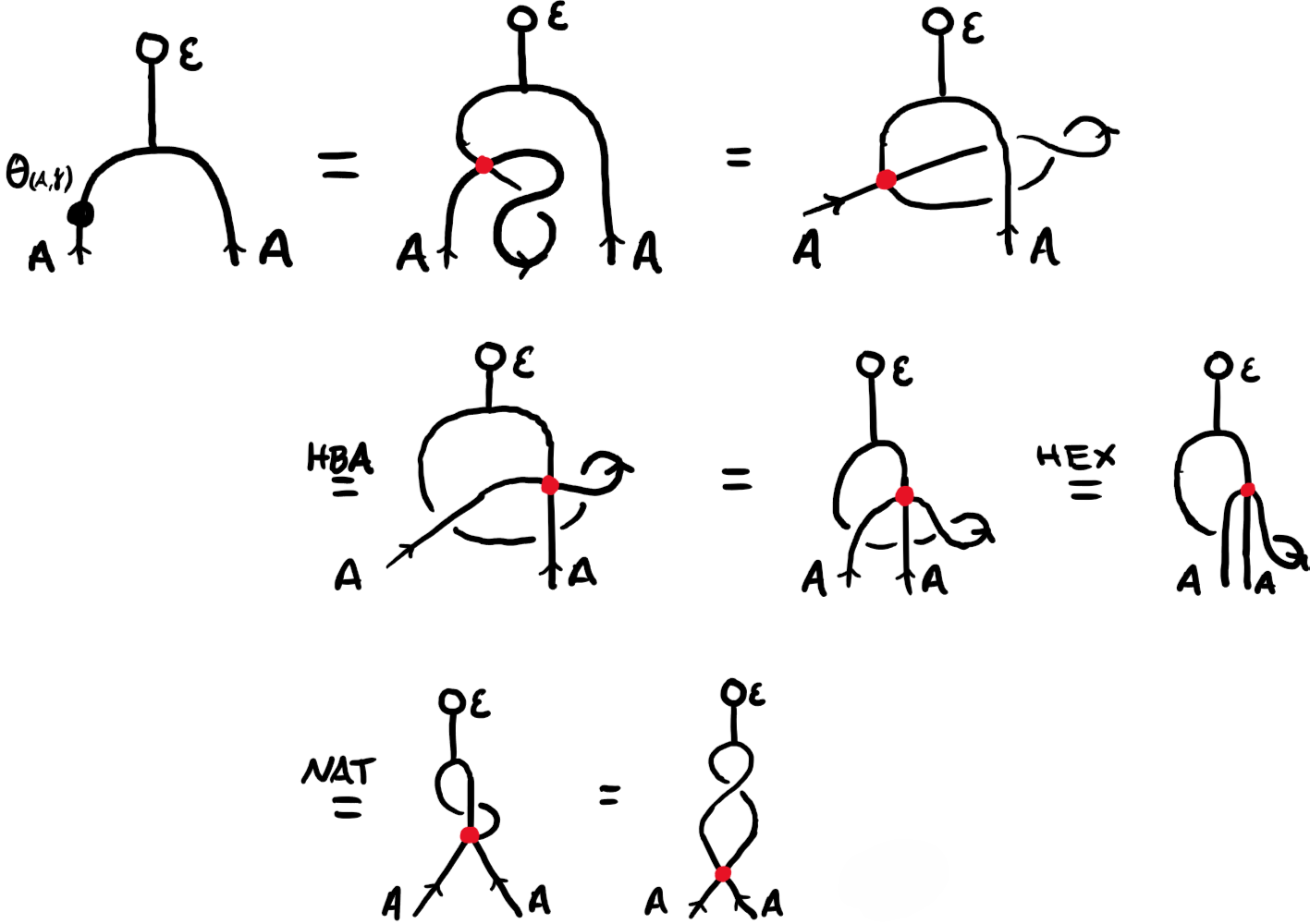}
    };
\end{tikzpicture}
\]
In particular, we get that $\mu_A \cdot (\theta_{(A,\gamma)} \otimes \operatorname{id}_A) = \mu_A \cdot \beta_{A,A}^{-1} \cdot \gamma_A$, which lets us rewrite our universal trace as
\[
\mathcal{E}_A = \operatorname{coeq}\left(
\begin{tikzcd}A \otimes A \arrow[rr, shift left = 1, "\mu_A\cdot \beta_{A,A}^{-1} \cdot \gamma_A"]\arrow[rr, shift right=1, "\mu_A \cdot \beta_{A,A}^{-1}"'] && A
\end{tikzcd}\right)
\]

Now lets get to the main story. In particular, from now on we'll consider separable half-braided algebras since $\operatorname{sHBA}(\mathcal{B}) \simeq \mathcal{Z}_1(\Sigma \mathcal{B})$. We want to be able to compute the Klein bottle invariant $\kappa$ on dualizable objects in some model of $\mathcal{Z}_1(\Sigma\mathcal{B})$, and it happens to not be too hard with half-braided algebras. If we have a half-braided algebra $(A,\gamma)\in \operatorname{sHBA}(\mathcal{B})$, and a self-duality isomorphism $\phi \colon A \rightarrow ^{*\hspace{-6pt}}A$, the invariant is given by
\[
\begin{tikzpicture}
    \node[draw] {
    \includegraphics[scale=0.25]{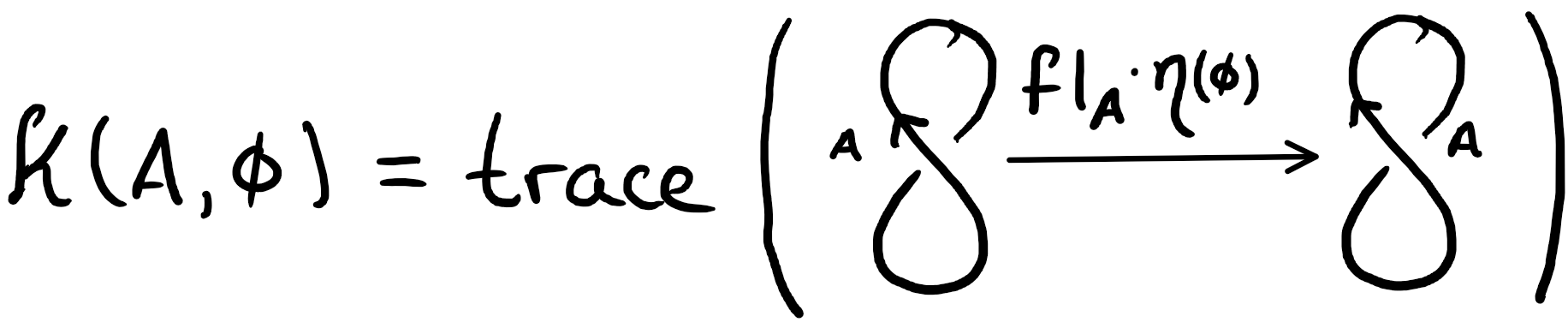}
    };
\end{tikzpicture}
\]
Now $\kappa$ is defined in terms of a trace, so maybe this suggests we should look at $\mathcal{E}_A$. This turns out to be useful, but not in the seemingly obvious way. Rather, Proposition 3.33 of \cite{JFR} implies we have the following isomorphism.
\begin{equation}
\label{eq:TraceIso}
    \mathcal{E}_A \cong \eta(A,\gamma)
\end{equation}
\noindent This way, the map $\operatorname{fl}_A \cdot \eta(\phi) \colon \mathcal{E}_A \rightarrow \mathcal{E}_A$ will be a ``transformation" of traces. More precisely, the induced map $\underline{\eta}(\phi) \cdot \underline{\operatorname{fl}}_A$ obtained by applying $\operatorname{Hom}_{\mathcal{B}}(-,\mathds{1})$ will be an automorphism of the vector space of twisted traces $\underline{\eta}(A,\gamma) \coloneqq \operatorname{Hom}_{\mathcal{B}}(\eta(A,\gamma),\mathds{1})$. $\kappa$ is then computed as the trace of this automorphism of traces. 

To obtain the isomorphism in (\ref{eq:TraceIso}), we need duality data for half-braided algebras. In Lemma 3.31 of \cite{JFR}, the authors show that these are given by the two ``opposite" algebras with multiplications given by composing $\mu_A$ with $\beta_{A,A}$ giving $\leftindex^{\text{op}}A = \leftindex^{*}A$ or with $\beta^{-1}_{A,A}$ giving $A^{\text{op}} = A^*$. There is a similar story for their half-braidings. Next, we make use of Corollary 3.39 of \cite{JFR}. This states that if $\phi \colon (A, \gamma) \rightarrow (\leftindex^{\text{op}}{A},\leftindex^{\text{op}}{\gamma})$ is an isomorphism of half-braided algebras, then the following diagram commutes in $\mathcal{B}$.
\[
\begin{tikzcd}
	A && A \\
	{\eta(A,\gamma)} && {\eta(A,\gamma)}
	\arrow["\phi", from=1-1, to=1-3]
	\arrow[two heads, from=1-1, to=2-1]
	\arrow[two heads, from=1-3, to=2-3]
	\arrow["{\operatorname{fl}_{A}\cdot\eta(\phi)}", from=2-1, to=2-3]
\end{tikzcd}
\]
In particular, applying $\operatorname{Hom}_{\mathcal{B}}(-,\mathds{1}) \colon \mathcal{B} \rightarrow \operatorname{Vec}$ to this diagram gives the following formula
\begin{align}
\label{eq:TraceTrace}
\begin{split}
    \underline{\eta}(\phi) \cdot \underline{\operatorname{fl}}_A \colon \underline{\eta}(A,\gamma) &\longrightarrow \underline{\eta}(A,\gamma)\\
    \varepsilon &\longmapsto \varepsilon \cdot \phi
\end{split}
\end{align}
And then $\kappa(A,\gamma)$ is the trace of this linear map.

Now that we know how to compute $\kappa(A,\gamma)$ for any fully dualizable object of $\operatorname{sHBA}(\mathcal{B})$, we can apply this to our favourite: the canonical Lagrangian half-braided algebra. This represents the distinguished Lagrangian object $\mathcal{L} \in \mathcal{Z}_1(\Sigma \mathcal{B})$ as a (separable) half-braided algebra. The underlying object is given by the coend 
\[
L = \int^{b\in \mathcal{B}}b\otimes b^*
\]
and so has some very nice properties. Specifically, $L$ is the universal object of $\mathcal{B}$ equipped with a dinatural transformation $\iota_b \colon b \otimes b^* \rightarrow L$. Thus any dinatural family of maps $\alpha_b \colon b \otimes b^* \rightarrow x$ factors uniquely through $\iota_b$. This means that to define a map $\alpha \colon L {\rightarrow} x$, it suffices to define components of a dinatural transformation, i.e. $\alpha_b$. If $f \colon b \rightarrow c$ is any morphism in $\mathcal{B}$, the dinaturality of $\alpha_b$ in $b$ amounts to the following equation
\[
\begin{tikzpicture}
    \node[draw] {
    \includegraphics[scale=0.25]{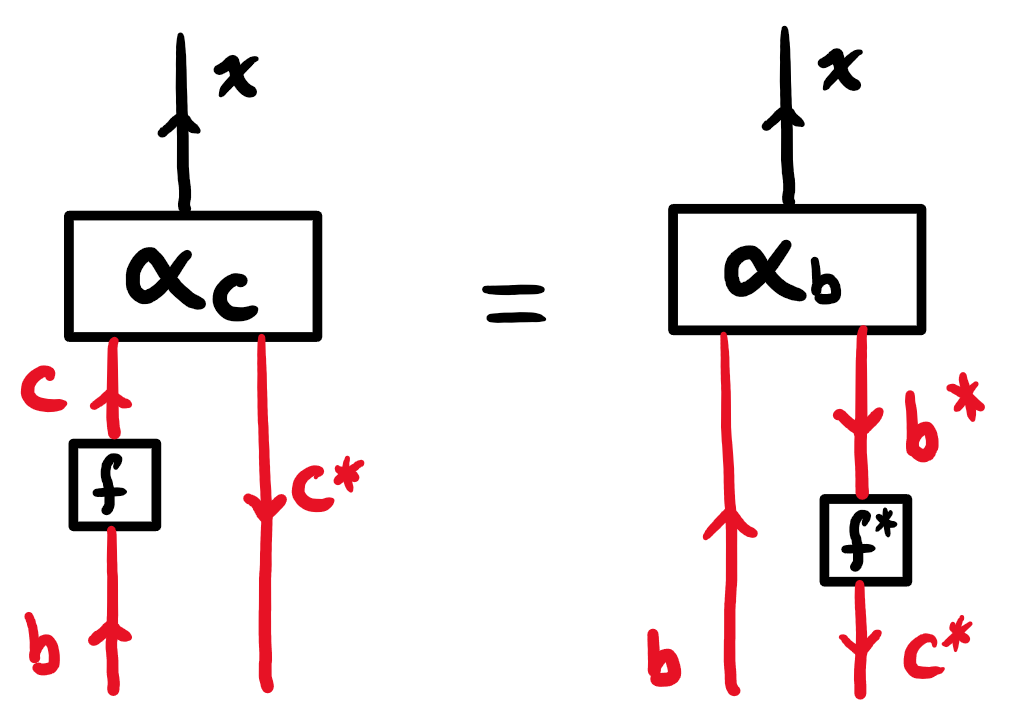}
    };
\end{tikzpicture}
\]
Going from the left hand side to the right, this feels like we're pushing the $c$ string through the $\alpha_c$ box, where it the comes out the other end, but dualized. When a morphism is attached to this $c$ string, it follows along for the ride and also gets dualized, but now the string hitting the box is labeled by the source of the morphism and so we change the label on the box to $\alpha_b$. As with the JFR paper, we use the convention that strings coloured in red are those that are dinatural. The unit of $L$ is given by $\iota_{\mathds{1}} \colon \mathds{1} \cong \mathds{1} \otimes \mathds{1}^* \rightarrow L$ and its product is given by the following dinatural map
\[
\begin{tikzpicture}
    \node[draw] {
    \includegraphics[scale=0.25]{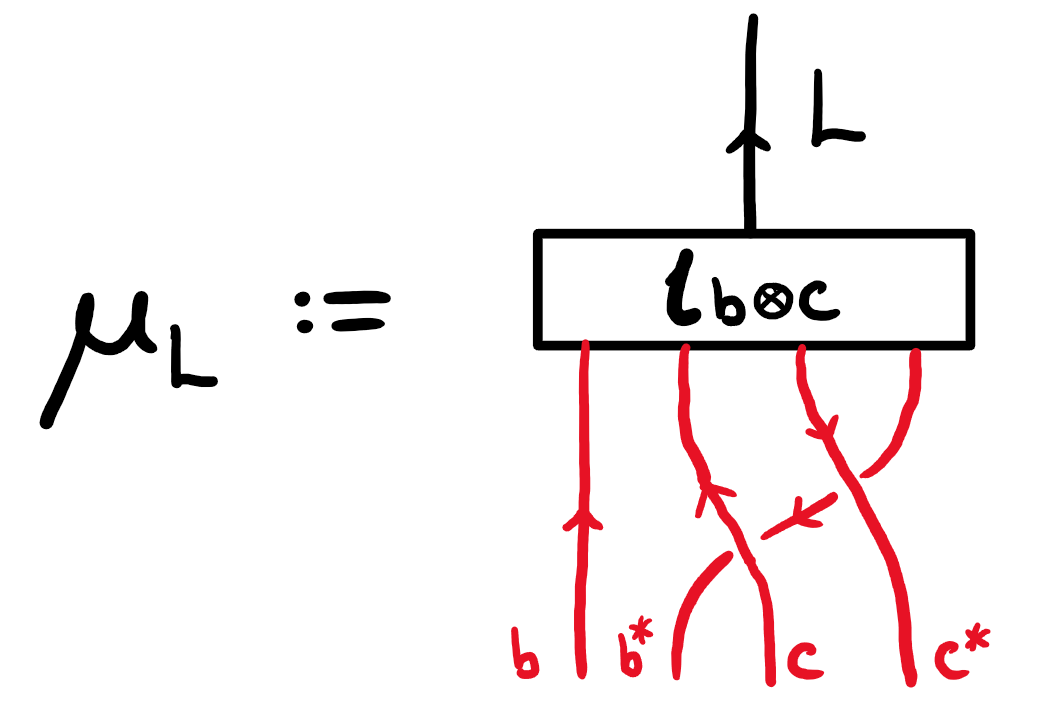}
    };
\end{tikzpicture}
\]
Since this is dinatural in $b \otimes c$, we can rewrite this by pushing the rightmost crossing through the box.
\begin{equation}
\label{eq:AltMult}
\begin{tikzpicture}
    \node[draw] {
    \includegraphics[scale=0.25]{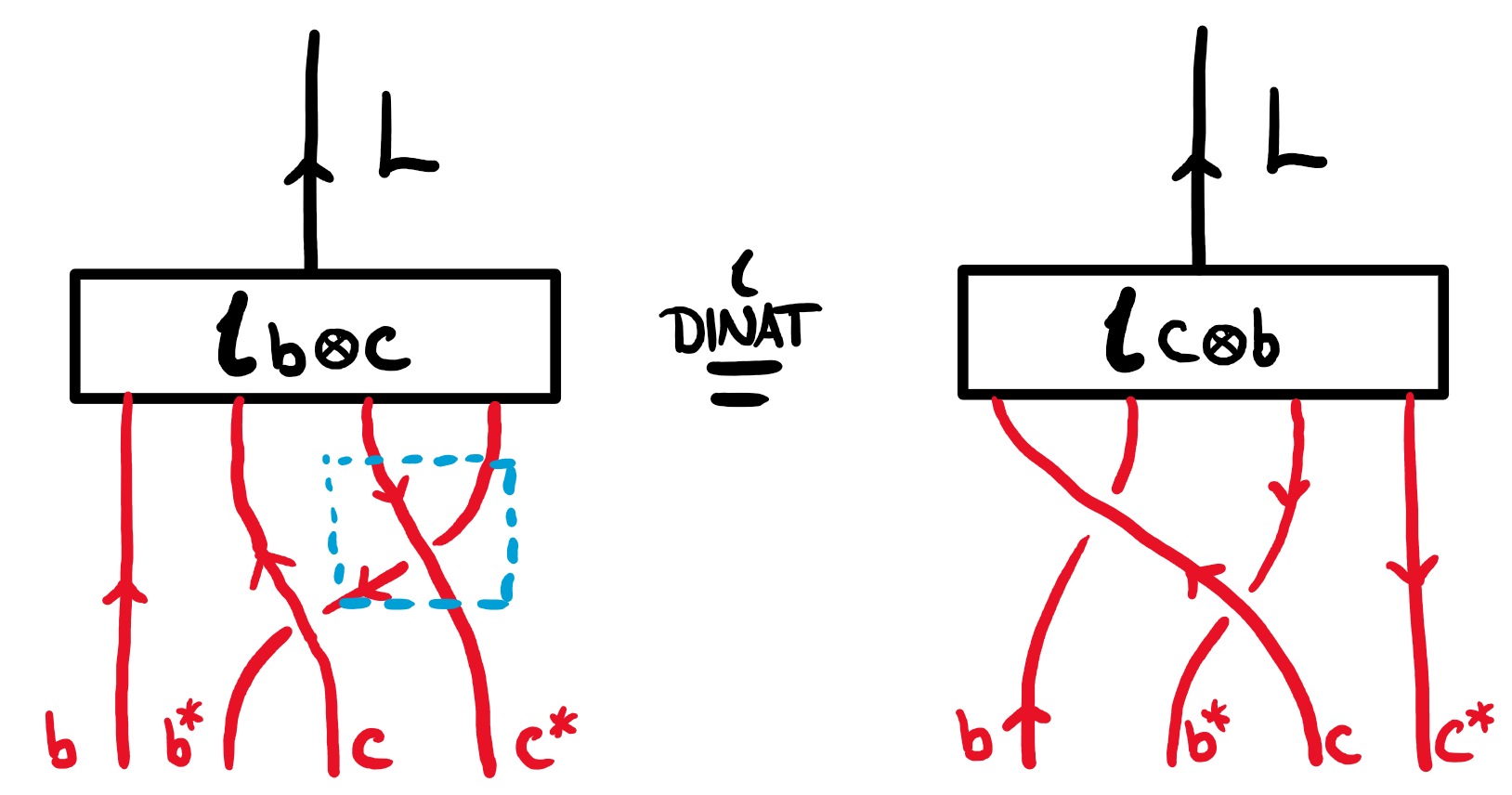}
    };
\end{tikzpicture}
\end{equation}
We'll use this equation later, but for now just let it serve as an example for the graphical calculus of dinatural transformations. Finally, the half-braiding on $L$ is defined via the following dinatural transformation
\[
\begin{tikzpicture}
    \node[draw] {
    \includegraphics[scale=0.25]{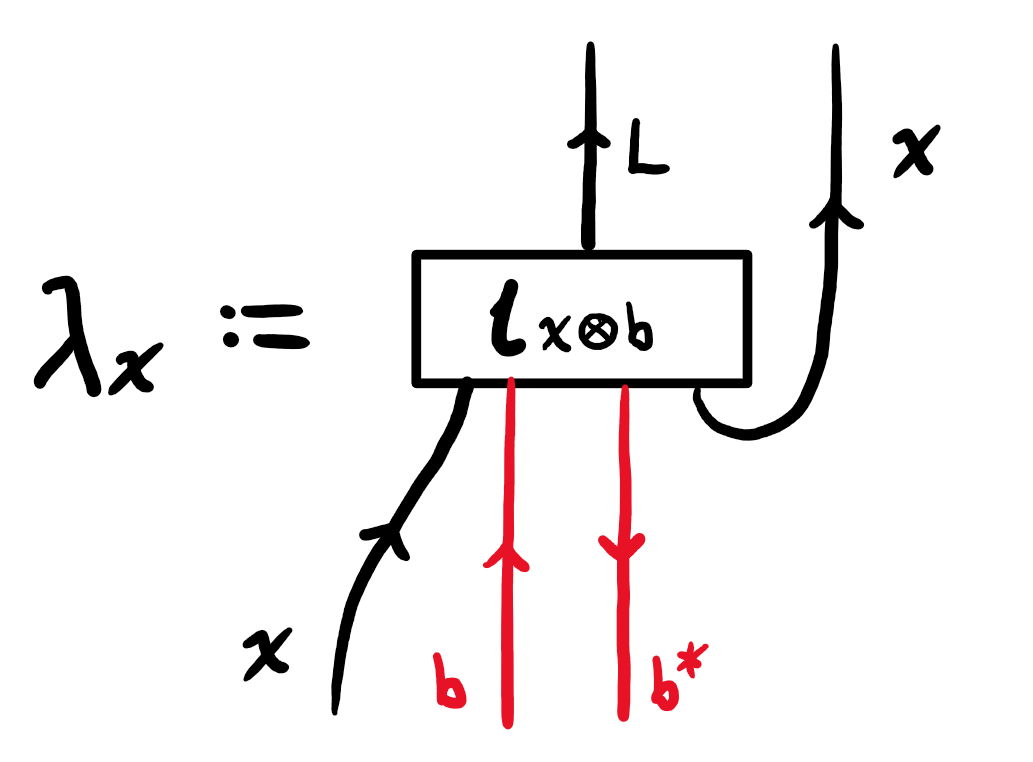}
    };
\end{tikzpicture}
\]
The dinaturality of $\iota$ then implies this is natural in $x$. For instance, consider what happens if we apply the dinaturality move using a morphism $f \otimes \operatorname{id}_b \colon x \otimes b \rightarrow y \otimes b$. Lemma 2.46 of \cite{JFR} shows (indirectly) that $(L,\lambda)$ is indeed separable.

Now that we've characterized $L$ as a half-braided algebra, let's turn to its twist. In Lemma 3.40 of \cite{JFR}, the authors claim that $\theta_{(L,\lambda)}$ is actually given by a dinatural twist. By the universal property of $L$ it suffice to compute $\theta_{(L,\lambda)}\cdot \iota_b$, and so their proof goes as follows.
\[
\begin{tikzpicture}
    \node[draw] {
    \includegraphics[scale=0.5]{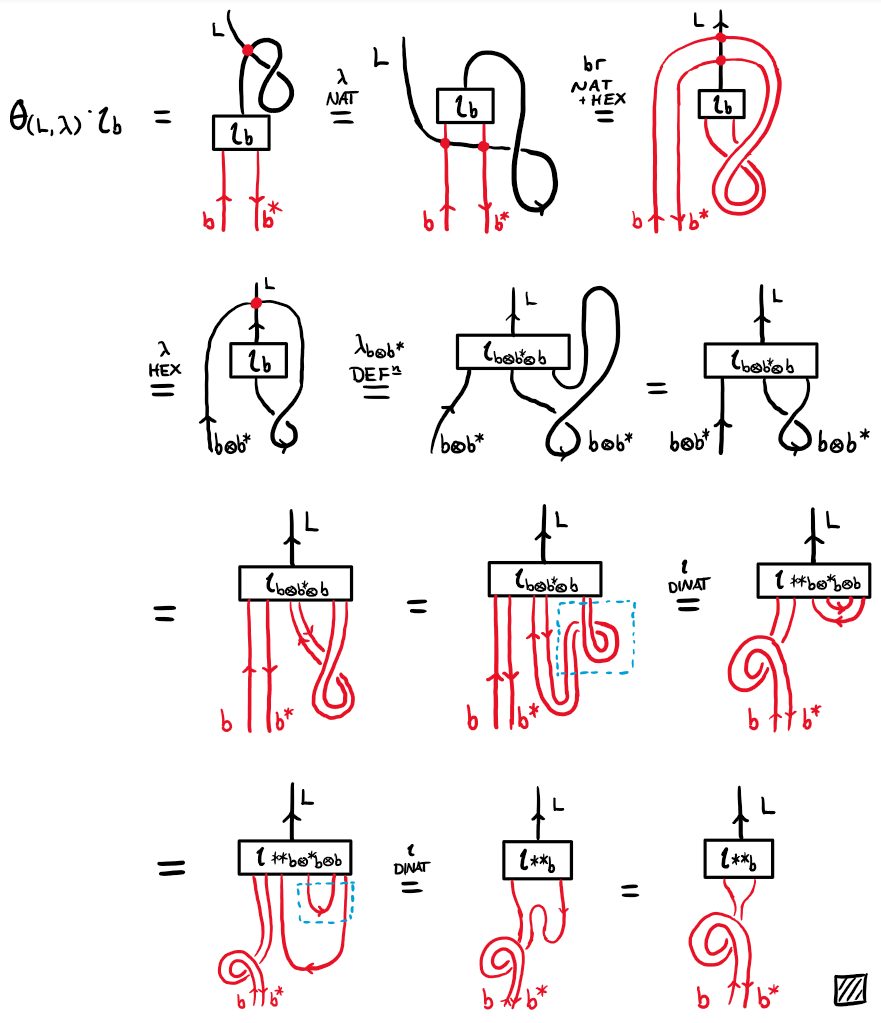}
    };
\end{tikzpicture}
\]
At the second line, we simply write $b \otimes b^*$ as a single object to make some of the moves cleaner. Now, In order to make use of the formula in (\ref{eq:TraceTrace}), we need an isomorphism of half-braided algebras $L \rightarrow \leftindex^{\text{op}}{L}$. This is supplied by the \textit{half-twist} which as a morphism $\psi \colon L \rightarrow L$ of underlying objects is given by the dinatural map
\[
\begin{tikzpicture}
    \node[draw] {
    \includegraphics[scale=0.25]{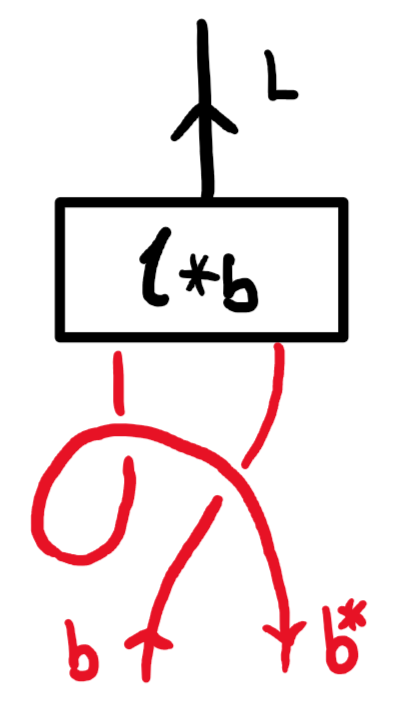}
    };
\end{tikzpicture}
\]
Notice that $\psi$, similarly to $\mu_L$, is not defined directly as a dinatural map $b \otimes b^* \rightarrow L$, but rather as a \textit{natural} map $\psi_b \colon b \otimes b^* \rightarrow \leftindex^*{b}\otimes b$ composed with the canonical dinatural map $\iota$. The main upshot is that it lets us easily compute composition of maps $L \rightarrow L$. In particular, $\psi^2$ via the following dinatural composite.
\[
\begin{tikzcd}
	{b\otimes b^*} & {^*b\otimes b} & {^{**}b\otimes\leftindex^*{b}} & L
	\arrow["{\psi_b}", from=1-1, to=1-2]
	\arrow["{\psi_{\leftindex^*{b}}}", from=1-2, to=1-3]
	\arrow["{\iota_{^{**}b}}", from=1-3, to=1-4]
\end{tikzcd}
\]
We call the map $\psi$ the half-twist because it squares to the twist. Indeed, it's not too hard to show graphically.
\[
\begin{tikzpicture}
    \node[draw] {
    \includegraphics[scale=0.25]{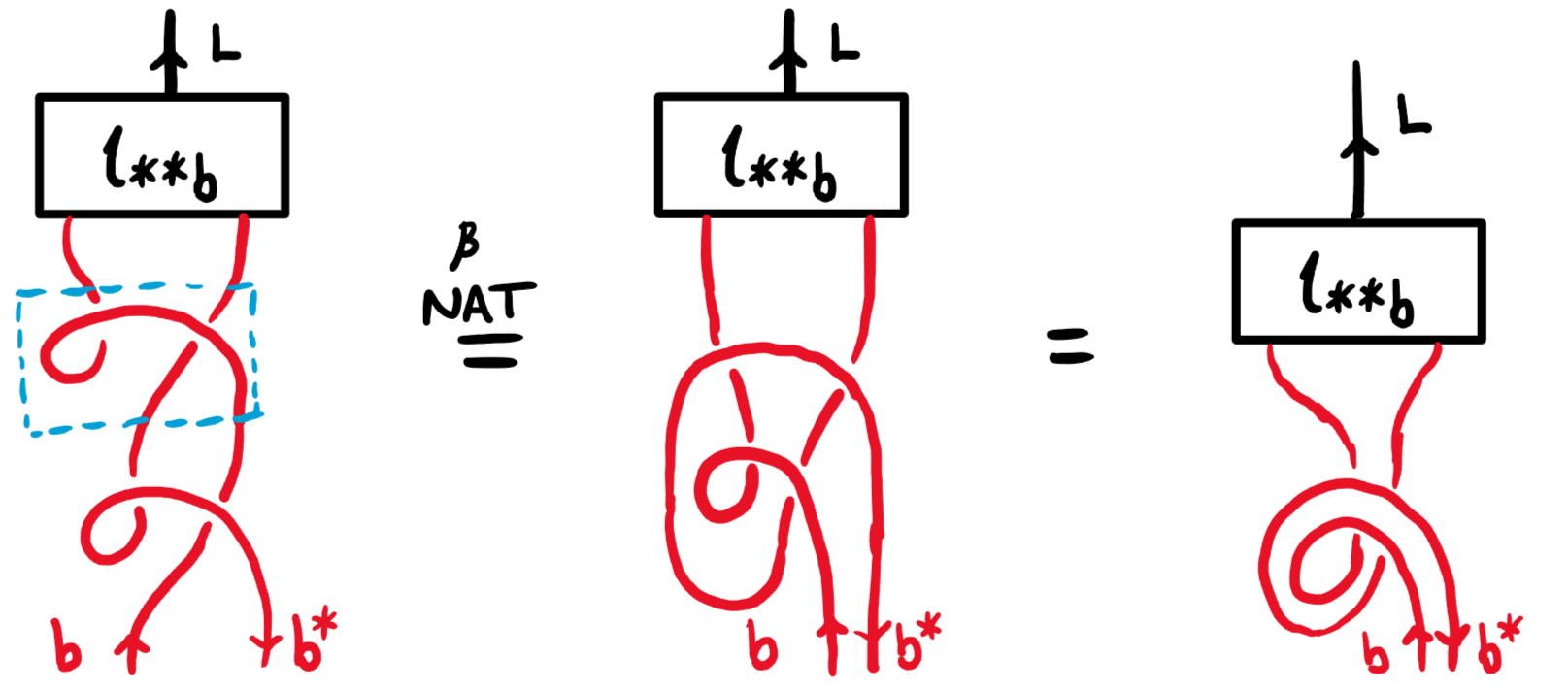}
    };
\end{tikzpicture}
\]
The first moves uses the naturality of the braiding $\beta$ in order to push the second half-twist under the loop of the first. Then we use the hexagon axiom to group strands together. With our self-duality in hand, the last thing we need to work out is the vector space of twisted traces of $(L,\lambda)$. This is done in Lemma 3.42 of \cite{JFR}, where they show that \textit{any} morphism $\varepsilon \colon L \rightarrow \mathds{1}$ is a twisted trace. As usual, define $\varepsilon$ in terms of a dinatural transformation $\varepsilon_b \colon b \otimes b^* \rightarrow \mathds{1}$. The proof is then by direct computation.
\[
\begin{tikzpicture}
    \node[draw] {
    \includegraphics[scale=0.5]{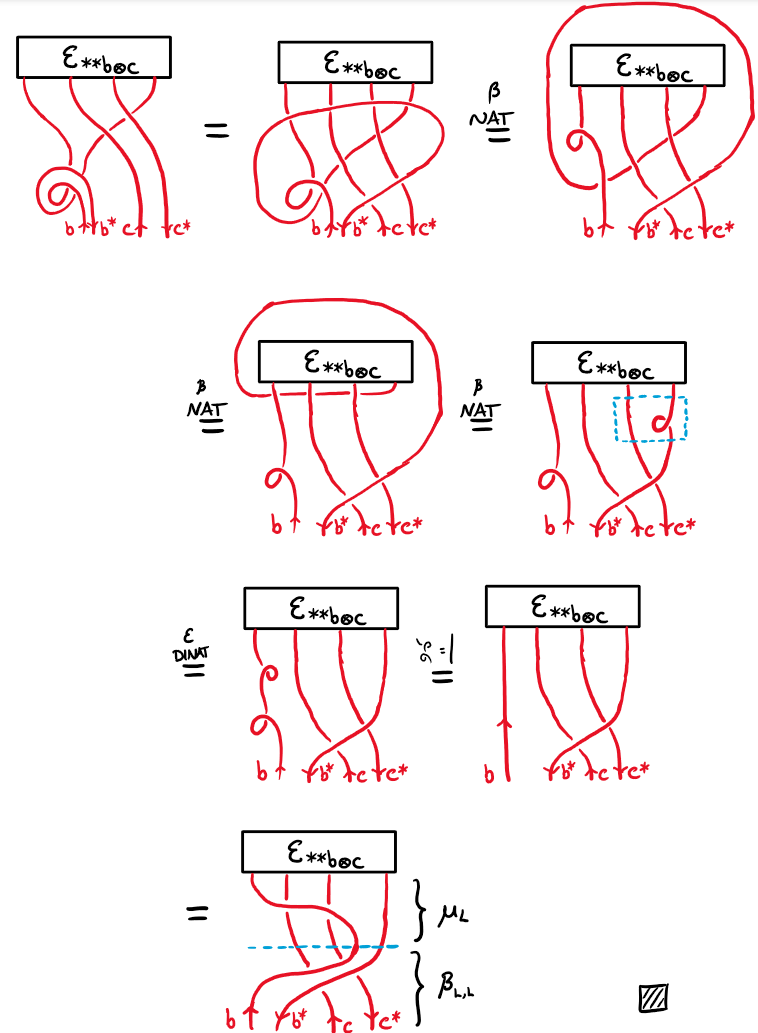}
    };
\end{tikzpicture}
\]
At the last step, we made use of the alternate formula for the multiplication of $L$ given by (\ref{eq:AltMult}). This then gives $\underline{\eta}(L,\lambda) \cong \operatorname{Hom}_{\mathcal{B}}(L,\mathds{1})$ and so the twisted traces of $(L,\lambda)$ are given by dinatural maps $\alpha_b \colon b \otimes b^* \rightarrow \mathds{1}$. This is precisely the same as a natural transformation $\operatorname{id}_{\mathcal{B}} \Rightarrow \operatorname{id}_{\mathcal{B}}$ via the following isomorphism.
\[
\begin{tikzpicture}
    \node[draw] {
    \includegraphics[scale=0.25]{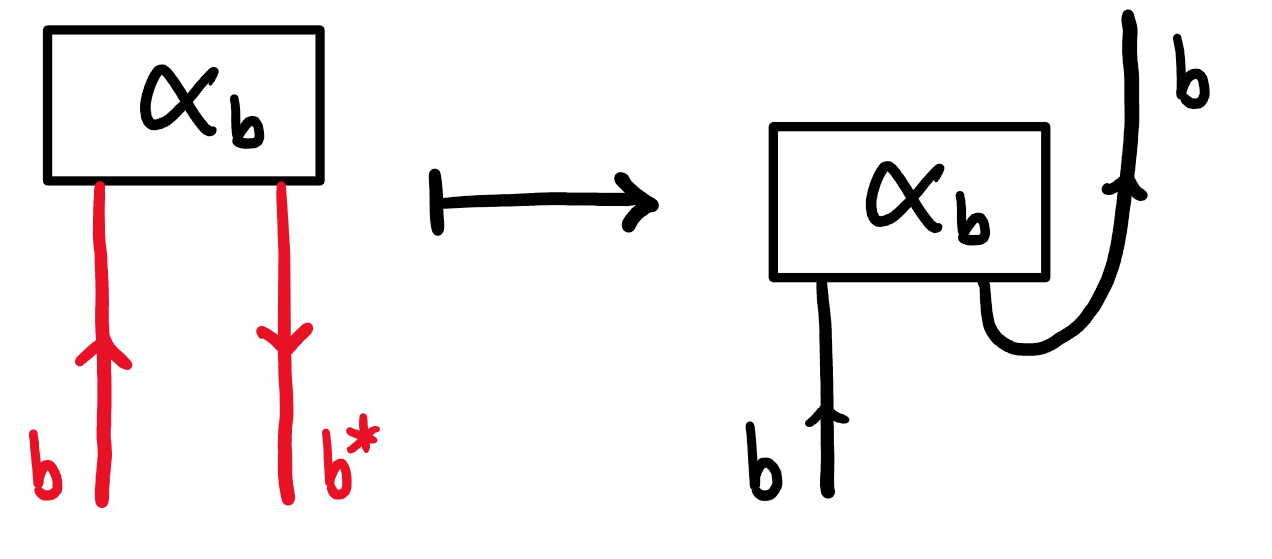}
    };
\end{tikzpicture}
\]
And so we have $\underline{\eta}(L,\lambda) \cong \operatorname{End}(\operatorname{id}_{\mathcal{B}})$, which by the semisimplicity of $\mathcal{B}$, is the same as the linear dual of the (linearized) Grothendieck ring of $\mathcal{B}$, i.e. 
\[
\underline{\eta}(L,\lambda)\cong (K_0(\mathcal{B})\otimes_{\mathbb{Z}} \mathds{k})^* \coloneqq \operatorname{Hom}_{\operatorname{Vec}}(K_0(\mathcal{B})\otimes_{\mathbb{Z}} \mathds{k},\mathds{k})
\]
We already have a formula for the map $\underline{\eta}(\psi)\cdot\underline{\operatorname{fl}}_L$, so we can translate it into an automorphism of $\operatorname{End}(\operatorname{id}_{\mathcal{B}})$ by chasing a natural transformation along the following diagram.
\[
\begin{tikzcd}
	{\underline{\eta}(L,\lambda)} && {\underline{\eta}(L,\lambda)} \\
	{\operatorname{End(id}_{\mathcal{B}})} && {\operatorname{End(id}_{\mathcal{B}})}
	\arrow["{\underline{\eta}(\psi)\cdot\underline{\operatorname{fl}}_L}", from=1-1, to=1-3]
	\arrow["{{\resizebox{3pt}{12pt}{$\wr$}}}", from=1-3, to=2-3]
	\arrow["\resizebox{3pt}{12pt}{$\wr$}", from=2-1, to=1-1]
	\arrow[from=2-1, to=2-3]
\end{tikzcd}
\]
Graphically, this amounts to the following map.
\[
\begin{tikzpicture}
    \node[draw] {
    \includegraphics[scale=0.3]{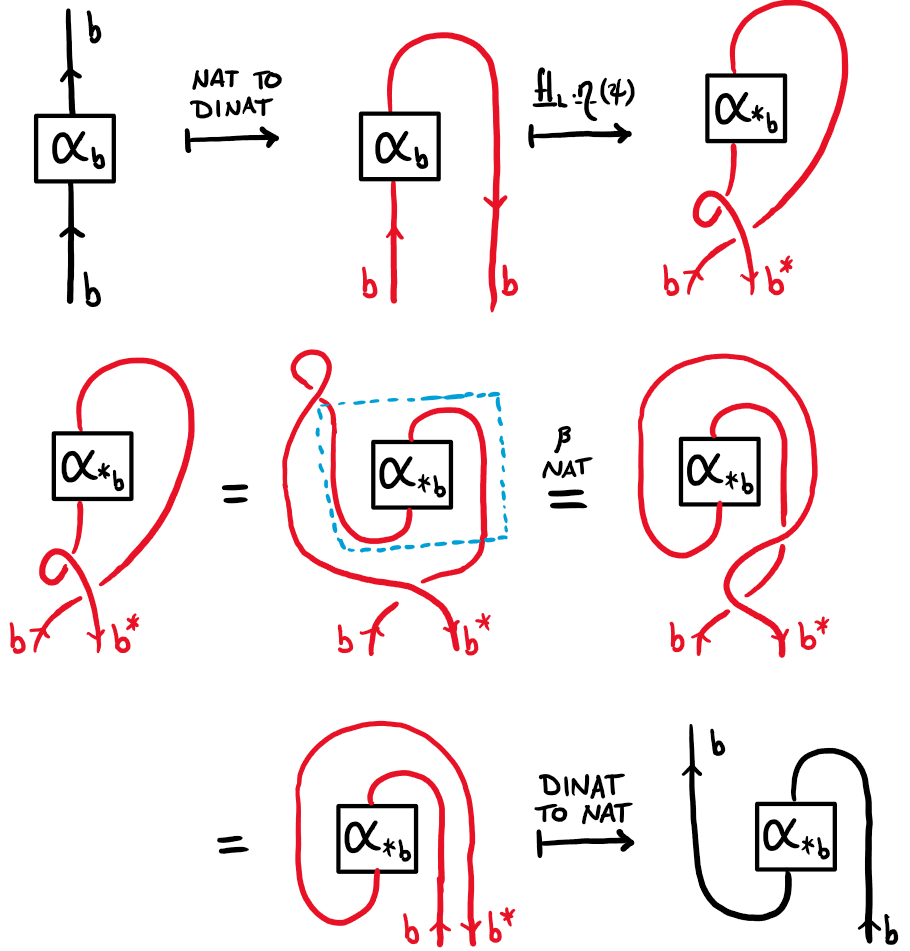}
    };
\end{tikzpicture}
\]
The map $\operatorname{End}(\operatorname{id}_{\mathcal{B}}) \rightarrow (K_0(\mathcal{B})\otimes_{\mathbb{Z}}\mathds{k})^*$ is given by $\alpha \mapsto (b \mapsto \langle \alpha_b \rangle)$ where $\alpha_b = \langle \alpha_b \rangle \operatorname{id}_b$. Thus, since $\langle \alpha_{\leftindex^*{b}}\rangle = \langle (\alpha_{\leftindex^*{b}})^* \rangle$ the automorphism $\underline{\eta}(\psi)\cdot\underline{\operatorname{fl}}_L$ translates into the linear dual of the map
\begin{align}
\label{eq:DualMap}
\begin{split}
        K_0(\mathcal{B})\otimes_{\mathbb{Z}} \mathds{k} &\longrightarrow K_0(\mathcal{B})\otimes_{\mathbb{Z}} \mathds{k}\\
        b &\longmapsto \leftindex^*{b}
\end{split}
\end{align}

We're finally at the home stretch. Let $\mathcal{B}$ be a slightly-degenerate braided fusion category, and recall that, the braided fusion 2-category $\mathcal{Z}_1(\Sigma\mathcal{B}) \simeq \operatorname{sHBA}(\mathcal{B})$ has two components. The identity component with two simples being the monoidal unit and the ``electron" $e$, and the ``magnetic" component, also with two simples. This electron then has a trivial full-braid with any object purely in the identity component, while for a purely magnetic object, the full braid with $e$ picks up a negative sign. We can decompose any half-braided algebra as its identity part plus magnetic part, which we denote by $(A,\gamma) = (A^+,\gamma^+) \oplus (A^-,\gamma^-)$. Moreover, the braiding with the electron gives us projectors $P^{\pm}_{A} \colon (A,\gamma) \rightarrow (A^{\pm},\gamma^{\pm})$ which are graphically described as follows.
\[
\begin{tikzpicture}
    \node[draw] {
    \includegraphics[scale=0.2]{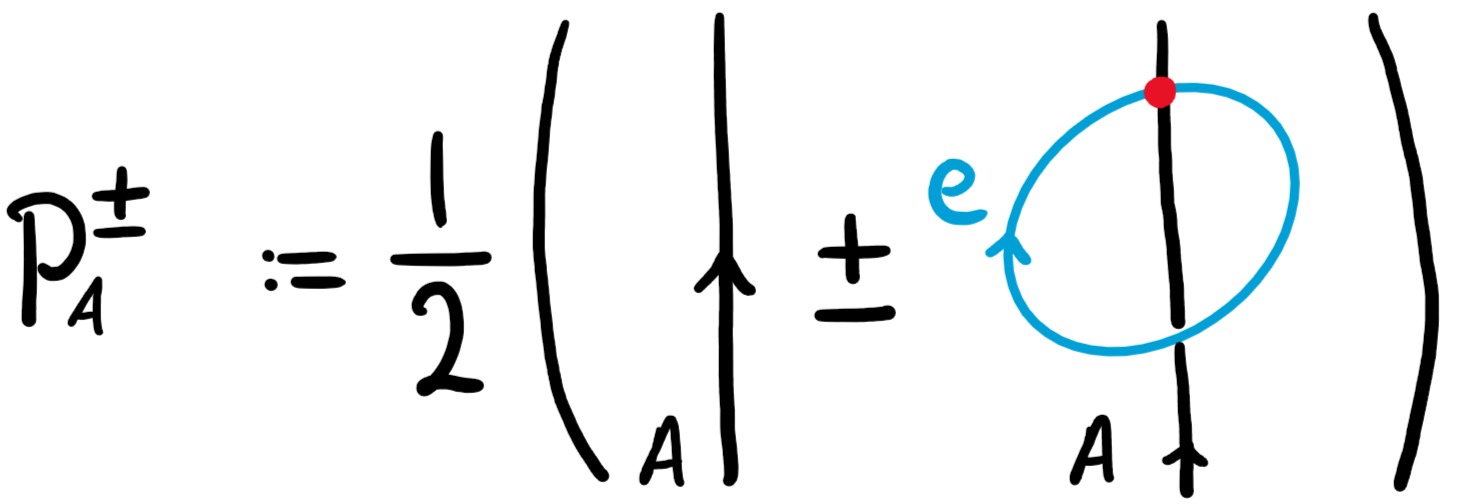}
    };
\end{tikzpicture}
\]
These projectors also let us decompose the vector space of twisted traces as $\underline{\eta}(A,\gamma) = \underline{\eta}(A^+,\gamma^+) \oplus \underline{\eta}(A^-,\gamma^-)$. Using a similar argument to the one used to obtain the map (\ref{eq:DualMap}), Lemma 3.45 of \cite{JFR} gives that the map $\underline{\eta}(P^{\pm}_{L})$ corresponds to the linear dual, i.e. transpose, of the map
\begin{align}
\label{eq:LDualProjector}
\begin{split}
        \underline{\eta}(P^{\pm}_{L})^T \colon K_0(\mathcal{B})\otimes_{\mathbb{Z}} \mathds{k} &\longrightarrow K_0(\mathcal{B})\otimes_{\mathbb{Z}} \mathds{k}\\
        b &\longmapsto \frac{1}{2}(1 \pm e) b
\end{split}
\end{align}
Now we're finally ready to compute $\kappa(L,\lambda)$. We do so by decomposing $L$ via the projectors $P^\pm_{L}$ and computing it on each component separately.
\begin{align*}
    \begin{split}
        \kappa(L^{\pm}, \psi^{\pm}) &= \operatorname{trace}\Big(\underline{\eta}(P^{\pm}_{L})\cdot\big( \underline{\eta}(\psi^{\pm})\cdot\underline{\operatorname{fl}}_L\big)\Big)\\
        &= \operatorname{trace}\Big(\big(\underline{\eta}(\psi^{\pm})\cdot\underline{\operatorname{fl}}_L\big)\cdot\underline{\eta}(P^{\pm}_{L})\Big)\\
        &=\operatorname{trace}\Big((b \mapsto \leftindex^*{b})^T\cdot\left(b \mapsto 1/2(1\pm e)b\right)^T\Big)\\
        &= \operatorname{trace}\Big(\left(b\mapsto 1/2(1\pm e)\leftindex^*{b} \right)^T\Big)\\
        &= \operatorname{trace}\Big(b\mapsto 1/2(1\pm e)\leftindex^*{b} \Big)\\
        &= \frac{1}{2} \operatorname{trace}\big(b\mapsto \leftindex^*{b}\big) \pm \frac{1}{2} \operatorname{trace}\big(b\mapsto e\leftindex^*{b}\big)\\
        &= \frac{1}{2}\Big|\{b \in \pi_0\mathcal{B}\ |\ b \cong \leftindex^*{b}\}\Big| \pm \frac{1}{2}\Big|\{b \in \pi_0\mathcal{B}\ |\ b \cong e\otimes\leftindex^*{b}\}\Big|\\
        &\geq \frac{1}{2}(1 \pm K)
    \end{split}
\end{align*}
for some $K \in \mathbb{N}$. We used many elementary properties of the trace of linear maps, as well as the fact that there is at least one simple object of $\mathcal{B}$ that is isomorphic to its left dual, namely $\mathds{1}$. The last piece of the puzzle is to show that $K = 0$, and we do so via some properties of the 1-categorical version of $\eta$. For one, we have that $\eta(x) = \eta(\leftindex^*{x})$ as $\operatorname{fl}_x$ decategorifies into a mere equality. Moreover, since $e$ is transparent, then $\eta(e \otimes x) = \eta(e)\eta(x)$ which can be graphically seen as follows.
\[
\begin{tikzpicture}
    \node[draw] {
    \includegraphics[scale=0.25]{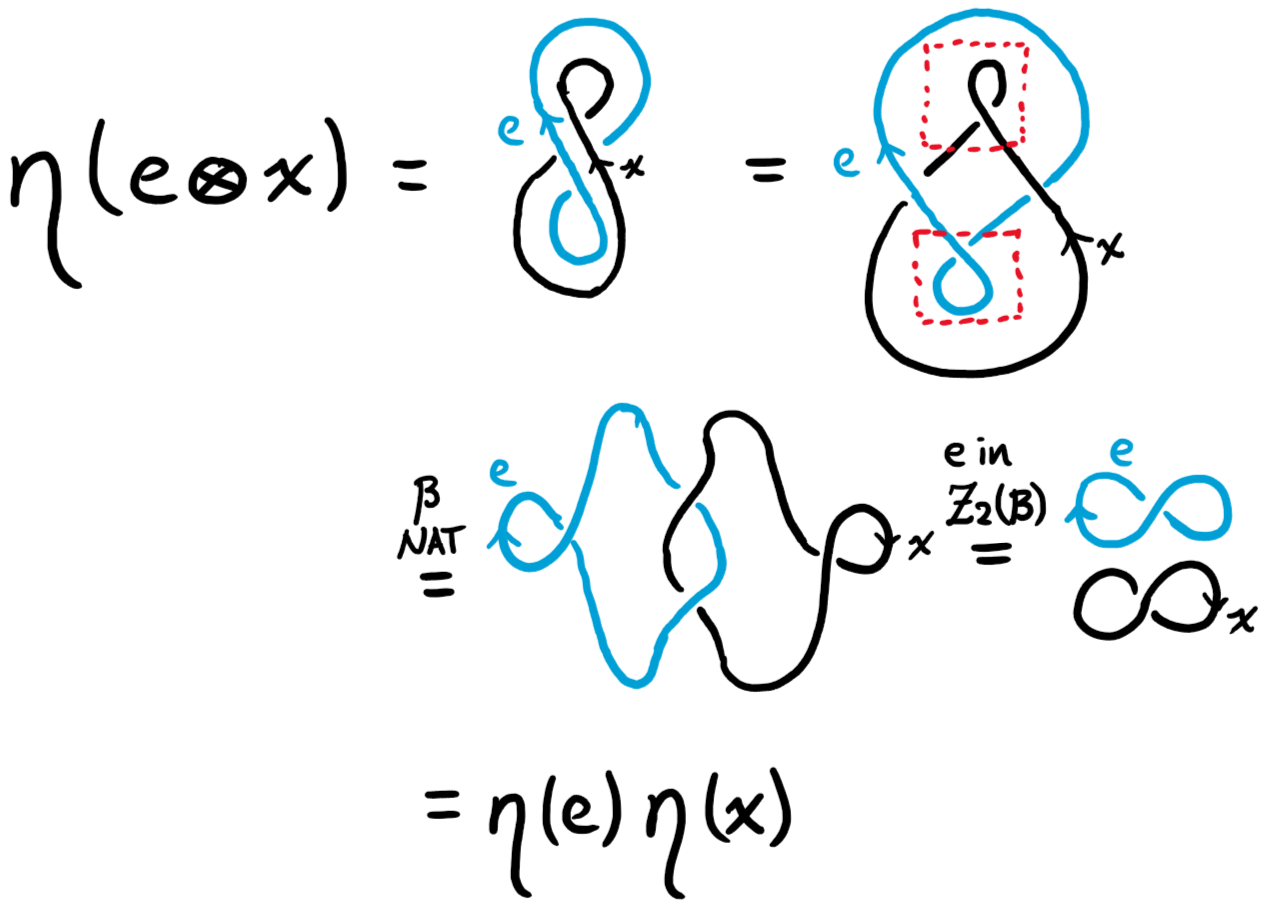}
    };
\end{tikzpicture}
\]
Since $e$ is the odd line in $\mathcal{Z}_2(\mathcal{B}) \simeq \operatorname{sVec}$, we have $\beta_{e,e} = -\operatorname{id}_{e \otimes e}$ and so $\eta(e) = - \operatorname{dim}_{\mathds{k}}(e) = -1$. Now suppose that $K > 0$, i.e. that there is a $b \in \pi_0 \mathcal{B}$ such that $b \cong e \otimes \leftindex^*{b}$. This implies that
\begin{align*}
    \eta(b) &= \eta(e\otimes \leftindex^*{b})\\
    &= \eta(e) \eta(\leftindex^*{b})\\
    &= -\eta(b)
\end{align*}
Thus we must have $\eta(b)=0$. However for any simple object $b$, we have $\eta(b) \neq  0$, and therefore by contradiction, $K = 0$. In particular, for the purely magnetic $(L^-,\lambda^-)$, we have $\kappa(L^-,\lambda^-) > 0$. By what we learned in the previous talks, we then have $\operatorname{sHBA}(\mathcal{B}) \simeq \mathcal{Z}_1(\Sigma\mathcal{B}) \simeq \mathcal{S}$, which means that $\mathcal{B}$ admits a minimal nondegenerate extension.

%% file: talks/5.2/5.2.tex




Talk by Matthew Yu, notes by Iordanis Romaidis.

\section{Perspective on fusion 2-categories}\label{sec:perspective}

A fusion 2-category (F2C) is called \textit{connected} if it only has a single component of the identity. 

\begin{example}
	Consider the F2C $2\!\operatorname{Rep}(G) = \Sigma \Rep (G)$ for a finite group $G$ \cite{greenough2010monoidal}. The objects are pairs $([H], \mu)$ where $ H$ is a subgroup of $G$, $[H]$ is an isomorphism class of $G/H$
	and $\mu \in H^2(G,\mathbb{C}^\times)$. The fusion rules are described by
	\begin{equation}
		([H], \mu) \cdot ([K], \sigma) = \sum_{[a]\in H\backslash G/K}{([H\cap a K a^{-1}], \mu \sigma^a)}~
	\end{equation}
	where $\sigma^a(x,y) := \sigma(axa^{-1},a y a^{-1})$. 
	
	We unpack this example for $G =\mathbb{Z}_2$: There are two objects $X^\mathrm{id}, X^{\mathbb{Z}_2}$ corresponding to the two subgroups as the corresponding group cohomology is trivial. The only non-trivial fusion rule is \[ (X^{\mathbb{Z}_2})^2 = X^{\mathbb{Z}_2} \oplus X^{\mathbb{Z}_2}~.\] 
	On the level of 1-morphisms, there are two simple endomorphisms of $X^{\mathrm{id}}$ corresponding to the two representations of $\mathbb{Z}_2$, two more endomorphisms of $X^{\mathbb{Z}_2}$ and 1-morphisms $X^{\mathrm{id}}\rightleftarrows X^{\mathbb{Z}_2}$. This is summarized in the following diagram: 
	\[
	\begin{tikzcd}
	X^{\mathrm{id}}
	 \arrow[out=150,in=210,loop, swap, "\Rep(\mathbb{Z}_2)" ] \arrow[r,->,shift left=1, "\Vect"] \arrow[r,<-,shift right=1, "\Vect" swap] & X^{\mathbb{Z}_2}\arrow[out=330, in = 30, loop, "\Vect_{\mathbb{Z}_2}" swap]
	\end{tikzcd}
	\]
	\[\Mod(\Vect_{A[0]}) = 2\!\Vect_{A[1]}\] 
	We can also take: 
	\begin{align}
		\Sigma \operatorname{sVect} &= 2\!\Vect_{\mathbb{Z}_2[1]}^\omega \\
		\Sigma \Vect_{\mathbb{Z}_2} &\rightarrow 2\!\Vect_{\mathbb{Z}_2[1]}
	\end{align}
	where $\omega$ arises from the braiding of $\operatorname{sVect}$.
\end{example}

A F2C $\mathfrak{C}$ with $\Omega \mathfrak{C} = \Vect$ or $\Omega\mathfrak{C} = \operatorname{sVect}$ is called \textit{strongly fusion}.

\begin{example}
	\begin{enumerate}
	\item	$2\!\Vect_{G}^\pi$  and  $2\!\Vect_{(G,z)}^\varpi = 2\!\operatorname{sVect}_{G}^\varpi$ where $\pi \in H^4(\mathrm{B}G,\mathbb{C}^\times)$ and $ \varpi\in \widetilde{SH}^4(G,z)$ is an element in the reduced supercohomology.
	\item Recall that the Tambara-Yamagami F1C associated to the group $\mathbb{Z}_2$ is the $\mathbb{Z}_2$-graded fusion category $(\Vect_{\mathbb{Z}_2})_0 \oplus (\Vect)_1$ where the simple $D$ in the non-trivial graded component $\Vect$ satisfies the fusion rules $a D = D a = D$ and $D^2 = 1 \oplus a$.
	
	The Tambara-Yamagami F2C is defined as
	\begin{equation}
		2\!\operatorname{TY}(\mathbb{Z}_2) := \operatorname{Bimod}_{2\!\Vect_{D^8}}(\Vect_{\mathbb{Z}_2})
	\end{equation}
	and represented by the following diagram:
	\[
	\begin{tikzcd}
	I
	 \arrow[out=150,in=210,loop, swap, "\Vect_{\mathbb{Z}_2}" ] \arrow[r,->,shift left=1] \arrow[r,<-,shift right=1] & X\arrow[out=330, in = 30, loop, "\Vect_{\mathbb{Z}_2}" swap]\\
	 J
	 \arrow[out=150,in=210,loop, swap, "\Vect_{\mathbb{Z}_2}" ] \arrow[r,->,shift left=1] \arrow[r,<-,shift right=1] & Y\arrow[out=330, in = 30, loop, "\Vect_{\mathbb{Z}_2}" swap] \\ 
	 D
	 \arrow[out=150,in=210,loop, swap, "\Vect" ],& D^2  = X\boxplus Y~.
	\end{tikzcd}
	\]
    This example was considered in \cite[Example 5.2.1]{D_coppet_2025}, see also the references therein for its appearance in the physics literature.
	\end{enumerate}
\end{example}

\begin{theorem}
	Let $\mathfrak{C}$ be a bosonic F2C with $Z_2(\Omega \mathfrak{C}) \simeq \Vect$. Then it is Morita equivalent to $2\!\Vect_G^\pi \boxtimes \Sigma \mathcal{B}$ where $\mathcal{B}$ is a braided non-degenerate fusion 1-category and $\pi\in H^4(\mathrm{B}G, \mathbb{C}^\times)$.
\end{theorem}

The statement is no longer true in the fermionic case where $Z_2(\Omega \mathfrak{C}) = \operatorname{sVect}$ and requires the choice of an MNE. 

\begin{theorem}
	Let $\mathfrak{C}$ be a fermionic F2C with $Z_2(\Omega \mathfrak{C}) = \operatorname{sVect}$.
	Then, $\mathfrak{C}$ is not canonically Morita to $2\!\operatorname{sVect}_G^\omega \boxtimes\Sigma\mathcal{B}$ for a slightly non-degenerate $\mathcal{B}$ and it depends on a choice of MNE. 
\end{theorem}

\section{Parametrizations of F2Cs}

\begin{figure}[h!]
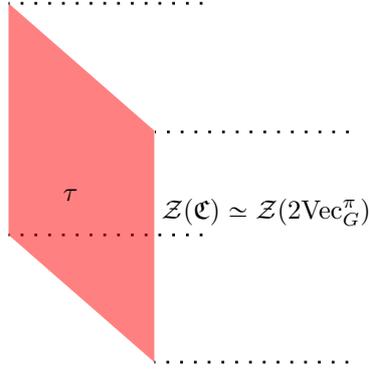

    \centering
    \begin{equation*}
    	\boxpic{0.9}{bulk-boudnary.pdf}{
    		\put (42,40) {$\mathcal{Z}(\mathfrak{C}) \simeq \mathcal{Z}(2\!\Vect_G^\pi)$}
    		\put (15,45) {$\large\tau$}
    	}
    \end{equation*}
    \caption{Boundary condition for the DW theory associated to $2\!\Vect_G^\pi$~.}
    \label{fig:bulk-bound}
\end{figure}

Before we arrive to the full classification, we present the physical perspective to see the necessary classifying data. The idea is to start with $\mathcal{Z}(\mathcal{C})$ in a bulk $(3+1)$D TQFT and $\mathfrak{C}$ as the symmetry of a $(2+1)$D theory on the boundary, see Figure~\ref{fig:bulk-bound}.  

One can look for boundary conditions of $\mathcal{Z}(2\!\Vect_G^\pi)$. There are minimal boundary conditions which are Dirichlet for $G$, while other come from gauging subgroups $H\subset G$.
In this dimension, there is also the choice to take a Dirichlet boundary condition stacked on by a $(2+1)$-dim TQFT associated to a non-degenerate BF1C $\mathcal{A}$ with $H$-symmetry, \ie $\rho: H \rightarrow \Aut_{br}(\mathcal{A})$ and gauging the diagonal $H$-action, while anomaly matching $O(\rho)$ and $\pi|_H$.

\begin{theorem}[\cite{class2fusion}]
Bosonic F2Cs are parametrized by the following data: 
\begin{itemize}
	\item A non-degenerate BF1C $\mathcal{A}$, 
	\item An inclusion of finite groups $H\hookrightarrow G$,
	\item A monoidal functor $\rho: H \rightarrow \Aut_{br}(\mathcal{A})$,
	\item A class $\pi \in H^4(\mathrm{B}G, \mathbb{C}^\times)$,
	\item and a homotopy between the anomaly of $\rho$ and $\pi|_H$. 
\end{itemize}
The anomaly of $\rho$ refers to the obstruction $O(\rho)$ of lifting $\rho$ to $\mathrm{Pic}(\mathcal{A})$.
\end{theorem}

Hence, a bosonic F2C is the data described by a commuting square: 

\begin{equation}
	\begin{tikzcd}
		\mathrm{B}H \arrow[d, "i",swap]
		  \arrow[r,->,"\rho"] & \mathrm{B}\Aut^\mathrm{br}(\mathcal{A}) \arrow[d, "\brack{-}"]\\
		 \mathrm{B}G\arrow[ur , Rightarrow, "\simeq"]\arrow[r,->,"\pi"]  & \mathrm{B}^4 \mathbb{C}^\times.
		\end{tikzcd}
	\end{equation}	

\begin{example}
	\begin{enumerate}
		\item If $H = G$, then the corresponding F2C is $\Sigma \mathcal{A}^G = \Mod(\mathcal{A}^G)$.
		\item If $H = 1$, then we have $ 2\!\Vect_G^\pi \boxtimes \Sigma\mathcal{A}$. 
	\end{enumerate}
\end{example}

For the fermionic case, recall that a finite supergroup $(G,z)$ is defined by:
\[
	1\rightarrow z \rightarrow (G,z) \rightarrow G \rightarrow 1\]
where $z$ is of order 2 and the extension determined by a class $\kappa \in H^2(G, \mathbb{Z}_2)$. Let $\Aut^\mathrm{br}_{\operatorname{sVect}}(\mathcal{A})$ denote the $\operatorname{sVect}$-linear braided automorphisms of $\mathcal{A}$. 

\begin{theorem}[\cite{class2fusion}]
Fermionic F2Cs are parametrized by: 
\begin{itemize}
	\item A slightly non-degenerate BF1C $\mathcal{A}$,
	\item An inclusion of finite supergroups $i:(H,z)\rightarrow (G,z)$,
	\item A monoidal functor of higher supergroups $\rho: (H,z) \rightarrow \Aut^{\mathrm{br}}_{\operatorname{sVect}}(\mathcal{A})$,
	\item A class in a torsor of the supercohomology $\varpi \in \widetilde{SH}^{4}(\mathrm{B}(G,z))$ 
	\item A homotopy between the anomaly of $\rho$ and $\varpi|_{(H,z)}$.
\end{itemize}
\end{theorem}

For the homotopy classification we build a functor from the space of fusion 2-cats into the space of commuting squares.